\documentstyle{amsppt}
\loadbold
\magnification=\magstep1
\overfullrule=0pt 
	\def\card{\mathop{\text{\rm card}}\nolimits} 
	\def\diam{\mathop{\text{\rm diam}}\nolimits} 
	\def\dist{\mathop{\text{\rm dist}}\nolimits} 
	\def\Im{\mathop{\text{\rm Im}}\nolimits} 
	\def\Int{\mathop{\text{\rm Int}}\nolimits} 
	\def\Lip{\mathop{\text{\rm Lip}}\nolimits} 
	\def\olim{\mathop{\overline{\text{\rm lim}}}\nolimits}
	\def\ulim{\mathop{\underline{\text{\rm lim}}}\nolimits}
	\def\oosc{\mathop{\overline{\text{\rm osc}}}\nolimits}
	\def\uosc{\mathop{\underline{\text{\rm osc}}}\nolimits}
	\def\os{\mathop{\text{\rm os}}\nolimits} 
	\def\osc{\mathop{\text{\rm osc}}\nolimits} 
	\def\Re{\mathop{\text{\rm Re}}\nolimits} 
	 
	\def\supp{\mathop{\text{\rm supp}}\nolimits} 
	\def\wt{\mathop{\text{\rm wt}}\nolimits} 
	\def\varep{\varepsilon}
	\def\complex{{\Bbb C}}
	
	\def\real{{\Bbb R}}
	\def\B{{\Cal B}}
	\def\D{{\Cal D}}

	\def\P{{\Cal P}}
	\def\U{{\Cal U}}
	\def\V{{\Cal V}}
	\def\W{{\Cal W}}
	
	\def\iitem{\itemitem}
	\def\dfeq{\mathop{\ \buildrel {\text{\rm df}}\over =\ }}
	\def\chix{{\raise.5ex\hbox{$\chi$}}}
	\def\To{\Rightarrow}
	\def\btheta{{\boldsymbol\theta}}
	\def\util{\tilde u}
	\def\uttil{\skew2\tilde\util}
	\def\vtil{\tilde v}
	\def\vttil{\skew2\tilde\vtil}
	\def\bone{\text{\bf 1}}

\topmatter
\title Differences of bounded semi-continuous functions, I\endtitle
\author Haskell Rosenthal\endauthor
\affil The University of Texas at Austin\endaffil
\address Department of Mathematics, The University of Texas at Austin, 
Austin, TX 78712-1082\endaddress 
\thanks This research was partially supported by NSF DMS-8903197 and 
TARP 235.\endthanks 
\abstract Structural properties are given for $D(K)$, the Banach algebra of 
(complex) differences of bounded semi-continuous functons on a metric 
space $K$. For example, it is proved that if all finite derived sets of $K$ 
are non-empty, then a complex function $\varphi$ operates on $D(K)$ 
(i.e., $\varphi\circ f\in D(K)$ for all $f\in D(K)$) if and only if 
$\varphi$ is locally Lipschitz. Another example: if $W\subset K$ and 
$g\in D(W)$ is real-valued, then it is proved that $g$ extends to a 
$\tilde g$ in $D(K)$ with $\|\tilde g\|_{D(K)} = \|g\|_{D(W)}$. Considerable 
attention is devoted to $SD(K)$, the closure in $D(K)$ of the set of simple 
functions in $D(K)$. Thus it is proved that every member of $SD(K)$ is a 
(complex) difference of semi-continuous functions in $SD(K)$, and that 
$|f|$ belongs to $SD(K)$ if $f$ does. An intrinsic characterization of 
$SD(K)$ is given, in terms of transfinite oscillation sets. Using the 
transfinite oscillations, alternate proofs are given of the results of 
Chaatit, Mascioni and Rosenthal that functions of finite Baire-index belong 
to $SD(K)$, and that $SD(K)\ne D(K)$ for interesting $K$. 
It is proved that the ``variable oscillation criterion'' characterizes 
functions belonging to $B_{1/4}(K)$, thus answering an open problem raised 
in earlier work of Haydon, Odell and Rosenthal. It is also proved that $f$ 
belongs to $B_{1/4}(K)$ (if and) only if $f$ is a uniform limit of simple 
$D$-functions of uniformly bounded $D$-norm iff $\osc_\omega f$ is bounded; 
the last equivalence has also been obtained by V.~Farmaki, using other 
methods. Elementary computations of the $D$-norm of some special simple 
functions are given; for example the $D$-norm of $\chix_A$ for a given 
set $A$ is computed precisely, in terms of $\partial^j A$, the $j$-th  
boundary of $A$, $j=1,2,\ldots$. The main structural results on $SD(K)$ 
and $B_{1/4}(K)$ are obtained using the finite oscillations of a given 
function. The higher order oscillations are exploited for the study of 
the transfinite analogues of $B_{1/4}(K)$, in subsequent work. 
\endabstract 
\endtopmatter 

\document

\head Table of Contents\endhead 

{\hsize=5.8truein
\line {\indent \hphantom{1.}  Introduction\dotfill 2}
\line {\indent 1. Preliminaries \dotfill 10}
\line {\indent 2. The $D$-norm of the characteristic function of 
a set\dotfill 29}
\line {\indent 3. The transfinite oscillations; properties and first 
applications\dotfill 38}
\line {\indent 4. Strong $D$-functions and a characterization of 
$B_{1/4}$\dotfill 57}
\line {\indent \hphantom{1.} References\dotfill 79}
\smallskip}

\baselineskip=18pt 
\head Introduction.\endhead 

Let $K$ be a fixed metric space. A function $f:K\to \complex$ is called a 
(complex) difference of bounded semi-continuous functions  if there exist 
bounded lower semi-continuous functions $\varphi_1,\ldots,\varphi_4$ on $K$ 
with $f= (\varphi_1-\varphi_2) + i(\varphi_3-\varphi_4)$. 
We denote the set of all such functions by $D(K)$; 
we also refer to functions in $D(K)$ as $D$-functions. 
A classical result of Baire (proved in Section~1 for completeness) yields 
that $f\in D(K)$ if and only if there exists a sequence $(\varphi_j)$ 
of continuous functions on $K$ with 
$$\sup_{k\in K} \sum |\varphi_j (k)| < \infty\ \text{ and }\ 
f= \sum\varphi_j\ \text{ pointwise.} 
\leqno(1)$$

Now defining $\|f\|_D = \inf \{\sup_k \sum |\varphi_j(k)| : (\varphi_j)$ 
is a sequence of continuous functions on $K$ satisfying (1)$\}$, it is 
easily seen that $D(K)$ is a Banach algebra, and of course $D(K)\subset 
B_1(K)$ where $B_1(K)$ denotes the first Baire class of bounded functions on 
$K$, i.e., the space of all bounded functions on $K$ which are the limit of a 
pointwise convergent sequence of continuous functions on $K$. 

The primary applications of $D(K)$ in analysis seem to occur in the case 
where $K$ is compact.  For example, a separable Banach space $X$ contains 
a subspace isomorphic to $c_0$ if and only if there is an $f$ in $X^{**} 
\sim X$ with $f\mid K$ in $D(K)$, where $K$ is the unit ball of $X^*$ 
in its $\omega^*$-topology (cf.\ \cite{HOR}, \cite{R1}). 
Using invariants of $D(K)$, it is proved in \cite{R1} that $c_0$ embeds 
in $X$ provided $X$ is non-reflexive and $Y^*$ is weakly sequentially 
complete for all subspaces $Y$ of $X$. 
For applications to spreading models in Banach spaces, see 
\cite{F1}, \cite{F2} and \cite{R3}. 

We are interested here in the intrinsic properties of $D(K)$, and 
compactness  or completeness 
of $K$ plays no role here; moreover if, e.g., $W$ is an open 
subset of $K$, then $D(W)$ plays a natural role in the study of $D(K)$ 
itself. 
We give several permanence properties of $D(K)$, which may be useful 
in further study, and obtain some results which hopefully illustrate 
the fascinating structure of this Banach algebra. 

For example, we obtain that if $K^{(n)}\ne \emptyset$ for all $n$ 
(where $K^{(n)}$ is the $n^{th}$ derived set of $K$), then the functions 
$\varphi$ on $\complex$ which operate on $D(K)$ are precisely those which 
are locally Lipschitz (Proposition~2.8 below). 
($\varphi$ operates on $D(K)$ if $\varphi\circ f\in D(K)$ for all 
$f\in D(K)$. If $K^{(n)} = \emptyset$ for some $n$, then {\it every\/} 
bounded function on $K$ belongs to $D(K)$.) 

We next give a certain oscillation invariant which gives a useful lower bound 
for $D$-norms. 
For $\varphi:K\to [-\infty,\infty]$ an extended real-valued function, 
let $U\varphi$ denote the upper semi-continuous envelope of $\varphi$; 
$(U\varphi)(x) = \olim_{y\to x} \varphi(y)$ for all $x\in K$. 
(We use non-exclusive lim~sups; thus equivalently, $U\varphi (x) = \inf_U
\sup_{y\in U} \varphi (y)$, the inf over all open neighborhoods of $X$.) 
Now for $f:K\to \complex$, we define $\uosc f$, the lower oscillation of $f$, 
by 
$$\uosc f(x)  = \olim\limits_{y\to x} |f(y)-f(x)|\ \text{ for all }\ x\in K\ . 
\leqno(2)$$ 
Then we define $\osc f$, the oscillation of $f$, by 
$$\osc f = U\uosc f\ . 
\leqno(3)$$ 

Now let $n\ge1$ and $(\varep_i)$ be a sequence of non-negative integers, 
either infinite or of length at least $n$. 
We define the oscillation sets $\os_j (f,(\varep_i))$ by induction as follows. 

First, for $\varep\ge0$, set $\os (f,\varep) = \{x\in K: \osc f(x)\ge\varep\}$; 
then let $\os_1(f,(\varep_i)) = \os (f,\varep_1)$. 
If $1\le j<n$ and $\os_j(f,(\varep_i)$ has been defined, let 
$\os_{j+1}(f,(\varep_i)) = \os(f\mid W,\varep_{j+1})$ where 
$W=\os_j (f,(\varep_i))$. 

If $\varep>0$ is given and $\varep_i =\varep$ for all $i$, we set 
$\os_n(f,\varep)=\os_n (f,(\varep_i))$ for all $n$. 
It is also useful to set $\os_0 (f,(\varep_i))=K$. 
We then have the following result, refining a similar lemma in \cite{HOR}. 

\proclaim{Lemma 1} 
Let $f:K\to\complex$ be a given $D$-function and suppose $n$ and 
$(\varep_i)_{i=1}^n$ are given with $\os_n (f,(\varep_i))\ne\emptyset$. 
Then $\sum_{i=1}^n \varep_i \le \|f\|_D$. 
(We prove this in Lemma~1.8 below.) 
\endproclaim 

It turns out that this invariant characterizes a larger class of functions,  
termed $B_{1/4}(K)$, and also yields the $D$-norm of all simple $D$-functions. 

\proclaim{Definition} 
$B_{1/4}(K)$ denotes the class of all functions $f:K\to\complex$ so that 
there exists a sequence $(f_n)$ in $D(K)$ and a $\lambda<\infty$ so that 
$$f_n\to f\ \text{ uniformly and }\ \|f_n\|_D\le\lambda
\ \text{ for all }\ n\ . 
\leqno(4)$$ 
\endproclaim 

For $f\in B_{1/4}(K)$, we define $\|f\|_{B_{1/4}}$ as the infimum of the 
$\lambda$'s satisfying (4) for some $(f_n)$ in $D(K)$. It is easily seen 
that this infimum is obtained; and that moreover $B_{1/4}(K)$ is a Banach 
algebra. $B_{1/4}(K)$ is introduced in \cite{HOR}, where it is shown that 
e.g., if $K= [0,1]$, $B_{1/4}(K)\sim D(K)\ne \emptyset$. 
In fact, Lemma~1 easily yields that if $f\in B_{1/4}(K)$ and 
$\os_n(f,\varep_i) \ne \emptyset$, then $\sum_{i=1}^n \varep_i 
\le \|f\|_{B_{1/4}}$. 

We obtain here the following characterization of $B_{1/4}(K)$. 

\proclaim{Theorem 2} 
Let $f:K\to\complex$ be a given function. Then the following are equivalent. 

{\rm (a)} $f\in B_{1/4}(K)$. 

{\rm (b)} There exists a sequence $(\varphi_n)$ of simple $D$-functions 
with $\varphi_n\to f$ uniformly and $\sup_n\|\varphi_n\|_D<\infty$.

{\rm (c)} There exists a constant $\lambda$ so that for any $n$ and 
sequence $(\varep_i)_{i=1}^n$ of non-negative numbers, 
$$\text{if }\ \os_n(f,(\varep_i)) \ne\emptyset,\ \text{ then }\  
\sum_{i=1}^n \varep_i \le\lambda\ .
\leqno(5)$$ 

Moreover, if $f$ is real-valued and $\beta$ is the best constant $\lambda$ 
satisfying $(5)$ for all $n$ and sequences $(\varep_i)$, then 
$$\tfrac12 \bigl( \|f\|_\infty +\beta\bigr) \le \|f\|_{B_{1/4}} 
\le \|f\|_\infty +3\beta\ . 
\leqno(6)$$ 
\endproclaim 

This result answers  Problem 8.2 in \cite{HOR} in the affirmative. 
(The results in \cite{HOR}, using an equivalent invariant given there, 
yield the necessity of the condition (5) above.) 
The proof of Theorem~2 uses a basic invariant for $D$-functions, the 
transfinite oscillations, which we now recall (cf.\ \cite{R1}). 

\proclaim{Definition} 
Let $f:K\to \complex$ be a given function, $\alpha$ an ordinal. 
We define the $\alpha^{th}$ oscillation of $f$, $\osc_\alpha f$, by 
induction, as follows: set $\osc_0f\equiv0$. 
Suppose $\beta>0$ is a given ordinal and $\osc_\alpha f$ has been defined 
for all $\alpha<\beta$. If $\beta$ is a successor, say $\beta=\alpha+1$, 
we define 
$$\widetilde{\osc}_\beta f(x) = \olim\limits_{y\to x} 
\left( |f(y)-f(x)| + \osc_\alpha f(y)\right)\ \text{ for all }\ 
x\in K\ . 
\leqno(7)$$ 

If $\beta$ is a limit ordinal, we set 
$$\widetilde{\osc}_\beta f= \sup_{\alpha <\beta} \osc_\alpha f\ . 
\leqno(8)$$ 
Finally, we set $\osc_\beta f= U\widetilde{\osc}_\beta f$.
\endproclaim 

Evidently we have that $\widetilde \osc_1 f= \uosc f$ and $\osc_1 f= \osc f$.  
The $\alpha^{th}$ oscillation is similar to a previous invariant, 
introduced by A.S.~Kechris and A.~Louveau \cite{KL}, 
for $f:K\to \real$ (and $K$ compact, 
which is not essential), which is denoted $v_\alpha (f)$. 
$v_\alpha (f)$ is defined in exactly the same way, but with the absolute value 
signs omitted in (7). 
The transfinite oscillations appear to be more appropriate than the 
$v_\alpha (f)$'s, for the study of the Banach space properties of $D(K)$ 
and related objects such as $D(K)$. 
In the present paper, we really only use the finite oscillations $\osc_1 f, 
\osc_2 f,\ldots$ and of course the natural limit of these, $\osc_\omega f$. 
The fundamental connection between the oscillation sets and the oscillation  
functions is given in Lemma~3.8, which immediately yields the following 
result. 

\proclaim{Lemma 3} 
Let $f:K\to\complex$ be a given function. 
Then $\|\osc_\omega f\|_\infty = \sup \sum_{i=1}^n \varep_i : n\ge1$ 
and $(\varep_i)_{i=1}^n$ are positive numbers with $\os_n (f,(\varep_i)) 
\ne \emptyset$. 
\endproclaim 

We thus obtain 

\proclaim{Corollary 4} 
Let $f:K\to \complex$ be a bounded function. 
Then $f$ belongs to $B_{1/4}(K)$ if and only if $\osc_\omega f$ is a 
bounded function. 
\endproclaim 

(This result is also obtained by V.~Farmaki (for $K$ compact) in \cite{F1} 
using different methods.) 
In a subsequent paper \cite{R2}, we will exploit the higher order 
oscillations to study the transfinite analogues of $B_{1/4}(K)$. 

\demo{Remark}
V.~Farmaki and A.~Louveau have recently proved the following remarkable 
identity for real-valued $f$ \cite{FL}: 
$$\|f\|_{B_{1/4}} = \big\|\,|f| + \widetilde{\osc}_\omega f\big\|_\infty\ .$$
See the Remark following Theroem 4.3 below, for some discussion of this 
result (which I learned about after writing the first draft of this article).
\enddemo 

The finite oscillations of a function are our basic tool here in studying 
the simple $D$-functions and the following natural class. 

\proclaim{Definition} 
Let $SD(K)$ be the closure of the set of simple $D$-functions in $D(K)$. 
Members of $SD(K)$ are called strong $D$-functions. 
\endproclaim 

Since the simple $D$-functions are an algebra, it follows immediately that 
$SD(K)$ is a Banach sub-algebra of $D(K)$. 
The next result itself lies rather below the surface (unlike the case of 
$D(K)$ itself). 

\proclaim{Theorem 5} 
Let $f:K\to\complex$ be a strong $D$-function. 
\smallskip
\iitem{\rm (a)} $f$ is a (complex) difference of semi-continuous strong 
$D$-functions. 
\iitem{\rm (b)} $|f|$ is a strong $D$-function. 
\endproclaim 

In fact, for (a), we prove (see Proposition 4.8 below) that 
{\it if $f$  is real-valued in $SD(K)$, and $\varep>0$ is given, there exist 
$u,v$ non-negative lower semi-continuous functions in $D(K)$ with 
$f=u-v$ and $\|u+v\|_\infty < \|f\|_D+\varep$.}  

It follows from results in \cite{HOR} that if e.g., $K$ is compact 
uncountable, then the $D$ and $B_{1/4}$-norms are not equivalent on $D(K)$. 
This does not occur on $SD(K)$, however. 
The techniques which prove Theorems~2 and 5 yield the following. 

\proclaim{Corollary 6} 
Let $f\in SD(K)$, with $f$ real-valued. 
Then $\|f\|_{B_{1/4}} = \|f\|_D$. 
\endproclaim

We now indicate the organization and contents of this work. 
The exposition is intended to be comprehensive; thus readers interested 
mainly in the proofs of the results stated so far, can skip quite a bit. 
For example, Theorem~2 follows from the results of Section~3, through 
Corollary~3.8, and the results of Section~4 through Theorem~4.3. 
The same applies to Lemma~3 and Corollary~4. 
Theorem~5 and Corollary~6 require the further development in Section~4, 
through Theorem~4.14. 
The development in Sections~1 and 2 does not require the transfinite 
oscillations, and is perhaps more elementary because of this. 

Section~1 consists of preliminary results, most of which are not 
explicitly used in the sequel, though they motivate much of what follows. 
Fix $K$ a metric space; we let $B_{1/2}(K)$ denote the uniform closure 
of $D(K)$ in the bounded functions on $K$. 
The development through Corollary~1.5 deals with the fact that $D(K)$ 
is a semi-simple Banach algebra with $K$ densely embedded in its 
maximal ideal space, $\Omega$. 
$B_{1/2}(K)$ can be identified with $C(\Omega)$, and then properties of 
$B_{1/2}(K)$, proved later on, yield that $\Omega$ is totally disconnected, 
with the simple functions in $C(\Omega)$ already belonging to $D(K)$. 

Proposition 1.6 deals with extension issues, and yields the (perhaps 
surprising) result that if $W$ is an arbitrary subset of $K$ and $g$ is 
real-valued in $D(W)$, then $g$ extends to a function  $\tilde g$ in 
$D(K)$, with $\|\tilde g\|_{D(K)} = \|g\|_{D(W)}$. 
If $W$ belongs to $\D$, the algebra of sets generated by the open subsets 
of $K$, then $g\chix_W$ already belongs to $D(K)$, and if $W$ is a 
difference of closed sets, then $\|g\cdot\chix_W\|_{D(K)} \le 
2\|g\|_{D(W)}$. 

Lemma 1.8 yields the lower bound for the $D$-norm given in Lemma~1. 
The proof is similar to an argument in \cite{HOR}; however 1.8, 
for real-valued functions, also follows from the independent development 
in Section~3. 
Lemma~1.8 is used to characterize $B_{1/2}(K)$ in Proposition~1.9, 
recapturing a result given in \cite{HOR}; and then Proposition~1.10 
yields that simple functions in $B_{1/2}(K)$ are already  $D$-functions. 

Proposition 1.4 yields the following localization principle: 
if a function locally belongs to $D$ with the local $D$-norms uniformly 
bounded then it belongs to $D$. 
The proof is achieved via partitions of unity; an alternative argument, 
involving transfinite oscillations, follows from the results of Section~3.  
We end Section~1 with a proof of Baire's famous theorem: 
every lower semi-continuous function is the limit of an increasing 
sequence of continuous functions (Proposition~1.18). 

Section 2 solves the following problem: 
{\it given a set $A$, find $\|\chix_A\|_D$\/}. 
It is easily seen (as shown in Section~1) that $\chix_A$ is a 
$D$-function if and only if $A$ belongs to $\D$; 
in turn, this happens if and only if $A$ is a finite disjoint union 
of differences of closed sets. 
We then obtain the solution to the following problem, as a by-product:  
given $A$ in $\D$, find the {\it smallest\/} integer $k$ with $A$ 
a union of $k$ disjoint differences of closed sets $W_1,\ldots,W_k$. 
(Actually, our results here also hold for arbitrary Hausdorff spaces.) 
The solution is as follows: 
for $A\subset K$, let $\partial'A=\partial A$, the boundary of $A$, 
and let $\partial^nA$, the $n^{th}$ boundary of $A$, equal the boundary 
of $A\cap \partial^{n-1}A$, {\it relative\/} to $\partial^{n-1}A$ (for $n>1$).
Now define $i(A)$, the Baire-index of $A$, as the largest $n$ with 
$\partial^nA\ne\emptyset$; if no such $n$ exists, set $i(A)=\infty$. 
Theorem~2.2 now yields that 
{\it $A\in \D$  iff $\chix_A \in D(K)$ iff $n \dfeq i(A)<\infty$. 
If $A\cap \partial^n A=\emptyset$, $\|\chix_A\|_D=n$; 
if $A\cap\partial^n A\ne\emptyset$, then $\|\chix_A\|_D = n+1$. 
Finally, there exist $k$ disjoint differences of closed sets 
$W_1,\ldots,W_k$ with $A= \bigcup_{i=1}^k W_i$ and 
$\|\chix_A\|_D = \sum_{i=1}^k \|\chix_{W_i}\|_D$ where 
$k= \bigl[\tfrac{n+1}2\bigr]$ if $A\cap\partial^n A =\emptyset$, 
$k= \bigl[\tfrac{n}2 +1\bigr]$ if $A\cap\partial^nA \ne \emptyset$. 
If $A$ is a disjoint union of $\ell$ differences of closed sets, 
then $\ell\ge k$.} 

We then show in Proposition 2.6 that if $K$ is any metric space with all its 
finite derived sets non-empty, then for all $n$, there exist subsets $A,B$ 
of $K$ with $i(A)= i(B)=n$ and $\|\chix_A\|_D =n$, 
$\|\chix_B\|_D= n+1$. 
It then follows, via Corollary~2.7, that $D(K)\ne B_{1/2}(K)\ne B_1(K)$ 
for any such $K$. 
We conclude Section~2 with the result that for such $K$, 
$\varphi :\complex\to\complex$ operates on $D(K)$ iff $\varphi$ is 
locally Lipschitz (Proposition~2.8). 
(The fact that locally Lipschitz functions operate on $D(K)$, is due 
jointly to F.~Chaatit and the author \cite{C}.) 

Section 3 treats some properties and first applications of the transfinite 
oscillations. After summarizing several elementary properties in 
Proposition~3.1, we give the basic structure theorem concerning the 
transfinite oscillations, in Theorem~3.2. 
The result yields the following information: 
{\it If $f:K\to \real$ is bounded, there is an ordinal $\gamma$ with 
$\osc_\gamma f = \osc_{\gamma+1}f$; 
we denote the least such ordinal by $i_Df$. 
Then if $\alpha = i_D f$, $f\in D(K)$ iff $\osc_\alpha f$ is bounded, and 
then $\|f\|_D = \|\, |f|+\osc_\alpha f\|_\infty$; 
moreover if then $\lambda = \|f\|_D$ and $u= {\lambda -\osc_\alpha f+f\over2}$,
$v= {\lambda-\osc_\alpha f-r\over2}$, 
$u,v$ are both non-negative lower semi-continuous, and of course $f=u-v$, 
$\|u+v\|_\infty = \|f\|_D$. 
Moreover $\|f\|_{qD} = \|\osc_\alpha f\|_\infty$, where 
$\|f\|_{qD} = \inf\{\|f-\varphi\|_D :\varphi \in C_b(K)\}$. 
}

If we define, for fixed $x\in K$, $\|f\|_{qD(x)} = \inf \|f\mid U\|_{qD(U)}$, 
the inf over all open neighborhoods $U$ of $x$, we obtain in 
Corollary~3.5 that $\osc_\alpha f(x) = \|f\|_{qD(x)}$ (where $\alpha=i_Df$). 
Moreover the ``quotient $D$-semi-norm'' $\|f\|_{qD}$ is attained; 
i.e.,  $\|f\|_{qD} = \|f-\varphi\|_D$ for some $\varphi \in C_b(K)$. 

The technical Lemma~3.6 characterizes the finite oscillations in terms of 
the oscillation sets $\os_k(f,(\varep_i))$, and yields Lemma~3 above, as 
an immediate corollary (Corollary~3.7). 
The development up to this point easily yields that functions of finite 
Baire-index are $D$-functions, as shown in \cite{CMR} by different methods. 
To formulate this concept, first define, for $\varep>0$, the 
$\varep$-Baire index of a bounded $f:K\to\complex$, $i_B(f,\varep)$, as 
the largest $n$ with $\os_n(f,\varep) \ne\emptyset$ (or set 
$i_B(f,\varep) = \infty$ if there is no such $n$). 
Proposition~1.9 yields that $f\in B_{1/2}(K)$ iff $i_B(f,\varep) <\infty$ 
for all $\varep>0$. 
Now set $i_B(f) = \sup_{\varep>0} i_B(f,\varep)$. 
We say $f$ is of finite Baire-index if $i_B(f)<\infty$. 
(It is easily seen that simple $D$-functions are of finite Baire-index; 
cf.\ the proof of Proposition~1.10). 
Now Corollary~3.8  yields that if $n=i_B(f)<\infty$, then $i_D(f)\le n$ 
and if $f$ is real valued, $\|f\|_D \le (2n+1)\|f\|_\infty$; 
this estimate is moreover best possible, in general 
(cf.\ the discussion following the proof of Proposition~3.10). 
A little more work is required to recapture the result of \cite{CMR} 
that functions of finite Baire-index are actually strong $D$-functions; 
this is done in Section~4. 

The rest of Section 3 illustrates the preceding general results with the 
computation of finite oscillations and $D$-norms of some simple functions. 
For example, Proposition~3.10 yields the precise description of 
$\osc_nf$ if $f=\chix_A$ with $n=i(A)$; it turns out that the 
boundaries $\partial^j A$, $0\le j\le n$, are the only invariant.  
Proposition~3.10 also yields another proof of the value of $\|\chix_A\|_D$ 
given by Theorem~2.2, as well as the fact that $\|\chix_A\|_{qD}=n$. 
Finally, Proposition~3.11 yields the $D$-norm of a certain natural 
class of simple $D$-functions, namely the functions $f$ so that letting 
$n=i_B(f)$ and defining $K_0=K$ and $K_{j+1} = \{x\in K_j: f\mid K_j$ 
is discontinuous at $x\}$, then $f\mid K_j\sim K_{j+1}$ is constant, for all 
$0\le j\le n$. 

Section 4 deals with the proofs of the results stated at the beginning.  
We show that continuous bounded functions belong to $SD$ in Proposition~4.1; 
it follows that if $K= K_0 \supset K_1\supset\cdots\supset K_n$,  
$\varphi :K\to\complex$ is bounded, and $\varphi\mid K_i\sim K_{i+1}$ 
is continuous for all $i$, then also $\varphi\in SD$ and $i_B(\varphi)\le n$ 
(Proposition~4.2). 
Theorem~4.3 then yields the characterizations of $B_{1/4}$ given by 
Theorem~2 and Corollary~4 (in virtue of Lemma~3). 
We then give further permanence properties of $SD$; thus if $f\in SD$, 
then $(\osc_nf)$ converges uniformly (Proposition~4.4), and 
hence $i_Df\le \omega$ and $\widetilde{\osc}_\omega f 
= \osc_\omega f$ (Corollary~4.5). 
We next obtain that the $B_{1/4}$ and $D$-norms coincide for real $f$ 
in $SD$ in Corollary~4.6 (Corollary~6 above), and then establish that 
every $SD$ function is a difference of lower semi-continuous $SD$-functions 
in Proposition~4.8 (Theorem~5a above). 

Next, we assemble some tools to prove Theorem~5b. 
We first recall the basic index result given in \cite{CMR}; 
{\it if $f,g$ are given bounded functions and $\varep>0$ is given, 
then $i_B(f+g,\varep)\le i_B(f,{\varep\over2}) + i_B(g,{\varep\over2})$}  
(Lemma~4.10),
and use this and the preceding development to recapture 
the result in \cite{CMR} that functions of finite index belong to 
$SD$ (Corollary~4.11). 
We then recall the class $B_{1/2}^0(K)$, given in \cite{CMR}, of bounded 
functions $f$ with $\lim_{\varep\to0}\varep i_B(f,\varep)=0$, and give the 
result of \cite{CMR} that $B_{1/2}^0$ is a linear space and complex lattice, 
containing $SD$, whose semi-continuous members belong to $SD$ 
(Proposition~4.12). 

We then prove that $f\in SD$ implies $|f|\in SD$, as follows: 
first we show, in Lemma~4.13, that for $f\in SD$, there is an 
upper semi-continuous strong $D$-function $F\ge0$ 
with $F+|f|$ upper semi-continuous. 
It follows that $F+|f|$ is strong $D$ since it belongs to $B_{1/2}^0$, 
whence $|f|$ is strong $D$. 
(At the end of Section~4, we give an intrinsic characterization of $SD$ 
which also yields Theorem~5b, in a perhaps more natural way.) 

We next give a variety of examples of functions failing the various 
criteria for $SD$ given above. 
For example, we construct in Example~3 a $D$-function $g$ with $i_D = 
\omega+1$ (so $g\notin SD$ by Corollary~4.5). 
In Example~5, we construct a function in $B_{1/2}^0 (K)\cap D(K)\sim 
SD(K)$ (for suitable $K$). 

The rest of Section 4 deals with an intrinsic criterion characterizing 
strong $D$. 
We generalize the sets $\os_n(f,\varep)$ to $\os_n (f,\alpha,\varep)$ for a 
given ordinal $\alpha$, by setting $\os_{n+1} (f,\alpha,\varep) = 
\{x\in L:\osc_\alpha f\mid L\ge \varep$ with $L= \os_n(f,\alpha,\varep)\}$, 
and then define $i(f,\alpha,\varep) = \sup \{n\ge 0: \os_n(f,\alpha,\varep) 
\ne\emptyset\}$. 
Corollary~4.18 yields that $\varep i(f,\alpha,\varep) \le \|\osc_{\alpha\omega} 
f\|_\infty$ for all $f,\varep>0$ and ordinals $\alpha$. 
We prove a generalization  of Lemma~4.10 in Lemma~4.20 (with essentially 
the same proof as that of 4.10 given in \cite{CMR}), 
and use this and preceding results to establish 

\proclaim{Theorem 4.19} 
A bounded function $f$ is strong $D$ if and only if $\lim_{\varep\to0} 
\varep i(f,\omega,\varep)=0$ and $i_D (f\mid W)\le \omega$ for 
all closed sets $W$. 
\endproclaim 

Most of the results given here were presented in topics courses at 
the University of Texas at Austin during the academic years 1991-1993. 
I am grateful to the course-participants for their patience and support 
of this work. 

\head \S1. Preliminaries.\endhead 

We begin with some elementary algebraic and lattice properties of $D(K)$. 

\proclaim{Proposition 1.1} 
Let $K$ be a metric space 

{\rm (a)} $D(K)$ is a commutative Banach algebra with identity. 

{\rm (b)} If $f\in D(K)$, then $|f|\in D(K)$ and $\|\, |f|\,\|_D\le \|f\|_D$.  

{\rm (c)} If $f$ is in $D(K)$ and $\inf_k|f(k)| >0$, then 
$1/f \in D(K)$. 
\endproclaim 

\demo{Proof} 
We omit the routine proof of (a). 
To see (b), let $\varep>0$, and choose $(f_n)$ in $C_b(K)$ with 
$f_n\to f$ pointwise and 
$$|f_1| + \sum |f_{n+1} - f_n| \le \|f\|_D +\varep 
\leqno(1)$$ 
But then $|f_n|\to |f|$ pointwise, and 
$$|f_1| + \sum \big|\ |f_{n+1}| - |f_n|\ \big| \le \|f\|_D +\varep 
\leqno(2)$$ 
by (1). 

Hence $\|\, |f|\,\|_D \le \|f\|_D+\varep$ for all $\varep>0$, so (b) 
is proved. 
To prove (c), suppose first that $f\ge0$, and say 
$\delta = \inf_{k\in K} f(k) >0$. 
Now choose $(f_n)$ in $C_b(K)$ with $f_n\to f$ pointwise, satisfying (1). 
We may obviously assume the $f_n$'s are real-valued. 
Now define $(g_n)$ by 
$$g_n = f_n \vee \delta\ \text{ for all }\ n\ . 
\leqno(3)$$ 
Then we have that $g_n\to f$ pointwise, and moreover 
$$|g_{n+1} -g_n|  \le |f_{n+1}-f_n|\ \text{ for all }\ n\ . 
\leqno(4)$$ 
But evidently $\frac1{g_n} \to \frac1f$ pointwise, and 
$$\sum_{n=1}^\infty \left| \frac1{g_{n+1}} - \frac1{g_n}\right| 
= \sum_{n=1}^\infty {|g_{n+1}-g_n| \over g_{n+1} g_n} 
\le {1\over \delta^2} \sum |f_{n+1} - f_n|\ . 
\leqno(5)$$ 
Evidently we then have 
$${1\over |g_1|} + \sum_{n=1}^\infty \Big| {1\over g_{n+1}} - {1\over g_n}
\Big| \le {1\over\delta} + {1\over\delta^2} 
\bigl( \|f\|_D + \varep\bigr)\ . 
\leqno(6)$$ 

Hence ${1\over f}\in D$ and in fact 
$$\Big\| {1\over f}\Big\|_D  \le {1\over\delta} +{1\over\delta^2} \|f\|_D\ . 
\leqno(7)$$
Finally, if $f$ is arbitrary, we have that $f\bar f = |f|^2$ is in 
$D(K)$ by (a), so by what we have proved, 
$1/(f\bar f) \in D(K)$, hence ${1\over f} = f{1\over f\bar f}$ belongs also 
to $D(K)$.\qed
\enddemo 

The following is an immediate consequence of 1.1(b). 

\proclaim{Corollary 1.2} 
$D(K)$ is a real function lattice. 
That is, if $f,g$ are real-valued in $D(K)$, then $f\vee g$, $f\wedge g$ 
belong to $D(K)$, and 
$$\eqalign{
\|f\vee g\|_D & \le \|f\|_D + \|g\|_D\cr 
\|f\wedge g\|_D & \le \|f\|_D + \|g\|_D\ .\cr}$$
\endproclaim 

\demo{Proof} 
This is immediate from 1.1(b) and the standard formulas 
$$\leqalignno{
f\vee g&= {(f+g)+ |f-g|\over 2} &(8)\cr 
f\wedge g &= {(f+g)-|f-g|\over 2}&(9)\cr}$$
\enddemo 

\demo{Remark} 
Evidently 1.1(b) yields that in $D(K)$, $f\to |f|$ is continuous at $f=0$, 
and so also the lattice  operations are continuous at $0$. 
However {\it none\/} of these operations are continuous at other points 
in $D(K)$, for general $K$. 
\enddemo 

We may define  an involution $*$ on $D(K)$ by $f^* = \bar f$ for all 
$f$ in $D$. Evidently we have that $\|f^*\|_D = \|f\|$ for all $f$. 
It then follows immediately from 1.1(c) that $D(K)$ is a ``completely 
symmetric ring'' as defined in \cite{N}. 
That is, we let $\Omega$ denote the ``structure space,'' 
or ``maximal ideal space,'' of $D(K)$, 
namely the set of all multiplicative linear  functionals on $D(K)$, 
endowed with the topology of pointwise convergence.
For $f$ in $D(K)$, we define $\hat f$ on $\Omega$ by $\hat f(\omega)  = 
\omega (f)$ for all $\omega\in\Omega$. 
Evidently $\hat f$ is then in $C(\Omega)$. 
($\hat f$ is sometimes called the Gelfand transform of $f$.) 

\proclaim{Corollary 1.3} 
$$\widehat{f^*} = \overline{(\hat f)}\ \text{ for all }\ f\ \text{ in }\ 
D(K)\ .$$ 
\endproclaim 

\demo{Proof} 
Since $\bone +f^*f$ is invertible in $D(K)$ for all $f$, by 1.1(c), this 
follows directly from a result in \cite{N}. 
The proof, however, is even simpler here than in \cite{N}, and is as follows.  
Since it is evident that $f$ belongs to $D(K)$ if and only if its real and 
imaginary parts belong, it suffices to prove that if $f$ in $D$ is 
real-valued (on $K$), then $\hat f$ is real-valued. 
(For then writing $f=u+iv$, with $u=\Re f$, $v=\Im f$, we have that 
$\widehat{f^*} = \widehat{u^*} - i\widehat{v^*} = \hat u- i\hat v
= \overline{(\hat f)}$.) 

By the general theory of commutative Banach algebras, the complex number 
$\lambda$ belongs to the range of $\hat f$ if and only if $f-\lambda\bone$ 
is not invertible. 
But if $\lambda =\alpha+i\beta$ say, with $\beta\ne0$, then 
$|f-\lambda\bone| \ge |\beta|>0$, hence $f-\lambda\bone$ is invertible 
in $D(K)$, by 1.1(c).\qed 
\enddemo 

It follows from 1.3 and the Stone-Weierstrass theorem that $\hat D$ 
is {\it dense\/} in $C(\Omega)$; indeed, $\hat D$ is a point-separating 
conjugation closed subalgebra of $C(\Omega)$ which contains the constants. 
We {\it define\/} $B_{1/2}(K)$ to be the uniform closure of $D$ in $K$, 
endowed with the sup-norm.
We show below that $B_{1/2}(K)$ can be canonically identified with 
$C(\Omega)$. 

We shall give below several alternate characterizations of $B_{1/2}(K)$, 
which yield the following result (see Propositions~1.9 and 1.10). 
(A function is called {\it simple\/} if its range is a finite set.) 

\proclaim{Proposition 1.4} 
Let $K$ be a given metric space. 

{\rm (a)} The set of simple functions in $B_{1/2}(K)$ is dense in $B_{1/2}(K)$. 

{\bf (b)} Every simple function in $B_{1/2}(K)$ belongs to $D(K)$. 
\endproclaim 

We now deduce some properties of $\Omega$. 

\proclaim{Corollary 1.5} 
Let $K$ be a metric space, $\Omega$ the structure space of $D(K)$, 
and $\tau :K\to\Omega$ the canonical map. 
Then $\tau K$ is dense in $\Omega$, and $\Omega$ is totally disconnected. 
In fact, given $A$, $B$ disjoint closed subsets of $\Omega$, there exists 
a $\{0,1\}$-valued $f$ in $D(K)$ with $\hat f\equiv 1$ on $A$ 
and $\hat f \equiv 0$ on $B$. 
\endproclaim 

\demo{Proof} 
Suppose it were false that $\overline{\tau K} = \Omega$. 
Then we could choose $g$ in $C(\Omega)$ and $\omega\in\Omega$ with $g\equiv0$ 
on $\tau K$ and $g(\omega)=1$. 
Since $\hat D$ is dense in $C(K)$, we can choose $f$ in $D(K)$ with 
$$\|\hat f-g\|_\infty < \tfrac13\ . 
\leqno(10)$$ 
But (10) yields that 
$$
|\hat f(\omega)| >\tfrac23\  
\text{ and }\ 
|\hat f(\tau (k))| < \tfrac13\ \text{ for all }\  k\in K\ . 
\leqno(11)$$ 
Hence, $|f-\hat f(\omega)\cdot \bone | \ge \frac13$ on $K$, so by 
Proposition~1.1(c), $f-\hat f(\omega)\cdot \bone$ is invertible. 
But by general Banach algebra theory, $f-\hat f(\omega)\cdot\bone$ is 
singular. This contradiction establishes our first assertion. 

Next, it follows that we may extend the map $\Lambda :D(K)\to C(\Omega)$, 
defined by $\Lambda (f) =\hat f$, 
to a bijection, also denoted $\Lambda$, between $B_{1/2}(K)$ and $C(\Omega)$. 
$\Lambda$ thus extended will be an algebraic isometry, in the sup-norms. 
Indeed, let $f$ in $B_{1/2}(K)$; we claim there is a continuous $g$ on 
$\Omega$ with $g(\tau k) = f(k)$, all $k\in K$. 
Once the claim is proved, we have of course that $g$ is unique, since 
$\tau K$ is dense in $\Omega$; so we define $\hat f=g$. 
The uniqueness of $g$ also yields that $\Lambda$  is an into-algebraic 
isometry, and the fact that $\hat D$ is dense in $C(\Omega)$ (as 
noted above) yields that $\Lambda :B_{1/2}(K)\to C(\Omega)$ is a 
surjection. To prove the claim, just choose $(f_n)$ in $D(K)$ with 
$f_n\to f$ uniformly on $K$. 
But then trivially $(\hat f_n)$ converges uniformly on $\tau K$, and 
since $(\hat f_n)\subset C(\Omega)$ and $\overline{\tau K} = \Omega$, 
$(\hat f_n)$ converges uniformly on $\Omega$, to a continuous $g$, proving 
the claim.  

Now in general, a compact Hausdorff space $X$ is totally disconnected if 
and only if the simple members of $C(X)$ are dense in $C(X)$. 
Thus we now obtain that $\Omega$ is totally disconnected, by 
Proposition~1.4(a). 

Finally, a standard compactness argument shows that if $A$ and $B$ are 
closed disjoint subsets of $\Omega$, there exists a clopen set $E$ with 
$A\subset E$ and $B\cap E=\emptyset$. 
Thus $\chix_E$ is continuous on $\Omega$, and so there is an $f$ in 
$B_{1/2}(K)$ with $\hat f= \chix_E$. 
$f$ is of course $\{0,1\}$-valued on $K$, since $\Lambda$ is an 
algebra-isomorphism, and thus $f$ belongs to $D(K)$ by 
Proposition~1.4(b).\qed
\enddemo 

We next treat the extension of $D$-functions. 
Given $A\subset K$ and $f:A\to\complex$, we use the notation $f\chix_A$ 
to denote the function on $K$ which is zero off $A$ and agrees with $f$ on 
$A$. Also, we say that $A\subset K$ is a difference of closed sets, or a 
DCS, if there exist closed subsets $A_1,A_2$ of $K$ with $A= A_1\sim A_2$. 
The class of all such sets is also denoted DCS. 

\proclaim{Proposition 1.6} 
Let $W\subset K$ be non-empty and $f$ in $D(W)$ be given. 

\item{\rm (a)}  There exists a $g$ in $D(K)$ with $g\mid W=f$. Moreover,

\iitem{\rm (i)} if $f$ is real-valued, $g$ may be chosen (real-valued) 
with $\|g\|_{D(K)} = \|f\|_{D(W)}$; 

\iitem{\rm (ii)} if $W$ is a {\rm DCS} and $\varep >0$ is given, then $g$ 
may be chosen with $\|g\|_{D(K)} \le \|f\|_{D(W)} +\varep$.  

\item{\rm (b)}  If $W$ is a {\rm DCS}, then $\|f\chix_W\|_{D(K)} 
\le 2\|f\|_{D(W)}$, while if $W$ is open, then $\|f\chix_W\|_{D(K)} 
= \|f\|_{D(W)}$.
\endproclaim 

The proof of the qualitative part of (a) is quite simple, and in fact it's 
easily seen that for $f$ real-valued, a $g$ may be chosen extending $f$ 
with $\|f\|_{D(K)} \le 2\|f\|_{D(W)}$. 
We show this first, then prove all the assertions of 1.6, for completeness. 

\demo{Proof of the qualitative part of (a)} 
Suppose first that $u$ is a non-negative lower semi-continuous function 
on $W$ and define $\util$ on $\bar W$ by 
$$\util (k) = \ulim\limits_{\scriptstyle w\to k\atop \scriptstyle w\in W} 
( \dfeq \sup \inf u(V\cap W)\ ,\text{ the sup over all open neighborhoods 
$V$ of $k$).}$$ 
Evidently $\util$ is then lower semi-continuous on $\bar W$, and of 
course $\util\mid W = u$. 

Next, if $f$ in $D(W)$ is real-valued, then by Theorem~3.5 of \cite{R1} 
(indicated also below), we may choose $u,v$ non-negative lower semi-continuous 
with $f=u-v$ and $\|u+v\|_\infty = \|f\|_{D(W)}$. 
Then defining $\util,\vtil$ as above, and setting $\tilde f = 
\util-\vtil$ (on $\bar W$), then evidently $\tilde f\mid W=f$ 
and for $k\in \bar W$, $\util(k) + \vtil(k) \le \olim_{w\to k,\, w\in W} 
u(w) +v(w) \le \|u+v\|_\infty$, 
whence $\|\tilde f\|_{D(\bar W)} = \|f\|_{D(W)}$. 
Now if we simply set $\lambda = \|f\|_{D(W)}$ and define $\uttil = 
\util \cdot \chix_{\bar W} + \lambda\cdot\chix_{\sim W}$, 
$\vttil = \vtil \cdot\chix_{\bar W} +\lambda\cdot \chix_{\sim W}$, 
then setting $g= \uttil - \vttil$, $g\in D(K)$, $g\mid W=f$, and 
$\|g\|_{D(K)} \le 2\|f\|_{D(W)}$, for we easily have that $\uttil$ and 
$\vttil$ are both lower semi-continuous. 
Of course the general complex-valued case now follows immediately. 
\enddemo 

To obtain the assertions (a)(i),(ii) in 1.6, we use the following result. 

\proclaim{Lemma 1.7} 
Let $W$ be a closed subset of $K$, $\lambda<\infty$, and $(\varphi_j)$ a 
sequence of continuous complex-valued functions on $W$ with 
$$\sum |\varphi_j (w)| \le \lambda\text{ for all } w\in W\ .
\leqno(12)$$ 
Then there exists $(\tilde \varphi_j)$ a sequence of continuous functions  
on $K$ with $\tilde\varphi_j \mid W=\varphi_j$ for all $j$ and 
$$\sum |\tilde\varphi_j (k)| \le \lambda\text{ for all } k\in K\ .
\leqno(13)$$ 
\endproclaim 

\demo{Proof of Lemma 1.7}
We require the following linear version of the Tietze extension theorem 
(cf.\ \cite{D}): 
\enddemo

\demo{Fact} 
Let $W$ be a closed subset of a metric space $X$. There exists a linear 
operator $T:C_b(W)\to C_b(X)$ satisfying 

\iitem{(a)} $T\bone_W = \bone_X$.
\iitem{(b)} $Tf|_W = f$ for all $f$.
\iitem{(c)} $\|Tf\|_\infty = \|f\|_\infty$ for all $f$.

\noindent We note that $T$ is then positive, i.e., $f\ge0$ implies 
$Tf\ge0$. It follows further that 
$$|T\varphi| \le T|\varphi|\text{ for all } \varphi\in C_b(W)\ .
\leqno(14)$$ 
Indeed, if $\varphi$ is real-valued, then $\varphi = \varphi^+ - \varphi^-$, 
so 
$$|T\varphi| = |T\varphi^+ - T\varphi^-| \le T\varphi^+ +T\varphi^- 
= T|\varphi|\ .$$ 
If $\varphi$ is complex-valued, fix $x\in X$ and choose $\lambda$ with 
$|\lambda|=1$ so that 
$$|T\varphi|(x) = \lambda T\varphi x 
= T\lambda\varphi (x) 
= (T\Re \lambda\varphi)(x) 
\le T|\Re\lambda\varphi|(x) 
\le T|\lambda\varphi|(x) 
= T|\varphi|(x)\ .$$

Now to obtain 1.7, let $T$ be as in the Fact, and simple set 
$\tilde\varphi_j = T\varphi_j$ for all $j$. Now fixing $n$, 
$$\leqalignno{
\sum_{j=1}^n |\tilde\varphi_j| = \sum_{j=1}^n |T\varphi_j| 
&\le \sum_{j=1}^n T|\varphi_j|\ \text{ by (14)}&(15)\cr 
&= T\biggl( \sum_{j=1}^n |\varphi_j|\biggr) \le \lambda\ ,\cr}$$
the last inequality holding by the positivity of $T$ and (12). 
Of course since $n$ is arbitrary, (13) holds, completing the proof of 1.7. 

We next prove the second assertion in 1.6(b). 
Assume $W$ is open. First note that we may choose a sequence $(\varphi_j)$ 
in $C_b(K)$ with 
$$\varphi_j\ge0\text{ for all } j \text{ and } \sum\varphi_j =\chix_W
\text{ pointwise.}
\leqno(16)$$ 
To see this elementary result, choose $K_1\subset K_2\subset\cdots$ 
closed subsets of $K$ with $W= \bigcup_{j=1}^\infty K_j$. 
Inductively choose a sequence $(f_j)$ in $C_b(K)$ as follows. 
First, choose $f_1:K_1 \to [0,1]$ continuous with $f_1= 1$ on $K_1$, and  
$\overline{\{x:f_1(x)\ne 0\}}\subset W$. 
Suppose $f_j$ chosen with $F_j \dfeq \overline{\{x:f_j(x)\ne0\}} \subset W$. 
Then choose $f_{j+1}:K\to [0,1]$ continuous with $f_{j+1}=1$ on $K_{j+1} 
\cup F_j$  and $\overline{\{x:f_{j+1}(x)\ne0\}} \subset W$. 
Now setting $\varphi_1 = f_1$, $\varphi_j = f_j-f_{j-1}$ for all $j>1$, 
then $(\varphi_j)$ satisfies (16). 

Now let $\varep>0$, and choose $(g_j)$ in $C_b(W)$ with $\sum |g_j| \le 
\|f\|_{D(W)} +\varep$ and $f=\sum f_j$ pointwise. 
But then it follows that $g_j\chix_W \varphi_i$ belongs to $C_b(K)$ for 
all $i$ and $j$, simply because $g_j\chix_W$ is continuous on $W$ and 
bounded on $K$ while $\varphi_i$ is continuous on $K$ and vanishes on 
$K\sim W$. Then $f\cdot\chix_W = \sum_{i,j} g_j \chix_W\varphi_i$ 
and $\sum_{i,j} |g_j\chix_W \varphi_i|\le \|f\|_{D(W)} +\varep$. 
Since $\varep>0$ is arbitrary, the assertion is proved. 

To obtain 1.6(a)(ii), first suppose $W$ is closed, and given $\varep>0$, 
choose $(\varphi_j)$ in $C_b(W)$ with $\sum |\varphi_j|\le \|f\|_{D(W)} 
+\varep$ and $f=\sum\varphi_j$-pointwise. 
Now applying Lemma~1.7, with $(\tilde\varphi_j)$ as in its statement, 
we obtain 1.6(a)(ii) by letting $g=\sum\tilde\varphi_j$ pointwise.
Finally if $W$ is a general DCS, it is easily seen that there is a closed set 
$A$ and an open set $U$ with $W=A\cap U$. 
Thus $W$ is a relatively open subset of $A$, and so letting $h=(f\cdot
\chix_W)\mid A$, $\|h\|_{D(A)} = \|f\|_{D(W)}$. 
Finally, for $\varep>0$, choose $g$ extending $h$ to $K$ as above, with 
$\|g\|_{D(K)} \le \|f\|_{D(W)} +\varep$. 
Now evidently we have that if $W$ is open non-empty, then 
$\|\chix_W\|_D=1$, whence if $W$ is closed, $\|\chix_W\|\le 2$ since 
$\chix_W=1$. Evidently then $\|\chix_W\|\le 2$ if also $W$ is a DCS, using 
the representation given above. 
Thus the first assertion of 1.6(b) follows from (1.6)(ii); 
for given $\varep>0$, choose $g$ in $D(K)$ with $g\mid W=f$ and 
$\|g\|_{D(K)} \le \|f\|_{D(W)} +\varep$. 
Then $g\cdot\chix_W=f\cdot\chix_W$ and so $\|f\cdot\chix_W\|\le 
2\|f\|_{DW} +2\varep$, but $\varep>0$ is arbitrary. 

Finally, for  (a)(i), suppose $f$ is real-valued. 
We proved above that then there exist $\util,\vtil$ non-negative lower 
semi-continuous functions on $\bar W$ with 
$$\|\util+\vtil\|_\infty = \|f\|_{D(W)} \text{ and } 
(\util-\vtil)\mid W=f\ .$$ 
It follows that we may choose $(\varphi_j)$ a sequence in $C_b(\bar W)$ 
with $\util+\vtil=\sum |\varphi_j|$ and $\util-\vtil = \sum\varphi_j$ 
pointwise. 
Now applying Lemma~1.7 with $(\tilde\varphi_j)$ as in its statement, 
we obtain that setting $g=\sum\tilde\varphi_j$, then 
$\|g\|_{D(K)} = \|f\|_{D(W)}$ and $f\mid W=g$, 
completing the proof.\qed
\enddemo 

\demo{Remark} 
Let $\D= \D(K)$ denote the algebra of sets generated by the closed subsets 
of $K$. A standard set-theoretic result yields that a set belongs to $\D$  
if and only if it is a finite disjoint union of members of DCS. 
Thus evidently we obtain that if $W\in \D$ and $f\in D(W)$, then 
$f\cdot\chix_W\in D(K)$. 
We show below that a simple function in $B_{1/2}$ is a $\D$-function, i.e., 
$\D$-measurable; of course it follows that in turn every simple 
$\D$-function belongs to $D(K)$. 
Finally, we note the following fact, whose proof is left to the reader. 
\enddemo 

\demo{Fact} 
Let $W\subset K$. The following are equivalent. 

\iitem{(a)} $W$ is a DCS.
\iitem{(b)} There exist subsets $A$ and $U$ of $K$ with $A$ closed, $U$ open, 
and $W= A\cap U$. 
\iitem{(c)} There exist closed subsets $A$ and $B$ of $A$ with $A\supset B$, 
$B$ nowhere dense relative to $A$, and $W= A\sim B$. 

\noindent The representation in (c) is unique, for then $A$ equals $\bar W$, 
while $B$ equals the boundary of $W$ relative to $A$. 
\enddemo 

We next give a fundamental lower bound for the $D$-norm, refining a 
similar result in \cite{HOR}. 
(For $f:K\to\complex$ a general function not in $D(K)$, we set 
$\|f\|_D=\infty$.) 

\proclaim{Lemma 1.8} 
Let $f:K\to \complex$ be a bounded function, $n$ a positive integer, and 
$(\delta_i)$ a sequence of positive numbers of length at least $n$. 
Then if $\os_n(f,(\delta_i))\ne\emptyset$, 
$$\|f\|_D \ge \sum_{i=1}^n \delta_i + \| f\mid \os_n (f,(\delta_i)\|_\infty\ .
\leqno(17)$$ 
\endproclaim 

\demo{Remark} 
We shall show below (in ......) that for $f$ a simple real-valued 
$D$-function, or more generally, for $f\in SD$, that 
the above estimate is exact. That is, 

\iitem{} {\it $\|f\|_D$ equals the sup of the right-hand side in the 
inequality $(17)$, over all $n$ and $(\delta_i)$ with $\os_n(f,(\delta_i))
\ne\emptyset$.} 
\enddemo 

\demo{Proof of Lemma 1.8} 
We may trivially assume $f\in D$, otherwise there is nothing to prove. 
Fix $(f_n)$ a sequence in $C_b(K)$ with $f_n\to f$ pointwise. 
We seek to estimate 
$$\tau \dfeq  \sup_{k\in K} |f_1| (k) 
+ \sum_{n=1}^\infty |f_{n+1} (k) - f_n(k)|\ . 
\leqno(18)$$ 
\enddemo 

The following tool easily yields Lemma 1.8. 

\proclaim{Sub-Lemma} 
Let $\U$ be an open  set in $K$ with $U\cap \os_n(f,(\delta_i)) \ne \emptyset$, 
and $0<\varep < \sum_{i=1}^n \delta_i$. 
There exist $m_1<m_2<\cdots <m_{2n}$ and $\V$ an open non-empty subset 
of $\U$ with 
$$\sum_{i=1}^n |f_{m_{2i}} - f_{m_{2i-1}}| > 
\sum_{i=1}^n \delta_i - \varep\ \text{ on }\ \V\ . 
\leqno(19)$$ 
\endproclaim 

\demo{Proof}
Let us first prove the Sub-Lemma, by induction on $n$. 
For $n=1$, let $\delta = \delta_1$; of course then $\os_1(f,(\delta_i)) 
= \os (f,\delta)$. 
We thus must find $\V$ and $i<j$ with $|f_j-f_i| >\delta-\varep$ on $\V$. 

First choose $u\in \U$ with $\osc f(u) \ge \delta$. 
Now set $\alpha= \varep/3$, and choose $k\in \U$ with 
$$\uosc f(k) >\delta-\alpha\ .
\leqno(20)$$ 
Then choose $i$ with $|f_i(k) - f(k)| <\alpha$. 
By continuity of $f_i$, choose $\W$ an open neighborhood of $k$ with 
$\W\subset \U$ and 
$$|f_i-f(k)| <\alpha \text{ on } \W\ .
\leqno(21)$$ 

Now using (20), we may choose $w\in \W$ with 
$$|f(w)- f(k)| >\delta-\alpha\ . 
\leqno(22)$$ 

Next, choose $j>i$ with $|f_j(w) - f(w)| <\alpha$. 
Again by continuity of $f_j$, we may choose $\V$ an open neighborhood of 
$w$ with $\V\subset \W$ and 
$$|f_j-f(w)| <\alpha \text{ on }\V\ .
\leqno(23)$$ 
Then evidently by (21)--(23), $|f_i-f_j| > \delta-3\alpha = \delta - \varep$ 
on $\V$. This  establishes the $n=1$ case. 
Now suppose the result proved for $n$, let $(\delta_i)$ a sequence 
of length at least $n+1$ be given, and set $Y = \os_n(f,(\delta_i))$. 
Then by definition, $\os_{n+1}(f,(\delta_i)) = \os (f\mid Y,\delta_{n+1})$. 
Now assuming $\U\cap Y \dfeq \tilde{\U}$ is non-empty, then since 
$\tilde\U$ is a relatively open subset of $Y$, by the $n=1$ case there exists 
$\tilde\W$ a non-empty relatively open subset of $Y$ with $\tilde\W \subset 
\tilde\U$ and $m_1<m_2$ with 
$$|f_{m_2}-f_{m_1}| > \delta - \varep/2\text{ on } \tilde\W\ .
\leqno(24)$$
Now let $\W = \{x\in K :|f_{m_2}(x) - f_{m_1}(x)| >\delta - \varep/2$. 
Then $\W$ is an open subset of $K$, and of course $\W\cap Y \supset \tilde\W$ 
is non-empty, so by the induction hypothesis applied to $(f_j)_{j>m_2}$, 
we may choose $m_3 < \cdots < m_{2(n+1)}$ with $m_3>m_2$ and $\V$ a 
non-empty open subset of $\W$ with 
$$\sum_{j=2}^{n+1} |f_{m_{2j}} - f_{m_{2j-1}}| 
> \sum_{i=1}^n \delta_i - {\varep\over2} \text{ on }\V\ . 
\leqno(25)$$ 
This completes the proof of the Sub-Lemma, for we have trivially that 
$|f_{m_2} - f_{m_1}|> \delta_{n+1} - \varep/2$ on $\V$ since 
$\V\subset\W$.\qed 
\enddemo  

To prove Lemma 1.8 itself, let $\varep>0$, set 
$\lambda=\|f\mid\os_n (f,(\delta_i))\|_\infty$ 
and choose $k\in \os_n(f,(\delta_i))$ with 
$$|f(k)| >\lambda -\varep\ . 
\leqno(26)$$ 
Next choose $m_0$ with $|f_{m_0} (k)|>\lambda-\varep$. 
Then let $\U= \{ x\in K: |f_{m_0} (k)| >\lambda-\varep\}$. 
$\U\cap \os_n(f,(\delta_i)) \ne\emptyset$, so by the Sub-Lemma (applied 
to $(f_j)_{j>m_0}$), we may choose $m_1<\cdots m_{2n}$ with $m_0<m_1$ 
and $v\in \U$ with 
$$\sum_{i=1}^n |(f_{m_{2i}} - f_{m_{2i-1}})(v)| > \sum_{i=1}^n \delta_i 
- \varep\ .$$ 
Hence since $v\in \U$, 
$$|f_{m_0} (v)| + \sum_{i=1}^n |f_{m_{2i}} - f_{m_{2i-1}}| (v) | 
> \lambda+ \sum_{i=1}^n \delta_i - 2\varep\ . 
\leqno(27)$$ 
But a collapsing series argument easily yields that 
$\tau \ge |f_{m_0} (v)| + \sum_{i=1}^n |f_{m_{2i}} - f_{m_{2i-1}}| (v)$. 
Hence $\tau \ge\lambda + \sum_{i=1}^n \delta_i -2\varep$; 
since $\varep>0$ is arbitrary, Lemma~1.8 is proved.\qed

\demo{Remark} 
For $f:K\to\complex$ a given function, set $\oosc f(x)= \olim_{y,z\to x} 
|f(y)-f(z)|$, for all $x\in K$. 
We term $\oosc f$ the {\it upper oscillation\/} of $f$. 
$\oosc f$ is usually defined as the oscillation of $f$; however our 
definition of $\osc f$ is  more appropriate for the study of $D(K)$. 
Now if $f$ is a real-valued function and $Lf$ is its lower semi-continuous 
``envelope,'' $Lf(x) \dfeq \ulim_{y\to x} f(y)$ for all $x\in K$, 
then we have that $\uosc f = \max \{Uf-f,f-Lf\}$ while 
$\oosc f= Uf-Lf$. 
$\oosc f$ is upper-semi-continuous but $\uosc f$ is not, in general. 
It is worth pointing out that for general $f:K\to \complex$, $\osc f\le 
\oosc f\le 2\uosc f$ and $\|\oosc f\|_\infty \le 2\|f\|_\infty$, while 
if $f$ is non-negative, then  $\|\osc f\|_\infty = \|f\|_\infty$. 
In \cite{HOR}, for a given $f$ and sequence $(\varep_i)$ of positive 
numbers, sets $K_n(f,(\varep_i))$  are defined inductively by letting 
$K_0(f,(\varep_i)) = K$ and $K_{n+1}(f,(\varep_i)) = \{x\in K_n :
\oosc f\mid K_n \ge \varep_{n+1}\}$. 
Then it follows easily that  $K_n(f,({\varep_i\over2})) \subset \os_n 
(f,(\varep_i)) \subset K_n(f,(\varep_i))$. 
Thus $K_n(f,(\delta_i))\ne \emptyset \To \sum \delta_i \le 
2\|f\|_D$ by Lemma~1.8.
\enddemo 

We next give a characterization of $B_{1/2}(K)$. This is implied by 
Proposition~2.3 of \cite{HOR}, except that we work here with arbitrary metric 
spaces. Our proof is somewhat different than the treatment in \cite{HOR}. 
We first define, for $\varep>0$, the (finite) Baire $\varep$-oscillation 
index of $f$, $i_B(f,\varep)$, as follows: 
$$i_B(f,\varep) = \sup \{n:\os_n (f,\varep) \ne \emptyset\}\ .$$ 
(We take the sup in $\{0,1,2,\ldots \} \cup \{\infty\}$.) 
Finally, we say that $f$ {\it is of finite Baire index provided there is 
an $n<\infty$ with $i_B(f,\varep) \le n$ for all $\varep>0$\/}; 
then we {\it define\/} $i_B(f)$, the Baire index of $f$, by 
$i_B(f) = \max_{\varep>0} i_B(f,\varep)$. 
(Thus $f$ is continuous iff $i_B(f)=0$.) 
It follows immediately from the above remark that $f$ is of finite 
Baire index if and only if for some $n$ and all $\varep>0$, 
$K_j (f,\varep) \ne\emptyset$ implies $j\le n$; in fact, as in \cite{HOR}, 
setting $\beta (f) =$ least $n$ with $K_n(f,\varep)=\emptyset$ for all 
$\varep>0$, then $i_B(f)+1= \beta (f)$. 
It is moreover easily seen (as shown in the proof of Proposition~1.10 below), 
that every simple $D$-function is of finite Baire index. 

\proclaim{Proposition 1.9} 
Let $f:K\to\complex$ be a bounded function. 
Then the following are equivalent. 

\iitem{\rm (a)} $f\in B_{1/2}(K)$.  

\iitem{\rm (b)} $i_B (f,\varep) <\infty$ for all $\varep>0$. 

\iitem{\rm (c)} $f$ is a uniform limit of simple $\D$-functions. 
\endproclaim 

\demo{Proof} 
(c) $\To$ (a) is trivial, since simple $\D$-functions belong to $D(K)$. 

(a) $\To$ (b).  This is a consequence of Lemma 1.8 and the following simple 
considerations. 

First, we note that for any functions $f$ and $g$, $\osc (f+g) \le \osc 
f+\osc g$, which implies that $|\osc f-\osc g| \le \osc (f-g)$. 
Now let $0<\eta <\varep/2$ and $f,g$ be such that $\|g-f\|_\infty \le\eta$. 
Then by the above, 
$$|\osc f- \osc g| \le 2\eta 
\leqno(28)$$ 

Now it follows easily that 
$$\os_n  (f,\varep) \subset \os_n (g,\varep-2\eta)\ \text{ for all }\ n\ . 
\leqno(29)$$ 
Indeed, for $n=1$, this simply says that $\osc f\ge \varep\To \os g\ge 
\varep -2\eta$, which is immediate from (28). 
Now if (29) is proved for $n$, let $L = \os_n (f,\varep)$ and 
$W= \os_n(g,\varep-2\eta)$. 
But then 
$$\eqalign{\os_{n+1}( f,\varep) 
& = \os (f|L,\varep) \cr 
&\subset \os (g|L,\varep-2\eta)\text{ by the $n=1$ case}\cr
&\subset \os (g|W,\varep -2\eta)\text{ since } L\subset W\cr
&= \os_{n+1} (g,\varep-2\eta)\ .\cr}$$ 

Now let $\varep>0$, and $f\in B_{1/2}(K)$. Choose $g\in D(K)$ with 
$$\|g-f\|_\infty \le {\varep\over3}\ . 
\leqno(30)$$ 
Now choose $n$ a positive integer with  
$$n{\varep\over 3} > \|g\|_D\ . 
\leqno(31)$$ 
It then follows from Lemma 1.8 that $\os_n(g,\varep/3) = \emptyset$. 
But then by (29), $\os_n(f,\varep)=\emptyset$; 
thus   $i_B(f,\varep) < n$, and (b) is proved. 
\enddemo

(b) $\To$ (c) follows easily from the following simple fact. 

\proclaim{Sublemma} 
Let $X$ be a metric space, $f:X\to\complex$ a bounded function, and 
$\varep >0$. Suppose $\osc f<\varep$ on $X$. 
Then given $\eta >0$, there is a simple $\D$-function $g$ on $X$ 
with $|g-f|<\varep+\eta$.
\endproclaim 

\demo{Proof} 
Since $f$ is bounded, we may choose $n$ and $c_1,\ldots,c_n$ distinct 
elements of $f(X)$ so that $\{c_1,\ldots,c_n\}$ is an $\eta$-net for $f(X)$. 
Now our hypothesis yields that given $x\in X$, there is an open 
neighborhood $\U$ of $x$ with $|f(y)-f(x)| <\varep$ for all $y\in \U$. 
But now choosing $i$ with $|c_i-f(x)|<\eta$, we have that 
$|f(y) - c_i| <\varep +\eta$ for $y\in \U$. 
Thus letting $W_i = \{x\in X: |f(x) -c_i| <\varep +\eta \}$ and 
$U_i$ be the interior of $W_i$, we have that $X= \cup U_i$. 
Now simply let $F_1=U_1$, $F_i = U_i\sim \bigcup_{j<i} U_j$ for 
$1< i\le n$. 
Then the $F_i$'s are in DCS, $X= \bigcup_{i=1}^n F_i$, so setting 
$g= \sum_{i=1}^n c_i \chix_{F_i}$, then $g$ is a simple $\D$-function 
with $|g-f| <\varep+\eta$.\qed 
\enddemo 

We now show (b) $\To$ (c). 
Let $\varep>0$, and $n= i_B (f,\varep)$. 
Thus by definition, $\os_{n+1} (f,\varep)=\emptyset$; we then have that if 
$X_i = \os_i (f,\varep) \sim \os_{i+1}(f,\varep)$ then 
$\osc f\mid X_i <\varep$ for all $0\le i\le n$. 
Thus by the Sub-Lemma, we may choose for each $i$, a simple $\D$-function 
$g_i$ on $X_i$ with $|g_i-f| <2\varep$ on $X_i$. 
Then letting $g= \sum_{i=0}^n g_i\chix_{X_i}$, $g$ is a simple 
$D$-function on $K$, with $|g-f|<2\varep$, so $\|g-f\|_\infty \le 2\varep$.\qed 
\medskip

The next result completes the proof of Proposition 1.4. 

\proclaim{Proposition 1.10} 
Let $f$ be a simple function in $B_{1/2}(K)$. 
Then $f$ is a simple $\D$-function.
\endproclaim 

\demo{Proof} 
Let $f$ be non-constant, and $c_1,\ldots,c_n$ be the distinct values of $f$ 
(so ${n\ge2}$). 
Let $\varep = \min \{|c_i-c_j| : i\ne j\}$. 
Then it follows that if $W\subset K$, $w\in W$, and $\osc f|_W (w) <\varep$, 
then $f$ is constant on a relative neighborhood of $w$ (in $W$). 
Thus if $\osc f|_W <\varep$, $f$ is continuous on $W$; if 
$d_1,\ldots,d_k$ are the distinct values of $f$ on $W$, then setting $W_i=
\{x\in W: f(x) = d_i\}$, then the $W_i$'s are relatively open, and hence 
relatively clopen subsets of $W$; of course then $f|_W = \sum d_i\chix_{W_i}$. 
If $W$ is itself a DCS in $K$, then we have that the $W_i$'s themselves 
are DCS's in $K$. 
\enddemo 

Now since $f\in B_{1/2}(K)$, $n\dfeq i_B(f,\varep) <\infty$. 
Setting $X_i = \os_i (f,\varep) \sim \os_{i+1} (f,\varep)$ for $0\le i\le n$, 
and fixing $i$, then $\osc f|X_i <\varep$ by definition. 
The above observation thus yields that $f|X_i$ is a $\D(X_i)$-function. 
Since $X_i$ itself is a DCS  in $K$, $f=\sum_{i=0}^n f\chix_{X_i}$ is 
a simple $\D$-function. 

\demo{Remark} 
The above argument yields a natural method for computing $i_B(f)$ for $f$ a 
simple $\D$-function.  
Define sets $K_0 = K\supset K_1 \supset K_2 \cdots$ inductively by letting
$K_{j+1} = \{x\in K_j :f|K_j$ is discontinuous at $x\}$. 
Then the above argument yields that if $\varep$ is as defined at the beginning 
of the proof, then $K_j = \os_j(f,\varep) = \os_j (f,\varep')$ for all 
$j$ and $0<\varep'\le \varep$; hence $i_B(f)$ is the largest $n$ with 
$K_n\ne \emptyset$. 
The argument also shows directly that if $W\subset K$ and $\chix_W$ is a 
$D$-function, then $W$ is a finite disjoint union of DCS's. 
\enddemo 

In the sequel, we shall find it convenient to introduce the following 
semi-norm on $D(K)$. 

\demo{Definition} 
Define $\|\cdot\|_{qD} = \|\cdot \|_{qD(K)}$ on $D(K)$ by 
$$\|f\|_{qD} = \inf_{\varphi \in C_b(K)} \|f-\varphi\|_{D(K)}\ ,\qquad 
\text{all }\ f\in D(K)\ .$$ 
For course this is really the quotient norm on $D(K)/C_b(K)$; that is, 
letting $\pi :D(K) \to D(K)/C_b(K)$ be the canonical map, then 
$\|\pi f\| = \|f\|_{qD}$. 
It is easily seen that for all $f\in D(K)$ 
$$\eqalign{
\|f\|_{qD} &= \inf \sup_{k\in K} \sum_{n=1}^\infty |(f_{n+1} - f_n)(k)| \ , \cr
&\qquad \text{ the infimum over $(f_n)$ in $C_b(K)$ with $f_n\to f$ 
pointwise.}\cr}
\leqno(32)$$

Several of our preceding results can also be formulated in terms of 
$\|\cdot\|_{qD}$. For example, the proof of Lemma~1.8 easily yields that 
$$\os_n(f,(\delta_i)) \ne\emptyset \to \|f\|_{qD} \ge \sum_{i=1}^n 
\delta_i\ . 
\leqno(33)$$ 
\enddemo 

The following result shows the simple connection between $\|\cdot\|_D $ 
and $\|\cdot\|_{qD}$. 

\proclaim{Proposition 1.11} 
For any $f\in D(K)$, $\|f\|_{qD} \le \|f\|_D \le \|f\|_\infty +
\|f\|_{qD}$. 
\endproclaim 

The first inequality is trivial. To prove the second one, we note the 
following elementary result, whose proof is left to the reader. 

\proclaim{Lemma 1.12} 
Fix $\lambda >0$, and for $z$ a complex number, let 
$$\cases \tilde z = \lambda {z\over |z|}&\text{if $|z|>\lambda$}\cr 
\tilde z= z&\text{if $|z| \le \lambda$}.\endcases
\leqno(34)$$ 
Then $z\to \tilde z$ is Lipschitz with constant one; that is, 
$$|\tilde z-\tilde w| \le |z-w|\ \text{ for all complex numbers }\ z,w\ .
\leqno(35)$$ 
\endproclaim 

\demo{Proof of Proposition 1.11} 
Using the formulation (32), let $\varep>0$ and choose $(f_n)$ in $C_b(K)$ 
with  $f_n\to f$ pointwise and
$$\sum_{n=1}^\infty |f_{n+1} -f_n| < \|f\|_{qD} +\varep\ .
\leqno(36)$$ 
Now let $\lambda = \|f\|_\infty$ and assume without loss of generality 
that $\lambda>0$. 
But then also $\tilde f_n\to f$ pointwise 
(where $(\tilde f_n)(k)= \widetilde{f_n(k)}$, as defined in (34)). 

Lemma 1.12 shows that $\tilde f_n$ is continuous for all $n$, and of course 
$$\eqalign{
\|f\|_D & \le \| \, |\tilde f_1| + \sum |\tilde f_{n+1}-\tilde f_n|\, 
\|_\infty\cr 
&\le \lambda + \|\sum |f_{n+1}-f_n|\, \|_\infty \quad\text{(by Lemma 1.12)}\cr 
&\le \lambda + \|f\|_{qD} +\varep\cr}$$ 
Since $\varep >0$ is arbitrary, 1.12 is proved.\qed
\enddemo 

The proof of 1.11 immediately yields the following. 

\proclaim{Corollary 1.13} 
Let $f$ be in $D(K)$, and $\varep>0$ be given. There exists a sequence 
$(f_n)$ in $C_b(K)$ with $f_n\to f$ pointwise, 
$|f_1| + \sum |f_{n+1}-f_n| < \|f\|_D +\varep$, and 
$\|f_n\|_\infty \le \|f\|_\infty$ for all $n$. 
\endproclaim 

We conclude this section with several applications of partitions of 
unity; in particular, we establish the following {\it Localization 
Principle\/}. 

\proclaim{Proposition 1.14}
Let $\lambda >0$, and $f:K\to\complex$ be a given bounded function. 

\iitem{\rm (a)} If for all $x\in K$, there is an open neighborhood $U_x$ 
of $x$ with $\|f\mid U_x\|_{D(U_x)} \le\lambda$, then $f\in D(K)$ and 
$\|f\|_{D(K)}\le\lambda$. 

\iitem{\rm (b)} If for all $x\in K$, there is an open neighborhood $U_x$ 
of $X$ with  $\|f\mid U_x\|_{qD(U_x)} \le\lambda$, then $\|f\|_{qD(K)}
\le\lambda$.
\endproclaim 

\demo{Remark} 
We give an alternate proof of the first part in Section~3, over the real 
scalars, using the transfinite oscillations rather than partitions of unity. 
\enddemo 

We first recall the basic results concerning the existence of partitions 
of unity; of course all results are valid for paracompact spaces in general. 

\demo{Definition}
An open cover $\U$ of $K$ is called {\it locally finite\/} if every point in 
$K$ is contained in some open set meeting only finitely many members of $\U$. 
\enddemo 

\demo{Definition} 
{\it Given $\V$ a locally finite open cover of $K$, a family 
${\P} = \{p_v :v\in \V\}$ of continuous functions on $K$ is called a 
partition of unity fitting $\V$ if for all $v\in \V$,} 
 
\iitem{(a)} $0\le p_v\le1$
\iitem{(b)} $\supp p_v \dfeq \{ k\in K:p_v (k)\ne0\} \subset v$
\iitem{(c)} $\sum_{v\in \V} p_v =1$.
\enddemo 

Now our needed topological-analytical result may be formulated as 
follows (cf. \cite{K}). 

\proclaim{Lemma 1.15. (The Partition of Unity Lemma)} 
\iitem{\rm (a)} Every open cover of a metric space has a locally finite 
refinement.
\iitem{\rm (b)} For every locally finite open cover $\V$ of a metric space, 
there exists a partition of unity ${\P}$ fitting $\V$.
\endproclaim 

It is convenient to isolate the next simple principle, which will be used 
several times in the sequel. 

\proclaim{Lemma 1.16} 
Let $\V$ be an open cover of $K$, and ${\P}$ a partition of unity  fitting 
$\V$. Suppose for each $v\in V$, there is given $\varphi_v\in C_b(V)$. 
Then $\varphi \dfeq \sum_{v\in \V} (\varphi_v \chix_v)p_v$ is a 
continuous function on $K$. If moreover for some $\lambda<\infty$, 
$\|\varphi_v\|_\infty\le\lambda$ for all $v$, then also 
$\|\varphi\|_\infty \le\lambda$.
\endproclaim 

\demo{Proof} 
We first note that fixing $v\in\V$, then $(\varphi_v\chix_v)p_v$ is 
continuous on $K$. Indeed, it is trivially continuous on $v$. 
But if $x\notin v$, $(x_n)$ in $K$, and $x_n\to x$, then since 
$p_v$ is continuous on $K$ and $p_v(x)=0$, $p_v(x_n)\to0$, so 
$(\varphi_v\chix_v)(x_n)p_v(x_n)\to 0$ since $\varphi_v$ is bounded. 

Next, given $x\in K$, choose $F$ a finite non-empty subset of $\P$ and $U$ 
an open neighborhood of $x$ with $U\cap V=\emptyset$ all $v\notin F$. 
But then by property (b) of the above definition, 
$\varphi\mid U = \sum_{v\in F} (\varphi_j \chix_v)p_v\mid U$, which is 
of course continuous on $U$. 
Finally, letting $\|\varphi_v\|_\infty \le\lambda$ for all $v\in\V$, 
then for $k\in K$, $|\varphi(k)|\le \sum_{v\in \V} 
\|\varphi_j\|_\infty p_v(k) \le\lambda \sum p_v (k)=\lambda$.\qed
\enddemo 

We now pass to the proof of the Localization Principle. 

\demo{Proof of 1.14}

(a) Let $\U \dfeq \{U_x:x\in K\}$ be as in the statement, and let $\V$ 
be an open locally finite refinement of $\U$; 
then let $\P$ be a partition of unity fitting $\V$. 
Let $\varep >0$. 
It follows that given $v\in\V$, we may choose $(\varphi_j^v)_{j=1}^\infty$ 
in $C_b(v)$ so that 
$$\left\{\eqalign{
&\sum_j |\varphi_j^v| < \|f\mid v\|_{D(v)} +\varep \le \lambda+\varep\cr 
&\text{and }\ f= \sum_j \varphi_j^v \text{ pointwise on } v\ .\cr}
\right.\leqno(37)$$ 
Now fix $j$ and define $\varphi_j$ on $K$ by 
$$\varphi_j = \sum_{v\in \V} (\varphi_j^v \chix_v) p_v\ . 
\leqno(38)$$
By Lemma 1.16, $\varphi_j \in C_b(K)$; we then have 
$$\eqalign{\sum_j |\varphi_j| & \le \sum_{v\in\V} \sum_j 
|\varphi_j^v | \chix_v p_v\cr
&\le \sum_{v\in \V} (\lambda+\varep) p_v = \lambda+\varep\text{ by (32)}\cr} 
\leqno(39)$$ 
Again using (37), 
$$\sum_j \varphi_j = \sum_{v\in\V} \sum_j \varphi_j^v \chix_v p_v
= \sum_{v\in\V} fp_v =f\ . $$ 
Thus we have proved $f\in D(K)$ with $\|f\|_D \le \lambda+\varep$. 
Since $\varep>0$ is arbitrary, (a) is proved. 

To prove (b), this time let $\varep >0$, fix $v\in\V$ and choose $g_v 
\in C_b (v)$ and $(\varphi_j^v)$ in $C_b(v)$ with 
$$\sum|\varphi_j^v| <\lambda +\varep \ \text{ and }\ 
f= g_v + \sum \varphi_j^v\text{ on } v\ .
\leqno(40)$$ 
But then $\|g_v\|_\infty \le\lambda +\|f\|_\infty +\varep$, 
so $g_v$ is bounded. Finally, define $g=\sum_{v\in \V} (g_v\chix_v)p_v$ 
and $\varphi_j$ by (38) for all $j$. 
By Lemma~1.16, $g$ is continuous and again, (39) holds. 
Finally, we check that 
$$ f= g+\sum \varphi_j\ . 
\leqno(41)$$ 
Indeed, 
$$\eqalign{ g+\sum \varphi_j & = \sum_{v\in\V} \sum_j (g_v+\varphi_j^v)
\chix_v p_v\cr 
&= \sum_{v\in\V} f\cdot p_v\ \text{ by (40)}\cr 
&= f\ .\cr}$$ 
Thus $\|f\|_{qD} \le \lambda+\varep$; since $\varep >0$ is arbitrary, 
$\|f\|_{qD} \le\lambda$.\qed 
\enddemo 

The next result shows that functions of small oscillation are close to 
continuous functions. 
The qualitative result is of course standard; however our quantitative version 
is essentially immediate from Lemmas~1.15 and 1.16 and also holds for 
the complex scalars. 

\proclaim{Proposition 1.17} 
Let $f:K\to\complex$ and $\varep>0$ be given, and suppose $\uosc f(x)<\varep$ 
for all $x\in K$. 
Then there is a continuous $\varphi$ on $K$ with $|\varphi -f|<\varep$. 
\endproclaim 

\demo{Comment} 
Of course if $f$ is bounded, $\varphi$ is also, with $\|\varphi\|_\infty 
\le \|f\|_\infty +\varep$.
\enddemo 

\demo{Proof} 
Given $x\in X$, choose $U_x$ an open neighborhood of $x$ so that 
$$|f(y)-f(x)| <\varep\text{ for all } y\in \U_x\ . 
\leqno(42)$$ 
Let $\V$ be a locally finite open refinement of $\U \dfeq \{U_x :x\in X\}$ 
and $\P = \{p_v :v\in\V\}$ a partition of unity fitting $\V$. 
For each $v\in\V$, choose $x$ with $v\subset U_x$; then set 
$\lambda_v = f(x)$. 

Now simply set $\varphi = \sum_{v\in\V} \lambda_v p_v$. 
By Lemma 1.16, $\varphi$ is continuous. 
Since $f=\sum f\cdot p_v$, we have that 
$$|\varphi-f| = |\sum_{v\in\V} (\lambda_v -f)p_v| 
\le \sum_{v\in\V} |\lambda_v -f| p_v <\varep\text{ by (42).} 
\eqno\qed$$ 
\enddemo 

We conclude this section with a proof of Baire's famous theorem. 

\proclaim{Proposition 1.18} 
Let $f:K\to\real$ be a lower semi-continuous function. 
There exists a sequence of real-valued continuous functions $(\varphi_j)$ 
on $K$ with $\varphi_1 \le\varphi_2 \le\cdots\le \varphi_j \le 
\varphi_{j+1} \le\cdots$ and $\varphi_j \to f$ pointwise.
\endproclaim 

\demo{Remark} 
Of course if $f\ge0$, we easily obtain that  we may choose the $\varphi_j$'s 
$\ge0$ also, by simply setting $\widetilde\varphi_j = \max \{\varphi_j,0\}$ 
for all $j$. 
>From this it follows that if $f$ is a non-negative bounded lower 
semi-continuous function, then $\|f\|_D = \|f\|_\infty$. 
Finally, if $f$ is a bounded semi-continuous function, then we obtain that 
$\|f\|_D\le 3\|f\|_\infty$. 
Indeed, let $\lambda = \|f\|_\infty$ and assume without loss of 
generality that $f$ is lower semi-continuous. 
But then $f+\lambda\bone$ is non-negative lower semi-continuous, 
hence $\|f +\lambda\bone\|_D = \|f+\lambda\bone\|_\infty \le 2\lambda$, 
whence $\|f\|_D \le3\lambda$. 
\enddemo 

\demo{Proof of Proposition 1.18} 

It suffices instead  to construct a sequence $(\varphi_n)$ of continuous 
real-valued functions on $K$ so that $\varphi_n\le f$ for all $n$ and 
$f\equiv \sup \varphi_n$ (pointwise). 
We then simply let $\widetilde\varphi_n = \max \{\varphi_1,\ldots,
\varphi_n\}$ for all $n$. 

Now fix $n$ a positive integer, and let $x\in K$. 
By the lower semi-continuity of $f$ we may choose an open neighborhood 
$U_x$ of $x$, of diameter at most $1/n$, so that 
$$f(u) > f(x) - {1\over n}\ \text{ for all $u$ in }\ U_x\ . 
\leqno(43)$$ 

Let $\V_n$ be an open locally finite refinement of $\U \dfeq 
\{U_x :x\in K\}$, and let $(p_v^n)_{v\in\V_n}$ be a partition of unity 
fitting $\V_n$. 
Given $v\in \V_n$, choose $U_x$ with $v\subset U_x$; then 
set $\lambda_v^n = f(x) - \frac1n$. 
Now define $\varphi_n$ by 
$$\varphi_n = \sum_{v\in\V_n} \lambda_v^n p_v^n \ .
\leqno(44)$$
(In the above, ``$n$'' is an index, not a power!) 

Then by Lemma 1.16, $\varphi_n$ is continuous. 
Now we have, for $v$ in $\V_n$, that $\lambda_v^n p_v^n \le fp_v^n$. 
Indeed, on $v$, this is obvious by the definition of $\lambda_v^n$. 
But off $v$, both sides of the inequality are zero. 
Thus, 
$$\sum_{v\in\V_n} \lambda_v^n p_v^n \le \sum_{v\in\V_n} fp_v^n =f\ .$$ 
Finally, we verify that $f=\sup_n\varphi_n$. 
Fix $x\in K$ and $\varep >0$. 
Then choose $W$ an open ball centered at $x$ of radius $r$, with 
$f>f(x) -\varep$ on $W$
(i.e., if $\rho$ is the metric on $K$, $W= \{y:\rho (y,x)<r\}$). 
Now let $n$ be chosen with $\frac1n <r$, $\frac1n <\varep$. 
Suppose $p_v^n (x)\ne0$. 
Thus $x\in v$. Choose $U_y \in \U$ with $v\subset U_y$ and $\lambda_v^n = 
f(y)-\frac1n$. 
Now $\diam U_y\le \frac1n$. 
Hence $\rho (x,y) \le\frac1n$, so $y\in W$. 
Thus $\lambda_v^n >f(x) -\varep -\frac1n > f(x)-2\varep$. 
But then $\lambda_v^n p_v^n (x) \ge (f(x)-2\varep)p_v^n(x)$. 
Of course this holds trivially if $p_v^n(x)=0$ as well. 
Thus $\varphi_n(x) = \sum \lambda_v^n p_v^n (x) \ge 
\sum (f(x) -2\varep) p_v^n (x) = f(x) - 2\varep$. 
That is, $\sup \varphi_n \ge f-2\varep$. 
Since $\varep>0$ is arbitrary, the proof is complete.\qed
\enddemo 

\head \S2. The $D$-norm of the characteristic function of a set.\endhead 

In this section, we give a topological method for computing the 
$D$-norm of a $\{0,1\}$-valued function, and apply our result to show that 
the locally Lipschitz functions are precisely those which operate on $D(K)$ 
for general $K$. 
We also solve the problem of representing a set $W$ in $\D$ as a union 
of $k$ disjoint DCS's with $k$ optimal; it turns out that 
$k\sim \frac{n}2$ where $n = \|\chix_W\|_D$. 

As always, $K$ is a fixed metric space. 
For $A\subset K$, $\partial A$ denotes the {\it boundary\/} of $A$. 
If $L\subset K$, then $\partial_L A$ denotes the {\it boundary\/} of 
$A\cap L$, {\it relative\/} to $L$. 
Thus, $x\in\partial_L A$ iff $x\in L$ and there exist sequences $(x_n)$ 
and $(y_n)$ converging to $x$ with $x_n$ in $L\sim A$ and $y_n$ 
in $A\cap L$, for all $n$. 

We then define $\partial^n A$, the $n^{th}$ boundary of $A$, as follows: 
$\partial^0 A= K$. 

If $\partial^n A$ has been defined, then 
$\partial^{n+1} A= \partial_L A$ where $L= \partial^n A$. 
It follows immediately by induction that $\partial^nA$ is closed for all $n$. 
As we shall see shortly, $A\in \D$ iff $\partial^n A=\emptyset$ for some $n$. 

\demo{Remark} 
The definition is easily extended to all ordinal numbers $\alpha$, 
rather than just the finite ones. 
Thus, if $\beta$ is a limit ordinal, set $\partial^\beta A = \bigcap_{\alpha 
<\beta} \partial^\alpha A$. 
If $\beta =\alpha+1$, set $\partial^\beta A= \partial_L A$ where 
$L=\partial^\alpha A$. 
Again $\partial^\beta A$ is closed for all $\beta$. 
If $K$ is a Polish space we then have that $\chix_A$ is a Baire-1 
function if and only if $\partial^\alpha A=\emptyset$ for some 
$\alpha <\omega_1$. 
Since $A\cap (\partial^\beta A\sim \partial^{\beta+1}A)$ is a 
clopen subset of $\partial^\beta A\sim \partial^{\beta+1}A$ for all $\beta$, 
it follows that $\bigcup_{\beta<\alpha} A\cap (\partial^\beta A\sim 
\partial^{\beta+1} A)$ represents $A$ as an $F_\sigma$ whose complement is 
also on $F_\sigma$; i.e., we obtain the classical fact (cf. \cite{H}) 
that $A$ is then both $F_\sigma$ and $G_\delta$. 
\enddemo 

We now have the following simple result. 

\proclaim{Proposition 2.1} 
Let $A\subset K$. Then 

\iitem{\rm (a)} $\osc\chix_A =\uosc\chix_A =\oosc\chix_A= \chix_{\partial A}$. 
\iitem{\rm (b)} $\partial^n A= \os_n (\chix_A,\varep)$ for all positive 
integers $n$, and $0<\varep \le1$. 
\endproclaim 

\demo{Proof} 

(a) If $x\notin \partial A$, then $\chix_A$ is continuous at $x$, so the 
various oscillations given in (a) are all zero at $x$. 
Now it is also evident that $0\le \oosc \chix_A \le1$. 
Suppose $x\in \partial A$, and choose $(x_n),(y_n)$ converging to $x$ 
with $(x_n)\subset\sim A$, $(y_n) \subset A$. 
If $x$ is in $A$, then $\uosc \chix_A (x) \ge \lim_{n\to\infty} 
|\chix_A(x_n) - \chix_A(x) | =1$. 
If $x$ is not in $A$, then $\uosc \chix_A (x) \ge \lim_{n\to\infty} 
|\chix_A (y_n) - \chix_A(x)| =1$. 
thus $\uosc \chix_A(x) = 1\le \osc \chix_A(x) \le \oosc \chix_A(x)\le1$. 
This proves (a). 
Thus we obtain that $\os (\chix_A,\varep) =\partial A$ for any 
$0<\varep\le1$. Then (b) follows easily by induction.  
Indeed, we have just seen its validity for $n=1$. 
Suppose proved for $n$. 
But then $\os_{n+1}(\chix_A,\varep) = \os (\chix_A |L,\varep)$ where 
$L=\partial^n A$; but this equals $\partial_LA$, again by 
the $n=1$ case.\qed 
\enddemo 

We now define the (finite) index, $i(A)$, of $A\subset K$ as follows: 

\demo{Definition} 
$i(A)$ equals the largest $n$ with $\partial^nA\ne\emptyset$, if there 
is such an $n$; otherwise $i(A) =\infty$. 

Of course we say that $A$ is a {\it set of finite index\/} if $i(A)<\infty$. 
Evidently Proposition~2.1 shows that $i(A) = i_B(\chix_A,\varep)$ for 
all $0<\varep\le1$. Thus $i(A) = i_B(\chix_A)$. 
\enddemo 

We may now formulate the main structural result of this section. 
Besides giving a formula for the exact computation of $D$-norm of $\chix_A$, 
for any given set in $\D$, we also obtain the rather surprising 
result that there exist $k$ and disjoint DCS's 
$W_1,\ldots,W_k$ with $\|\chix_A\|_D = \sum_{i=1}^k \|\chix_{W_i}\|_D$. 

\proclaim{Theorem 2.2} 
$A$ belongs to $\D$ if and only if $A$ is of finite index. 
Suppose this is the case, and let $n=i(A)$. 

\item{\rm (a)} $\|\chix_A\|_{qD} = n$.  
\item{\rm (b)} {\rm (i)} If $A\cap \partial^n A=\emptyset$, then  
$\|\chix_A\|_D =n$. 
\item{} {\rm (ii)} If $A\cap \partial^n A\ne\emptyset$, then $\|\chix_A\|_D
= n+1$. 
\item{\rm (c)} If $A\cap \partial^n A=\emptyset$, then letting $k= 
[\frac{n+1}2]$, $A$ is a union of $k$ disjoint {\rm DCS's}; moreover one 
of these sets may be chosen open in case $n$ is odd. 
\item{\rm (d)} If $A\cap \partial^n A \ne\emptyset$, then letting 
$k= [\frac{n}2 +1]$, $A$ is a union of $k$ disjoint {\rm DCS's}; 
moreover one of these may be chosen open in case $n$ is even. 
\endproclaim 

The proofs of (c), (d) are constructive, and also yield the optimal number of 
disjoint DCS's where union is $A$. 

\proclaim{Corollary 2.3} 
Let $A$, $n$ be as above, and let $k$ be as in case (c) or (d). 
Suppose $W_1,\ldots,W_m$ are disjoint {\rm DCS's} with $A=\bigcup_{i=1}^m 
W_i$. Then $m\ge k$. 
\endproclaim 

\demo{Proof} 
We have that 
$$\|\chix_A\|_D \le \sum_{i=1}^m \|\chix_{W_i}\|_D \le 2m\ . 
\leqno(1)$$ 
Thus in case (c), 
$n\le 2m$ by (b)(i), so $m\ge \frac{n}2$, hence $m\ge [\frac{n+1}2]$. 
The proof for case (d) is identical, using (b)(ii).\qed 
\enddemo 

Let us next dispose of the parts of Theorem 2.2 that follow immediately 
from our previous results. 
To see the first assertion, suppose $A\in \D$. 
Then $\chix_A\in D$, hence by Proposition~1.9 and the previous result, 
$A$ is of finite index. 
But if $n=i(A)<\infty$, then $A=\bigcup_{i=0}^n A\cap (\partial^iA\sim 
\partial^{i+1}A)$, and for each $i$, we have that setting $L= \partial^i A 
\sim \partial^{i+1}A$, then $A\cap L$ is a relatively clopen subset of $L$, 
hence $A\cap L$ is a DCS, so $A\in\D$. 
Next, as noted in Section~1, 
if $\os_n(f,\varep)\ne\emptyset$, then the proof of Lemma~1.8 easily yields 
that $\|f\|_{qD} \ge n\varep$. 
Hence since $\os_n(\chix_A,1) = \partial^n A$ by Proposition~2.1, 
$\|\chix_A \|_{qD} \ge n$. 
We delay the proof of the reverse inequality until 
Section~3, where we shall see that the natural properties of transfinite 
oscillations render this transparent. 
Evidently we have immediately that $\|\chix_A\|_D\ge n$ in case (b)(i), 
and in case (b)(ii), Lemma~1.8 gives 
that if $f=\chix_A$, $\|f\|_D \ge n+ \|f|\os_n(f,1)\|_\infty 
= n+ \|f|\partial^n A\|_\infty = n+1$. 

Finally, $\|\chix_A\|_D\le n$ in (b)(i), $\|\chix_A \|_D\le n+1$ in (b)(ii), 
follows immediately from (c) and (d). 
(2.2(b)may also be deduced, alternatively, from the general structural 
result in Section~3.) 
Indeed, let $A= \bigcup_{i=1}^k W_i$ with the $W_i$'s disjoint DCS's, 
$k$ as in (c) or (d). 
Now in case (c), if $n$ is even, then $\|\chix_A\|_D \le 
\sum_{i=1}^k \|\chix_{W_i}\|_D \le 2k=n$. 
If $n$ is odd, then one of the $W_i$'s may be chosen open, say $W_1$ is open. 
But then 
$$\|\chix_A\|_D \le 1+ 2(k-1) = 1+ 2\Bigl({n+1\over2}\Bigr) - 2=n\ .$$ 
Similarly, in case (d), we again obtain that 
$\|\chix_A\|_D \le n+1$. 
Thus it remains to construct the representations in (c) and (d). 
We first require two simple results. 

\proclaim{Lemma 2.4} 
Let $A\subset K$ and $i(A) \le1$. 
Then $A$ is a {\rm DCS}.
\endproclaim 

\demo{Proof} 
If $i(A)=0$, $A$ is a clopen set, since then $\partial A=\emptyset$, 
so this is trivial. 
Suppose $i(A) =1$. Now if $A\cap \partial A=\emptyset$, then $A$ is an open 
set, so suppose finally that $A\cap \partial A\ne\emptyset$. 
Since $\partial^2 A=\emptyset$, $A\cap \partial A$ is a relatively clopen 
subset of $\partial A$, so there exists an open set $\V$ with $\V\cap 
\partial A= A\cap \partial A$. 
Let $\U= \V\cup \Int A$. 
We claim that $A= \bar A\cap \U$, (hence $A$ is a DCS). 
First, if $x\in A$ and $x\notin \Int A$, then $x\in \partial A$, 
hence $x\in A\cap \partial A\subset \U$, so $x\in \bar A\cap \U$. 
But if $x\in \bar A\cap \U$ and $x \notin \Int A$, then 
$x\in \V$, $x\in \bar A \To x\in \partial A 
\To x\in \V\cap \partial A\subset A$.\qed 
\enddemo 

\demo{Remark} 
It is easily seen that if $W$ is a DCS, then $i(W) \le 2$. 
Indeed, choose $A\supset B$ closed with $B$ nowhere dense in $A$, and 
$W= A\sim B$. 
Then $B\subset \partial W \subset A$. 
But then $W\cap \partial W$ is a relatively open subset of $\partial W$, hence
its index {\it relative\/} to $\partial W$ is at most 1. 
In fact, $i(W)=2$ if and only if $W\cap \partial W$ is relatively open 
in $\partial W$ but not a closed set. 
\enddemo 

Given  sets $A$ and $B$, and $j$ an integer, let 
$\partial^j (A|B)$ denote  the $j^{th}$ boundary of $A\cap B$, relative 
to $B$. (Thus $\partial (A|B) = \partial_B(A)$.) 

\proclaim{Lemma 2.5}  
\iitem{\rm (a)} Let $B,\U$ be given sets with $\U$ open in $K$, and 
$j\ge0$. Then $\partial^j (B|\U) = (\partial^j B)\cap \U$. 
\iitem{\rm (b)} Let $A$ be any set, $i$ any integer. Then 
$A\cap (\partial^i A\sim \partial^{i+2}A)$ is a {\rm DCS}. 
\endproclaim 

\demo{Proof} 

(a) The statement is evident for $j=1$. 
In fact, we have for any set $M$ that $\partial_{M\cap \U} B = 
(\partial_M B) \cap \U$. 
But then if the result is proved for $j$ and $M= \partial^j B$, we have 
that $(\partial^{j+1}B)\cap \U = \partial_M B\cap \U$ and 
$\partial^{j+1} (B|\U) = \partial_{M\cap \U} B$ by induction hypothesis. 

(b) Let $X= \partial^i A$, $\U = \partial^i A\sim \partial^{i+2}A$, 
$B=A\cap X$. 
Now we compute indices relative to $X$. 
Evidently, $\U$ is open in $X$. 
By Lemma~2.4, it suffices to prove that $i(B\cap \U|\U)$ is at most one. 
For then $B\cap \U$ is a relative DCS in $\U$, so it's a DCS in $K$. 
Now by (a), $\partial_X^2 (B|\U) = (\partial_X^2 B)\cap \U 
= (\partial^{i+2} A)\cap \U=\emptyset$, proving the assertion.\qed 

We now complete the proof of Theorem 2.2, 
(except for the equaltiy in (a)) proving (c) and (d). 
Suppose then $A\cap \partial^n A=\emptyset$, and first assume $n$ is even, 
so $k=\frac{n}2$. 
Then setting $W_i = A\cap (\partial^{2i} A\sim \partial^{2(i+1)}A)$, 
we have that $W_i$ is a DCS by Lemma~2.5(b), and of course 
$A=\bigcup_{i=0}^{k-1} W_i$, so this case is proved. 
If $n$ is odd, then $n=2k-1$ ($k$ as in (c)). 
This time set $W_i = A\cap (\partial^{2i-1}A\sim \partial^{2i+1}A)$ for 
$1\le i<k$ and $W_0 = A\cap \sim \partial^1 A$. 
Then $W_0$ is open, $W_i$ is a DCS for all $1\le i<k$ as before, and 
$A= \bigcup_{i=0}^{k-1} W_i$. 
Thus (c) is proved. 
Finally, for case (d), suppose first $n$ is even. 
Thus $n=2k-2$ ($k$ as in (d)). 
Again set $W_i = A\cap (\partial^{2i-1} A\sim \partial^{2i+1}A)$ for 
$1\le i<k$ and $W_0 = A\cap \sim \partial^1 A$. 
Thus $W_0$ is open, the $W_i$'s are DCS's, and 
$A= \bigcup_{i=0}^{k-1} W_i$ as before. 
(Note that $\partial^{2k-1}A=\emptyset$ since $i(A) = 2k-2$.) 
Finally, suppose $n$ is odd. 
Then $n=2k-1$, and thus we simply set $W_i = A\cap (\partial^{2i} A
\sim \partial^{2i+2}A)$ for $0\le i<k$. 
Again using that $\partial^{2k}A=\emptyset$, 
$A=\bigcup_{i=0}^k W_i$, completing the proof.\qed 
\enddemo 

We next give examples of the phenomena described in the Theorem. 
For $K$ a metric space, $\alpha$ an ordinal, $K^{(\alpha)}$ denotes 
the $\alpha^{th}$ derived set of $K$. 
Thus, $K^{(0)}=K$ by definition, and $K^{(1)}$ denotes the set 
of cluster points of $K$. 
If $K^{(\alpha)}$ has been defined for all $\alpha<\beta$ and $\beta$ is 
a successor, say $\beta = \alpha+1$, then $K^{(\beta)} = (K^{(\alpha)})^{(1)}$. 
Otherwise, $K^{(\beta)} = \bigcap_{\alpha<\beta} K^{(\alpha)}$.  

\proclaim{Proposition 2.6} 
Let $K$ be a metric space such that $K^{(n)} \ne \emptyset$ for all 
$n=1,2,\ldots$. 
Then for every positive integer $n$, there exist sets $A,B\subset K$, with 
$i(A)=n$, $i(B)=n$, $\|\chix_A\|_D=n$ and $\|\chix_B\|_D= n+1$. 
\endproclaim 

\demo{Proof} 
The hypotheses imply that we may choose closed subsets 
$K= X^0 \supset X^1 \supset X^2 \supset \cdots \supset X^n \supset X^{n+1} 
\cdots$ so that 
$$X^{n+1}\ \text{ is nowhere dense in $X^n$ for all $n$.} 
\leqno(2)$$ 
Indeed, if $K$ has no perfect subset, we may simply set $X^n = K^{(n)}$ 
for all $n$. 
Otherwise, we may apply the following  topological fact: 
{\it if $X$ is a perfect metric space (i.e., $X^{(1)} = X\ne\emptyset$), 
there exists a closed perfect nowhere-dense subset $Y$ of $X$.} 
We may then simply choose $Y$ a closed perfect subset, and choose 
$Y = Y^0 \supset\cdots\supset Y^n \supset \cdots$ with 
$Y^{n+1}$ nowhere-dense in $Y^n$ for all $n$; then set 
$X^n = Y^n$ for $n\ge1$. 

Now fix $n$, and suppose first that $n$ is odd, say $n=2k+1$. 
Let $A= \bigcup_{j=0}^k X^{2j} \sim X^{2j+1}$. 
We then claim: 
$$\left\{ \eqalign{
&\partial^i A = X^i\text{ for all } 0\le i\le n\ ,\text{ and }\cr
&\partial^{n+1}A =\emptyset\ .\cr}\right.
\leqno(3)$$ 
Suppose this is proved for $0\le i\le n$. 
Now if $i<n$, then if $i$ is even, say $i=2j$, $A\cap \partial^i A\supset 
X^{2j} \sim X^{2j+1}$, and since then $\overline{X^{2j}\sim X^{2j+1}} 
\supset X^{i+1}\sim X^{i+2}$, and $A\cap (X^{i+1}\sim X^{i+2})=\emptyset$, 
$\partial^{i+1}A\supset \overline{X^{i+1}\sim X^{i+2}} = X^{i+1}$, 
we clearly have that $\partial^{i+1} A\supset X^{i+1}$.  
But $X^{2j} \sim X^{2j+1}$ is a relatively open subset of $\partial^i A$, 
so $\partial^{i+1} A = X^{i+1}$. 
On the other hand, if $i$ is odd, say $i=2j-1$, then $A\cap \partial^i A
\cap (X^{2j-1}\sim X^{2j})=\emptyset$. 
But since $X^{2j-1}\sim X^{2j}$ is a relatively open dense subset of 
$\partial^i A$, $\partial^{i+1}A \subset X^{2j}$. 
Since $X^{2j}\sim X^{2j+1} \subset A\cap \partial^i A$ and 
$X^{2j}\sim X^{2j+1}\subset \overline{X^{2j-1} \sim X^{2j}}$, 
$\partial^{i+1}A \supset X^{2j}\sim X^{2j+1}$. 
Thus $\partial^{i+1} A\supset \overline{X^{2j}\sim X^{2j+1}} = X^{2j}$. 
Finally, if $i=n$, we obtain that $A\cap \partial^n A=\emptyset$, proving (3). 
Hence we have that $i(A) = n$, and by Theorem~2.2, 
that $\|\chix_A\|_D= n$. 

Next, let $B= (\bigcup_{j=1}^k X^{2j-1} \sim X^{2j}) \cup X^{2k+1}$. 
Again, by a proof about identical to the above, we have 
$$\partial^i B= X^i\text{ for all } 0\le i\le n\text{ and } \partial^{n+1} 
B = \emptyset\ . 
\leqno(4)$$
Thus again $i(B)=n$.
Evidently now $B\cap \partial^n B= X^n \ne \emptyset$, so 
$\|\chix_B \|_D = n+1$ by Theorem~2.2. 

Finally, if $n$ is even, say $n=2k$, let $A= \bigcup_{j=1}^k X^{2j-1} \sim 
X^{2j}$ and $B= (\bigcup_{j=0}^{k-1} X^{2j} \sim X^{2j+1})\cup X^{2k}$. 
We again have that (3) and (4) hold, so $i(A) = i(B) =n$ and 
$\|\chix_A\|_D = n$, $\|\chix_B\|_D = n+1$ by Theorem~2.2.\qed 
\enddemo 

\proclaim{Corollary 2.7} 
Let $K$ satisfy the hypotheses of Proposition 2.6. 
Then $D(K) \ne B_{1/2} (K)$. 
Moreover there exists a set $A\subset K$ with $\chix_A \in B_1(K) \sim 
B_{1/2}(K)$.  
\endproclaim 

\demo{Proof} 
Since $D(K)\subset B_{1/2}(K)$ and the $D$-norm is stronger than the 
sup-norm, were $D(K) = B_{1/2}(K)$, the norms would be equivalent by 
a theorem of Banach. 
But of course Proposition~2.6 shows they are not. 
Finally, letting the sets $X^j$ be as in the proof of 2.6, 
the set $A= \bigcup_{j=0}^\infty X^{2j}\sim X^{2j+1}$ has the 
property that $i(A) = \infty$, hence $\chix_A \notin B_{1/2}(K)$; 
since $\partial^{\omega+1}A =\emptyset$, $\chix_A$ belongs to $B_1(K)$.\qed 
\enddemo 

\demo{Remarks} 
1. We may also easily construct explicit functions which belong to 
$B_{1/2}(K)$ but not to $D(K)$. 
Thus, if $K$ is as in 2.6, we may choose disjoint open subsets 
$U_1,U_2,\ldots$ of $K$, and for each $n$, a subset $A_n$ of $U_n$ 
with $i(A_n) =n$ and $\|\chix_{A_n}\|_D = n$. 
Now let $f= \sum_{n=1}^\infty {1\over \sqrt n} \chix_{A_n}$. 
Then $f$ is clearly the uniform limit of $D(K)$-functions, hence is 
in $B_{1/2}(K)$. 
However $\|f\|_D \ge \|f|\U_n\|_D\ge {n\over \sqrt n} = \sqrt{n}$ 
for all $n$, hence $f\notin D(K)$. 

2. We obtain, in the next section, that if however $X^{(n)}=\emptyset$ 
for some $n$, then every bounded function on $K$ belongs to $D(K)$; 
thus $D(K) = B_{1/2}(K) = B_1 (K) = \ell^\infty (K)$. 
\enddemo 

Our final result shows that a function $\varphi :\complex\to\complex$ 
operates on $D(K)$ (for interesting $K$) precisely when $\varphi$ is 
locally Lipschitz; that is, $\varphi$ is Lipschitz on compact subsets 
of $\complex$. 
(If $W\subset \complex$ and $\varphi :W\to\complex$ is Lipschitz, we 
set $\|\varphi\|_{\Lip} = \sup \{ {|\varphi (z_1)-\varphi(z_2)| \over 
|z_1-z_2|} : z_1\ne z_2$, $z_1,z_2\in W\}$. 
$\|\varphi\|_{\Lip}$ is called the Lipschitz constant of $\varphi$.) 

\proclaim{Proposition 2.8} 
Let  $\varphi :K\to\complex$ be a given function. 

{\rm (a)}  If $\varphi$ is locally Lipschitz, then $\varphi\circ f\in D(K)$ 
for all $f\in D(K)$. 

{\rm (b)} If conversely $K$ satisfies the hypotheses of Proposition~2.6 and 
$\varphi\circ f \in D(K)$ for all $f\in D(K)$, then $\varphi$ is locally 
Lipschitz. 
\endproclaim 

\demo{Remark} 
(a) is due jointly to F.~Chaatit and the author (cf.\ \cite{C}). 
We note also that the proof of (a) yields that for $f\in D(K)$ and $\varphi$ 
locally Lipschitz
$$\|\varphi\circ f\|_{qD} \le \|\varphi\mid W\|_{\Lip} \|f\|_{qD}\text{ where } 
W= \{z:|z| \le \|f\|_\infty\}\ .$$ 

(a) Let $\lambda = \|f\|_\infty$ and suppose that $\varphi$ has Lipschitz 
constant at most $M$ on $\{z:|z|\le \lambda\}$. 
Then we claim that 
$$\|\varphi \circ f\|_D \le \|\varphi\|_\infty + M\|f\|_D\ . 
\leqno(5)$$ 
Let $\varep>0$. By Corollary 1.13, we may choose $(f_n)$ in $C_b(K)$ 
with $f_n\to f$, $\|f_n\|_\infty \le \lambda$ for all $n$, and 
$$|f_1| + \sum |f_{n+1} -f_n| \le \|f\|_D +\varep\ .
\leqno(6)$$ 

Since $\varphi$ is continuous, $\varphi (f_n) \to \varphi (f)$ and by 
the definition of $M$,   
$$|\varphi (f_{n+1}) - \varphi (f_n)| \le M|f_{n+1} -f_n|\text{ for all } n\ .
\leqno(7)$$ 
Hence by (6) and (7), 
$$|\varphi (f_1)| + \sum |\varphi (f_{n+1}) - \varphi (f_n)| 
\le \|\varphi\|_\infty + M(\|f\|_D +\varep)\ . 
\leqno(8)$$ 
Of course (8) yields the estimate in (5), as well as the estimate in the 
above remark. 

(b) Suppose to the contrary that $\varphi$ operates on $D(K)$ but $\varphi$ 
is not locally Lipschitz. 
Suppose first that $\varphi$ is continuous. 
Then by compactness, it follows that there exist a scalar $z_0$ and 
for all $n$ scalars, $u_n,v_n$ with $u_n\ne v_n$ so that 
$$u_n,v_n\to z_0\ \text{ and }\ {|\varphi (u_n)-\varphi(v_n)| \over 
|u_n-v_n|} \to \infty \text{ as } n\to\infty\ . 
\leqno(9)$$ 
Now we may assume without loss of generality that 
$z_0 = \varphi(z_0) = 0$. 
Indeed, simply replace $\varphi$ by $\psi$, where $\psi (z) = \varphi (z+z_0) 
- \varphi (z_0)$, if necessary; since $\varphi$ operates on $D$, so does 
$\psi$.  

Next, suppose $a,b$ are distinct complex numbers, $n$ is a positive 
integer, and $K= X_0\supset X_1\supset X_2\supset\cdots\supset X_n$ are 
closed non-empty subsets of $K$ with $X_i$ nowhere dense in $X_{i-1}$ 
for all $i$, $1\le i\le n$. 
Let $X_{n+1}= \emptyset = X_m$ for all $m>n$ and define $g$ by 
$g=0$ on $X_0\sim X_1$ and 
$$\left\{\eqalign{ 
&\quad g= a\text{ on } X_{2i-1} \sim X_{2i}\cr 
&\quad g= b\text{ on } X_{2i}\sim X_{2i-1}\cr
&\text{for all }i = 1,2,\ldots\ ;\cr}
\right. 
\leqno(10)$$ 
we then have 
$$n|b-a| -2|b| \le \|g\|_D \le (n+1) |b-a| + 2|b| 
\leqno(11)$$ 

Indeed, the proof of Proposition 2.2 yields that 
$n|b-a| \le \|g-b\chix_{X_1}\|_D \le (n+1) |b-a|$, 
since $\frac1{b-a} (g-b\chix_{X_1}) = \chix_A$ where $i(A) = n$. 

Now by (9), and the continuity of $\varphi$, for each $j$ we may choose 
distinct complex numbers $a_j$ and $b_j$ with 
$$\leqalignno{
&|\varphi (a_j)|,|\varphi(b_j)| \le 1 &(12)\cr 
&|a_j|,|b_j|\le {1\over 2j} &(13)\cr 
&|\varphi (b_j) - \varphi (a_j)| > j|b_j-a_j|\ .&(14)\cr}$$ 
Then choose $n_j$ a positive integer with 
$${1\over n_j+1} < |b_j-a_j| \le {1\over n_j}\ .
\leqno(15)$$ 
By (13), we have $n_j\ge j$ for all $j$.

Next, we may choose disjoint open sets $U_1,U_2,\ldots$ in  $K$ with 
$U_j^{(n_j)} \ne\emptyset$ for all $j$. 
Finally, fix $j$, let $n=n_j$ and (by the proof of Proposition~2.2), 
choose sets $X_1\supset \cdots\supset X_n$ as above, with $X_1$ a 
closed nowhere dense subset of $U_j$. 
Then define $g$ on $U_j$ by (10), where $a=a_j$, $b=b_j$. 
Of course we simply define $g$ to be zero off the union of the $U_j$'s; 
again fixing $j$, we have by (11), (13) and (15) that 
$$\eqalign{ 
\| g\mid U_j\|_D &\le  (n_j+1) |b_j-a_j| + 2|b_j|\cr 
&\le 2+ {1\over j} \le 3\ .\cr}
\leqno(16)$$ 
thus, $g\in D$, with $\|g\|_D\le 3$. 

Now again, fixing $j$, we have by (11), (12) and (14) that 
$$\eqalign{ 
\|(\varphi\circ g)\mid U_j\|_D 
&\ge n_j |\varphi (b_j) - \varphi (a_j)| - 2|\varphi (b_j)|\cr 
& > jn_j |b_j -a_j| -2\cr 
& > j {n_j\over n_j+1} -2\ \text{ by (15).}\cr} 
\leqno(17)$$ 
Thus $\|\varphi\circ g\|_D \ge \sup_j j {n_j\over n_{j+1}} -2 
= \infty$, so $\varphi \circ g \notin D$. 
This contradiction completes the proof of (b), 
in case $\varphi$ is continuous. 
Now if $\varphi$ is not continuous, but $\varphi$ still operates, then 
without loss of generality (by making the same modification as above), 
we may assume that $\varphi (0)=0$ and $\varphi$ is discontinuous at $0$. 
Thus we may choose numbers $a_1,a_2,\ldots$ and $\delta>0$ so that 
for all $n$, 
$$|a_n| < {1\over 2^n} \ \text{ and }\ |\varphi (a_n)| >\delta\ . 
\leqno(18)$$ 

Next, let the sets $K= X^0 \supset X^1 \supset X^2 \cdots$ be chosen 
as in the proof of 2.6, and define $f$ by 
$$f= \sum_{n=1}^\infty a_n \chix_{X^{2n-1} \sim X^{2n}}\ . 
\leqno(19)$$ 
Since $\|a_n \chix_{X^{2n-1}\sim X^{2n}}\|_D \le {1\over 2^{n-1}}$ for all 
$n$, and $D(K)$ is a Banach space, $f\in D(K)$. 
However 
$$\varphi \circ f = \sum_{n=1}^\infty \varphi (a_n) \chix_{X^{2n-1}\sim 
X^{2n}}\ . 
\leqno(20)$$ 
It then follows from (18), by the same argument as in the proof of 2.6, that 
$$X_n = \os_n (\varphi \circ f,\delta)\ \text{ for all }\ n\ . 
\leqno(21)$$ 

Thus by Proposition 1.9, $\varphi \circ f$ does not even belong to 
$B_{1/2}(K)$. 
This contradiction completes the proof.\qed
\enddemo 

\head \S3. The transfinite oscillations; properties and first applications. 
\endhead 

In this section, we introduce the transfinite oscillations as a tool in 
studying $D(K)$. 
As noted in the introduction, our main applications (given in Section~4), 
really only involve the finite oscillations. 
However we give some initial general results here; deeper applications 
involving arbitrary ordinals are given in \cite{R3}. 

We begin by listing some useful permanence properties of the transfinite 
oscillations. 

\proclaim{Proposition 3.1} 
Let $f,g$ be given complex-valued functions on $K$, $t$ a complex number, 
and $\alpha,\beta$ non-zero ordinals. 

\iitem{\rm (a)} $\osc_\alpha f$ is an upper semi-continuous 
$[0,\infty]$-valued function: if $\alpha\le\beta$, then $\osc_\alpha f 
\le \osc_\beta f$. 

\iitem{\rm (b)} $\osc_\alpha tf = |t| \osc_\alpha f$ and $\osc_\alpha (f+g) 
\le \osc_\alpha f +\osc_\alpha g$.

\iitem{\rm (c)} $\osc_\alpha fg\le U|f|\osc_\alpha g + U|g| \osc_\alpha f$.  

\iitem{\rm (d)} $\osc_\alpha |f| \le \osc_\alpha f$.  

\iitem{\rm (e)} If $\osc_\alpha f= \osc_{\alpha+1} f$, then $\osc_\alpha f 
= \osc_\beta f$ for all $\beta >\alpha$. Moreover if $f$ is real-valued, 
this happens if and only if $\osc_\alpha f\pm f$ are both upper 
semi-continuous functions. 

\iitem{\rm (f)} If $f$ is semi-continuous, then $\osc_\alpha f= \osc f= 
\uosc f = \oosc f$. 

\iitem{\rm (g)} $\osc_\alpha f = \osc_\alpha (f+\varphi)$ for all 
continuous $\varphi :K\to \complex$. 
\endproclaim 

\demo{Proof} 
The assertions up to the ``moreover'' statement in (e), as well as the 
assertion (g), are easily proved by transfinite induction. 
For example, to see the assertion in (c), (where we define $\infty \cdot 0 
= \infty$), suppose $\beta$ is a non-zero ordinal and the inequality 
proved for all ordinals $\alpha <\beta$. 
Now suppose first that $\beta$ is a successor ordinal, say $\beta = \alpha+1$. 
Then we have for $x,y$ in $K$ that 
$$\eqalign{
|f(y) g(y) - f(x)g(x)| 
& \le |f(y)| \, |g(y)-g(x)| + |g(x)|\, |f(y)-f(x)|\cr 
& \le U|f|(y) |g(y)- g(x)| + U|g|(x) |f(y)-f(x)|\ .\cr} 
\leqno(3)$$ 

Now suppose none of the terms $U|f|(x)$, $\widetilde{\osc}_\beta f(x)$, 
$U|g|(x)$, and $\osc_\beta g(x)$ are equal to infinity. 
Then by induction hypothesis, we have 
$$\eqalign{
\widetilde{\osc}_{\alpha+1} fg(x) 
& \le \olim_{y\to  x} U|f| (y) |g(y)-g(x)| + Uf(y) \osc_\alpha g(y)\cr 
&\qquad + \olim_{y\to x} U|g| (x) |f(y) - f(x)| + U|g|(x)\osc_\alpha f(y)\cr 
&\qquad + \olim_{y\to x} (U|g|(y) - U|g|(x)) \osc_\alpha f(y)\ .\cr}
\leqno(4)$$ 
The last term above is at most zero, since 
$\olim_{y\to x} [U|g|(y) - U|g|(x)] 
\olim_{y\to x} \osc_\alpha f(y)=0$, by the upper semi-continuity of $U|g|$. 
Thus (4) yields  
$$\widetilde{\osc}_\beta f g(x) \le U|f| (x) \widetilde{\osc}_\beta g(x) 
+ U|g| (x) \widetilde{\osc}_\beta f(x)\ .
\leqno(5)$$ 
Of course (5) holds trivially by our convention, if any of these terms 
equal infinity. 

Again, if $\beta$ is a limit ordinal, we obtain that (5) holds, by simply 
taking sups. 
The desired inequality (c) now follows for $\beta$, by taking upper 
semi-continuous envelopes. 

Next, we pass to the ``moreover'' assertion in (e), and the proof of (f). 
(These may also be found in \cite{R1}, but for completeness we give 
them again, here.) 

To prove the ``moreover'' assertion in (e), we first note that $\osc_\alpha f
\le \widetilde{\osc}_{\alpha+1} f\le \osc_{\alpha+1} f$. 
It then follows that 
$$\osc_{\alpha+1} f = \osc_\alpha f \text{ if and only if } 
\widetilde{\osc}_{\alpha+1} f= \osc_\alpha f 
\leqno(6)$$ 
(for if the latter equality holds, then since $\osc_\alpha f$ is upper 
semi-continuous, $\osc_{\alpha+1} f= U\widetilde{\osc}_{\alpha+1} f 
= U\osc_\alpha f = \osc_\alpha f$). 

Now assume $f$ is real valued, and suppose first that $\osc_\alpha f = 
\osc_{\alpha+1} f$. 
To see that $\osc_\alpha f+f$ is upper semi-continuous, let $x\in K$ and 
$(y_n)$ be a sequence in $K$ with $y_n\to x$. 
Then 
$$\eqalign{
&\olim\limits_{n\to\infty} \osc_\alpha f(y_n) + f(y_n) - f(x)\cr 
&\qquad \le \olim\limits_{n\to\infty} \osc_\alpha f(y_n) + 
|f(y_n) -f(x)|\cr 
&\qquad \le \widetilde{\osc}_{\alpha+1}  f(x) = \osc_\alpha f(x) 
\text{ by (6).}\cr}$$ 
Hence 
$$\olim\limits_{n\to\infty} \osc_\alpha f(y_n) + f(y_n) \le 
\osc_\alpha f(x) + f(x)\ ,$$ 
proving $\osc_\alpha f +f$ is upper semi-continuous. 
Since $\osc_\beta f = \osc_\beta -f$ for all $\beta$ by 3.4(b), it follows 
immediate upon replacing $f$ by $-f$ that also $\osc_\alpha f-f$ is upper 
semi-continuous. 

Now suppose conversely that $\osc_\alpha \pm f$ are upper semi-continuous, 
yet $\osc_{\alpha+1}f\ne \osc_\alpha f$. 
Then by (6) we may choose $x\in K$ so that $\widetilde{\osc}_{\alpha+1}f(x) 
> \osc_\alpha f(x)$. 
But $\widetilde{\osc}_{\alpha+1} f(x) = \olim_{y\to x} |f(y) -f(x)| 
+ \osc_\alpha f(y) = \max\{ \olim_{y\to x} (f(y)-f(x)) + \osc_\alpha f(y)$, 
$\olim_{y\to x} (f(x) - f(y)) + \osc_\alpha f(y)\}$. 
Thus either 
$$\leqalignno{
&\olim\limits_{y\to x} f(y) - f(x) + \osc_\alpha f(y)>\osc_\alpha f(x) 
&\text{\rm (7)(i)}\cr 
\noalign{\hbox{\rm or}} 
&\olim\limits_{y\to x} f(x) - f(y) + \osc_\alpha f(y) >\osc_\alpha f(x)\ . 
&\text{\rm (7)(ii)}\cr}$$
But if (7)(i) holds, $f+\osc_\alpha f$ is not upper semi-continuous, 
while if (7)(ii) holds, $(-f)+\osc_\alpha f$ is not upper semi-continuous. 

Finally, to prove (f), suppose without loss of generality that $f$ is 
upper semi-continuous. 
(For if $f$ is lower semi-continuous, $-f$ is upper semi-continuous, and 
$\osc_\alpha f= \osc_\alpha -f$.) 
But then $f= Uf$ and hence
$$\uosc f = f - Lf = \osc f\qquad (= \oosc f)\ .
\leqno(8)$$ 
But then $\osc f +f=f-Lf +f=2f-Lf$ and $\osc f-f = -Lf$; thus 
$\osc f\pm f$ are both upper semi-continuous, so (f) follows from (e).\qed 
\enddemo 

\demo{Remarks} 
1. Actually, an appropriate version of (e) holds for 
complex valued $f$ as well. 
The result: $\osc_\alpha f= \osc_{\alpha+1}f$ 
{\it if and only if $\osc_\alpha f+ \Re \mu f$ is upper semi-continuous 
for all scalars $\mu$ with $|\mu| =1$.} 
To see this, suppose first $\osc_\alpha f= \osc_{\alpha+1}f$. 
Since $\osc_\alpha f= \osc_\alpha \mu f$ if $|\mu| =1$, it obviously 
suffices to prove (in general) that $\osc_\alpha f+\Re f$ is upper 
semi-continuous. 
Were this false, we could choose $x$ and $(x_n)$ converging to $x$ with 
$$\lim_{n\to\infty} \Re f(x_n) + (\osc_\alpha f)(x_n) > \Re f(x) + 
\osc_\alpha f(x)\ .$$ 
But then 
$$\eqalign{\widetilde{\osc}_{\alpha+1} f(x) & \ge \lim_{n\to\infty} 
|f(x_n) - f(x)| + \osc_\alpha f(x_n)\cr 
&\ge \lim_{n\to\infty} \Re \bigl[ f(x_n)- f(x)\bigr] + \osc_\alpha f(x_n) 
> \osc_\alpha f(x)\ ,\cr}$$ 
contradicting our hypothesis. 
Conversely, were it false that $\osc_\alpha f= \osc_{\alpha+1} f$, choose 
$x\in K$ and $(x_n)$ converging to $x$ so that 
$$\widetilde{\osc}_{\alpha+1} f(x) = \lim_{n\to\infty} 
|f(x_n) -f(x)| + \osc_\alpha f(x_n) > \osc_\alpha f(x)\ .$$ 
We may assume without loss of generality that $L \dfeq \lim_{n\to\infty} 
|f(x)-f(x_n)|$ exists; thus also $M \dfeq \osc_\alpha f(x_n)$ exists, 
and since $M\le \osc_\alpha f(x)$, $L>0$ (and $M<\infty$). 
But then it follows that we may choose $\mu$ with $|\mu|=1$ so that 
setting $z_n = f(x_n) - f(x)$ for all $n$, then after pursuing to a 
subsequence if necessary, we have that $L= \lim_{n\to\infty} \Re \mu z_n$. 
(If $L=\infty$, $\mu = \pm1$ or $\pm i$ works; otherwise let $\mu_n z_n = 
|z_n|$, $|\mu_n|=1$, and after passing to a subsequence, suppose $\mu_n\to 
\mu$.) 
But then we have that 
$$\widetilde{\osc}_{\alpha+1} f(x) = \lim_{n\to\infty} \Re 
\bigl[ \mu f(x_n) - \mu f(x)\bigr] + \osc_\alpha f(x_n) > \osc_\alpha f(x)\ ,$$ 
which shows that $\osc_\alpha f+ \Re \mu f$ is not upper semi-continuous. 

2. For $\alpha$ an arbitrary ordinal, we let $D_\alpha (K)$ denote the set 
of all bounded $f:K\to \complex$ with $\osc_\alpha f$ bounded (for 
$\alpha\ge \omega$; if $0\le\alpha <\omega$, let $D_\alpha (K) = B_{1/2}(K)$). 
It follows from 3.1 (b), (c), and transfinite induction that $D_\alpha (K)$ 
is a Banach algebra under the norm  $\|f\| \dfeq \|f\|_\infty + 
\|\osc_\alpha f\|_\infty$. Corollary~4 of the Introduction yields that 
$B_{1/4}(K)=D_\omega (K)$. 
If $K$ is separable, Theorem~3.2 below yields that 
$\bigcap_{\alpha<\omega_1} D_\alpha (K) = D(K)$. 
The Banach algebras $D_\alpha (K)$, as well as the transfinite analogues 
of $B_{1/4}(K)$, are studied in \cite{R2}. 
\enddemo 

We may now formulate a fundamental structural result for $D(K)$. 
(Recall that for a metric space $K$, $\wt K$ denotes the smallest possible 
cardinality of a base for the open subsets of $K$; thus $\wt K=\omega$ 
if $K$ is infinite separable; of course we identify cardinals with 
``initial'' ordinals.) 

\proclaim{Theorem 3.2} 
Let $K$ be an infinite metric space and $f:K\to\complex$ be a bounded 
function. 
Let $\mu$ be the least cardinal number with $\mu>\wt K$. 
There exists an ordinal $\alpha$ with $\alpha <\mu$ so that 
$\osc_\alpha f= \osc_\beta f$ for all $\beta >\alpha$. 
Letting $\tau$ be the least such $\alpha$, then $f$ is in $D(K)$ if and 
only if $\osc_\tau f$ is bounded. 
When $f$ is real valued and this occurs, then 
$$\|f\|_D = \|\, |f| + \osc_\tau f\|_\infty
\leqno(9)$$ 
and 
$$\|f\|_{qD} = \|\osc_\tau f\|_\infty\ . 
\leqno(10)$$ 
Moreover setting $\lambda = \|\, |f| + \osc_\tau f\|_\infty$, 
$u= {\lambda -\osc_\tau f+f\over 2}$ and 
$v= {\lambda -\osc_\tau f-f\over 2}$, $u,v$ are non-negative lower 
semi-continuous functions with $f=u-v$ and $\|f\|_D= \|u+v\|_\infty$. 
\endproclaim 

The proof is a minor modification of the one given in \cite{R1}. 
We summarize the ingredients, but refer to \cite{R1} for certain details. 
We also note that the qualitative part of Theorem~3.2 (i.e., the first 
four sentences in its statement) follows from the earlier work of 
A.L.~Kechris and A.~Louveau \cite{KL}. 

\proclaim{Lemma 3.3} 
Let $u,v$ be non-negative bounded lower semi-continuous functions  on $K$. 
Then for all ordinals $\alpha$, 
$$\osc_\alpha (u-v) \le \osc (u+v)\ . 
\leqno(11)$$ 
\endproclaim 

For the proof, see Lemma 3.6 of \cite{R1}. 
The next stability result appears in \cite{KL} for compact metric spaces; 
the generalization to arbitrary metric spaces presents no problem. 

\proclaim{Lemma 3.4} 
Let $\mu$ be as in Theorem 3.2, and $(\varphi_\alpha)_\alpha <\mu$ be a 
family of upper semi-continuous extended real-valued functions defined 
on $K$ so that $\varphi_\alpha \le\varphi_\beta$ for all $\alpha<\beta$. 
Then there is an ordinal $\alpha <\mu$ so that $\varphi_\alpha =\varphi_\beta$ 
for all $\beta >\alpha$. 
\endproclaim 

\demo{Proof} 
This is essentially the same as the argument for Lemma~3.7 of \cite{R1}, 
but we give the argument for the sake of completeness. 
Suppose not. 
Then by renumbering, we may assume that 
$$\varphi_\alpha \ne \varphi_{\alpha+1}\ \text{ for all }\ \alpha <\mu\ . 
\leqno(12)$$ 

Now let $\B$ be a base for the open subsets of $K$ with $\card \B =\wt K$. 
Fix $\alpha <\mu$; by (12), we may choose $x= x_\alpha \in K$ with  
$\varphi_\alpha (x) <\varphi_{\alpha+1}(x)$. 
Then by the upper semi-continuity of $\varphi_\alpha$, choose $U_\alpha \in\B$ 
so that $x\in U_\alpha$ and 
$$\lambda_\alpha \dfeq \sup \varphi_\alpha (U_\alpha ) < \varphi_{\alpha+1} 
(x)\ .
\leqno(13)$$ 

Now we may choose an uncountable subset $\Gamma$ of $\mu$ so that 
$$U_\alpha = U_\beta \dfeq U\ \text{ for all }\ \alpha,\beta\in\Gamma\ . 
\leqno(14)$$ 

Indeed, for each $U\in \B$, let $\Gamma_U = \{\alpha <\mu : U_\alpha = U\}$. 
Then $\mu = \bigcup_{U\in \B} \Gamma_U$. 

Now in fact, we must have that $\card \Gamma_U = \mu$ for some $U\in\B$; 
for otherwise letting $\kappa = \card \B$, $\mu$ would be  
at most $\card \kappa \times\kappa = \kappa$. 

Finally, we have that for $\Gamma$ satisfying  (14), 
$$\lambda_\alpha <\lambda_\beta \ \text{ if }\ 
\alpha<\beta\ ,\ \alpha,\beta\in\Gamma\ . 
\leqno(15)$$ 
Indeed, fixing $\alpha<\beta$ in $\Gamma$ and letting $x=x_\alpha$ as above, 
we have that $\lambda_\alpha <\varphi_{\alpha+1} (x) \le \varphi_\beta (x) 
\le \sup \varphi_\beta (U) = \lambda_\beta$. 
But of course since $\Gamma$ is uncountable, (15) is impossible.\qed 
\enddemo 

\demo{Proof of Theorem 3.2} 

The first assertion follows immediately from the preceding Lemma and 3.1(a). 
Next, assume $f$ is real-valued. 
Then if $f$ is in $D(K)$, $f=u-v$ for some lower semi-continuous bounded 
functions $u$ and $v$; then $\osc_\tau f\le \osc (u+v)$ by Lemma~3.3, so 
$\osc_\tau f$ is bounded. 
Conversely, if $\osc_\tau f$ is bounded, let $\lambda, u$ and $v$ be as 
in the statement of 3.2. 
Then trivially $f= u-v$ and $u,v$ are non-negative. 
Since $\osc_\tau f= \osc_{\tau+1} f$, $\osc_\tau f\pm f$ are upper 
semi-continuous by Proposition~3.1(e), which implies the lower 
semi-continuity of $u$ and $v$; since $u,v$ are bounded, we have that 
$f$ is in $D(K)$. 

Finally, for the norm identity, we first note that 
$$\|f\|_D \le \|u+v\|_\infty = \|\lambda-\osc_\tau f\|_\infty \le\lambda\ . 
\leqno(15)$$ 
For the reverse inequality, let $\varep >0$, and choose $g,h$ non-negative 
lower semi-continuous with $f=g-h$ and $\|g+h\|_\infty \le \|f\|_D+\varep$. 

Then 
$$\eqalign{ |f| + \osc_\tau f
& = |g-h| + \osc_\tau (g-h)\cr 
&\le |g-h| + \osc (g+h)\ \text{ by Lemma 3.3}\cr 
&= |g-h| + U(g+h) - (g+h)\cr 
& \le U(g+h)\ .\cr}$$ 
Hence $\lambda = \|\,|f| + \osc_\tau f\|_\infty \le \|U(g+h)\|_\infty 
= \|g+h\|_\infty \le \|f\|_D +\varep$. 
Since $\varep >0$ is arbitrary, (9) is proved. 

To prove (10), we first observe that 
$$\|\osc_\tau f\|_\infty \le \|f\|_{qD} \ . 
\leqno(16)$$ 
Indeed, if $\varphi \in C_b(K)$, then by what we have already shown and 
Proposition~3.1(g), for any ordinal $\alpha$, 
$$\|\osc_\alpha f\|_\infty = \|\osc_\alpha (f+\varphi) \|_\infty 
\le \|f+\varphi \|_D\ .$$
Thus taking the inf over $\varphi$ in $C_b(K)$ and letting $\alpha=\tau$, 
we obtain (16). 

For the other inequality, let $\varep>0$. 
We shall prove that for all $x\in K$, there exists $U$ an open neighborhood 
of $x$ so that 
$$\|f\|_{qD(U)} \le \osc_\tau f(x) +\varep\ . 
\leqno(17)$$
Once this is proved, we have by the localization principle, 
Proposition~1.14(b), that 
$$\|f\|_{qD} \le \|\osc_\tau f\|_\infty +\varep\ . 
\leqno(18)$$ 

Of course since $\varep>0$ is arbitrary, we then obtain the opposite 
inequality  in (16) as well, so (10) follows. 

Now fixing $x\in K$, since $\widetilde{\osc}_\tau f(x) = \osc_\tau f(x)$, 
we may choose $U$ an open neighborhood of $x$ so that 
$$|f(y) -f(x)| + \osc_\tau f(y) \le\osc_\tau f(x) +\varep 
\ \text{ for all }\ y\in U\ . 
\leqno(19)$$ 
Now setting $\beta = f(x)$, we have proved that 
$$\sup_{u\in U} |f(u)-\beta| + \osc_\tau f(u) \le\osc_\tau f(x)+\varep\ . 
\leqno(20)$$ 
But of course $\osc_\tau (f-\beta) = \osc_\tau f= \osc_{\tau+1}f = 
\osc_{\tau+1} (f-\beta)$. 
In particular, this holds on $U$. 
Thus, by (9) (i.e., the first norm assertion of the theorem), we have 
by (20) that 
$$\|f-\beta\|_{D(U)} \le \osc_\tau f(x) +\varep\ ,$$ 
which of course yields (17). 

The theorem is thus established for real-valued functions. 
Now suppose $f$ is complex-valued. 
Then it is easily established by transfinite induction that if $g=\Re f$ or 
$\Im f$, then 
$$\osc_\alpha g\le \osc_\alpha f\ \text{ for all ordinals }\ \alpha\ . 
\leqno(21)$$ 
Thus we obtain that $\osc_\beta g\le \osc_\beta f= \osc_\tau f$ for 
all $\beta >\tau$ (where $\tau$ is as in the statement of the Theorem). 
Hence if $\osc_\tau f$ is bounded and $\beta$ is such that $\osc_{\beta+1}g 
= \osc_\beta g$ for both $g=\Re f$ and $g= \Im f$, then $\osc_\beta \Re f$, 
$\osc_\beta \Im f$ are both bounded, whence $f$ is in $D(K)$ since its real 
and imaginary parts belong to $D(K)$. 
Of course {\it if\/} $f$ is in $D(K)$, then we trivially have that $\Re f$, 
$\Im f$ belong to $D(K)$, and then $\osc_\tau f\le \osc_\tau \Re f + 
\osc_\tau \Im f$ by Proposition~3.4b; thus $\osc_\tau f$ is bounded. 
This completes the proof of Theorem~3.2.\qed
\enddemo 

We next give several complements and remarks concerning Theorem~3.2. 

Let $f:K\to\complex$ be a general function. 
We define the $D$-index of $f$, denoted $i_D(f)$, to be the least ordinal 
$\alpha$ so that $\osc_\alpha f= \osc_{\alpha+1}f$. 
We show in \cite{R2} that if $K$ is an uncountable compact metric space, 
then for $f\in D(K)$, $i_Df$ may be any countable ordinal. 
(An analogous index and result were previously obtained in \cite{KL}.) 

Evidently if $n=i_D(f)<\infty$, and $f$ is real-valued, we obtain from 
Theorem~3.2 that $\|f\|_D \le (2n+1) \|f\|_\infty$, while if $f$ is 
non-negative, $\|f\|_D \le (n+1)\|f\|_\infty$ (since then 
$\|\osc_nf\|_\infty \le n\|f\|_\infty$). 
In particular, by Proposition~3.1, we recapture our observation at the end of 
Section~2 that if $f$ is semi-continuous, $\|f\|_D \le 3\|f\|_\infty$; 
moreover if $f$ is non-negative, then $\|f\|_D \le 2\|f\|_\infty$. 

Suppose $\alpha = i_Df$ and $f\in D(K)$, with $f$ real-valued. 
Since $\osc_\alpha f\pm f$ are both upper semi-continuous, it follows that 
$\osc_\alpha f+|f| = \max \{\osc_\alpha f+f$, $\osc_\alpha f-f\}$ is upper 
semi-continuous, and since $f^+ = {|f|+f\over2}$, also $\osc_\alpha f+f^+$, 
$\osc_\alpha f+f^-$ are upper semi-continuous. 
Thus we obtain $f= u-v$ where $u= \osc_\alpha f+f^+$, $v= \osc_\alpha f+f^-$; 
$u,v$ are non-negative upper semi-continuous and again 
$\|u+v\|_\infty = \|f\|_D$. 

We also note that for $f\in D(K)$ complex valued and $\alpha = i_D(f)$, 
we have 
$$\tfrac12 \|f\|_D \le \|\, |f| + \osc_\alpha f\|_\infty \le 2\|f\|_D\ . 
\leqno(22)$$ 
Indeed, let $\beta = \max\{ \alpha, i_D\Re f,\ i_D\Im f\}$. 
Then 
$$\eqalign{ \|f\|_D 
& \le \|\Re f\|_D + \|\Im f\|_D\cr 
& = \|\, |\Re f| + \osc_\beta \Re f\|_\infty 
+ \|\, | \Im f| + \osc_\beta \Im f\|_\infty \cr 
&\le 2\|\, |f| + \osc_\beta f\|_\infty \ \text{ by (21))} \cr
& = 2\|\, |f| + \osc_\alpha f\|_\infty\ .\cr}$$ 
On the other hand, $|f|+\osc_\alpha f\le |\Re f| + \osc_\alpha \Re f 
+ |\Im f| + \osc_\alpha \Im f$ by Proposition~3.1(b). Hence 
$$\eqalign{ 
\|\, |f| + \osc_\alpha f\|_\infty 
& \le \|\, |\Re f| + \osc_\alpha \Re f\|_\infty 
+ \|\, |\Im f| + \osc_\alpha \Im f\|_\infty\cr 
& \le \|\Re f\|_D + \|\Im f\|_D \ \text{ by Theorem 3.5}\cr 
&\le 2\|f\|_D\ .\cr}$$ 

\demo{Remark} 
Actually, we may obtain a decomposition 
of an arbitrary complex-valued $D$-function 
into a linear combination of semi-continuous functions, without passing to the 
possibly higher indices of its real and imaginary parts; 
also the considerations 
about absolute values hold as well. 
Thus, suppose $f\in D(K)$ and $\alpha = i_Df$. 
We then have (by the Remark following the proof of Proposition~3.2) that 
if $F=\osc_\alpha f$, then 
$$F + \Re \mu f\ \text{ is upper semi-continuous for all scalars $\mu$ 
with $|\mu|=1$.}
\leqno(*)$$ 
Thus in particular, $F\pm \Re f$ and $F\pm \Re if$ are upper semi-continuous, 
so setting  $u= {F+\Re f\over2}$, $v= {F-\Re f\over2}$, 
$\tilde u= {F-\Re if\over 2}$, $\tilde v= {F+\Re if\over2}$, 
then $u,v,\tilde u,\tilde v$ are all upper semi-continuous, and 
$f= (u-v)+i(\tilde u-\tilde v)$. 
Finally, we note that if 
{\it $F$ is any non-negative bounded upper semi-continuous function 
satisfying $(*)$, then it follows that $F+|f|$ is upper semi-continuous.} 
Thus $\osc_\alpha f+|f|$ is upper semi-continuous. 
\enddemo 

To see the above claim, suppose to the contrary that $F+|f|$ is not 
upper semi-continuous. Then choose $x$ and $(x_n)$ converging to $x$ with 
$$\lim_{n\to\infty} F(x_n) + |f|(x_n) > F(x) + |f|(x)\ .$$ 
By passing to subsequences, we may assume that $L \dfeq \lim_{n\to\infty}  
F(x_n)$ and $M = \lim_{n\to\infty} |f|(x_n)$ exist. 
Since $L\le F(x)$ by the upper semi-continuity of $F$, we have that $M>0$.
Choose $\mu_n$ with $|\mu_n|=1$ and $|f|(x_n) =\Re \mu_n f(x_n)$ 
for all $n$. 
By passing to a further subsequence, we may assume that $\mu\dfeq \lim \mu_n$ 
exists. Since $f$ is bounded, $\Re (\mu_n-\mu) f(x_n)\to 0$, whence 
$M= \lim_{n\to\infty} \Re \mu f(x_n)$. 
But then 
$$\lim_{n\to\infty} \bigl[ F(x_n) + \Re \mu f(x_n)\bigr] > 
F(x) + |f|(x) \ge F(x) + \Re \mu f(x)\ ,$$ 
contradicting $(*)$.  

The next result yields an interpretation of the function $\osc_\alpha f$, 
for $\alpha = i_Df$, and also shows the quotient norm $\|f\|_{qD}$ is 
always attained (for real-valued $f$). 

\demo{Definition} 
Given $f:K\to \complex$ a bounded function and $x\in K$, set 
$$\|f\|_{qD(x)} 
= \inf \{\|f\mid U\|_{qD} :U\text{ is an open neighborhood of }x\}\ .$$
\enddemo 

\proclaim{Corollary 3.5} 
Let $f:K\to\real$ be a bounded function and $\alpha = i_Df$. 

\iitem{\rm (a)} $\osc_\alpha f(x) = \|f\|_{qD(x)}$ for all $x\in K$. 
\iitem{\rm (b)} If $f\in D(K)$, there exists a $\varphi\in C_b(K)$ with 
$\|f\|_{qD} = \|f-\varphi\|_D$. 
\endproclaim 

\demo{Remark} 
The proof of (b) yields an alternate proof of (10) in Theorem~3.2; the proof 
doesn't use partitions of unity. 
\enddemo 

\demo{Proof} 
(a) Let $x\in K$. If $\osc_\alpha f(x) <\infty$, then by the upper 
semi-continuity of $\osc_\alpha f$, given $\varep>0$, there is an open 
neighborhood $U$ of $x$ with $\osc_\alpha \mid U<\osc_\alpha f(x) +\varep$, 
and hence by Theorem~3.2, $f\in D(U)$; moreover $i_D (f\mid U)\le\alpha$, 
so also we have that 
$$\osc_\alpha f(x) \le \|f\mid U\|_{qD} = \|\osc_\alpha f\mid U\|_\infty 
\le \osc_\alpha f(x) +\varep
\leqno(22)$$ 
(where the equality follows by (10) of Theorem 3.2). 
Thus (a) follows when $\osc_\alpha f(x)<\infty$, since $\varep>0$ is arbitrary. 
Again if $f\in D(U)$ for some open neighborhood $U$ of $x$, we obtain that 
$\osc_\alpha f(x) \le \|f\mid U\|_{qD} <\infty$, and this establishes (a), 
for we also get that then $\osc_\alpha f(x)=\infty$ iff $f\notin D(U)$ 
for every open neighborhood $U$ of $x$. 
To prove (b), we first recall the standard result 
(the Hahn interposition theorem): 
{\it given $u,\ell$ upper, lower semi-continuous functions respectively 
on a metric space $K$, with $u\le \ell$, there exists a continuous 
function $\varphi$ on $K$ with $u\le\varphi\le \ell$.} 

Now let $f\in D(K)$, and $\beta = \|\osc_\alpha f\|_\infty$. 
It suffices to prove there exists a continuous $\varphi$ on $K$ with 
$$\osc_\alpha f-\beta \le \varphi-f \le\beta -\osc_\alpha f\ . 
\leqno(23)$$ 

Indeed, then $\varphi$ satisfies 
$$|f-\varphi| + \osc_\alpha f\le \beta 
\leqno(24)$$ 
Of course then $\varphi$ is bounded, since $f$ is, and moreover $i_D(f-
\varphi) = i_Df$, so by Theorem~3.2, since also $\osc_\alpha f=\osc_\alpha 
(f-\varphi)$, 
$$\|f-\varphi\|_D = \|\, |f-\varphi +\osc_\alpha (f-\varphi)\|_\infty 
\le\beta \ .
\leqno(25)$$ 

However it follows from (9) of 3.2 also that 
$\|\osc_\alpha f\|_\infty \le \|f\|_{qD}$, whence $\|f-\varphi\|_D = 
\|f\|_{qD}$. 

Now set $u= f+\osc_\alpha f-\beta$, 
$\ell = f-\osc_\alpha f+\beta$; then $u,\ell$ are upper, lower semi-continuous 
respectively, by Theorem~3.2, and of course $u\le\ell$, for this just 
says $\osc_\alpha f\le\beta$. 
Hence by the Hahn interposition theorem, there is a continuous $\varphi$ 
with $u\le \varphi\le\ell$; then $\varphi$ satisfies (23), completing 
the proof.\qed 
\enddemo 

We recall that for $f:K\to\real$ {\it bounded\/}; $\|\oosc f\|_\infty 
= \inf_{\varphi \in C_b(K)} \|f-\varphi\|_\infty$, and again the infimum 
is attained. 
Thus also $\oosc f(x)$ may be obtained as the local distance (at $x$) 
from $f$ to $C_b(K)$, just as we have done for $D(K)$. 
Now suppose $f\in D(K)$ and $i_Df=1$. 
Then we obtain that $\|\osc f\|_\infty = \inf_{\varphi\in C_b(K)} 
\|f-\varphi\|_D$, (and the infimum is attained). 
thus $\osc f$ in this case, is the appropriate measure for the 
$D$-distance to $C_b(K)$, while $\oosc f$ gives the measure for the 
sup-distance.  

We next give a basic tool for computing the finite oscillation functions. 
For example, this result, combined with Theorem~3.2, implies Lemma~1.8 for 
real functions; we shall see shortly that it also yields immediately 
(in combination with 3.2) that functions of finite Baire-index belong to $D$. 

\proclaim{Lemma 3.6} 
Let $f:K\to\complex$ be bounded, $n$ a positive integer, and $x\in K$ 
given. Then 
$$\osc_n f(x) = \sup \biggl\{ \sum_{i=1}^k \varep_i : 1\le k\le n\ ,\ 
\varep_i >0\ \text{ for all }\ 1\le i\le k\ ,\ \text{ and }\ 
x\in \os_k (f,(\varep_i))\biggr\}\ . 
\leqno(26)$$ 
\endproclaim 

\demo{Remark} 
We interpret the sup of the empty sum to be zero. 
Now let $W_n(x)$ be the term on the right side of the equality in (26). 
Then it is obvious that 
$$W_n(x)  = \sup \biggl\{ \sum_{i=1}^n \varep_i : \varep_i\ge0
\text{ for all $i$ and } x\in \os_n(f,(\varep_i))\biggr\}\ .$$ 
Indeed, suppose $n$, $(\varep_i)_{i=1}^n$ are given, and 
$i_1<\cdots < i_k$ are the indices $j$ with $\varep_j>0$. 
Then $\os_n(f,(\varep_i)) = \os_k (f,(\varep_j))$. 
\enddemo

\demo{Proof of 3.6} 
It is most convenient to prove the two relevant inequalities by induction 
on $n$. Let us first then show (by induction on $k$) that for $\varep_i>0$ 
all $1\le i\le k$, 
$$x\in \os_k(f,(\varep_i)) \text{ implies } \osc_k f(x)\ge 
\sum_{i=1}^k \varep_i\ .
\leqno(27)$$ 
Now this is trivial for $k=1$. 
Suppose proved for $k$, and let then $x\in \os_{k+1}(f,(\varep_i))$. 
Set $Y= \os_k (f,(\varep_i))$ and $g= f\mid Y$. 
Since then $\osc g(x) \ge \varep_{k+1}$, given $0<\varep <\varep_{k+1}$, 
we may choose $(x_m)$ in $Y$ with $x_m\to x$ and 
$$\uosc g(x_m) > \varep_{k+1}-\varep\text{ for all } m\ .
\leqno(28)$$
Now fixing $m$ we may choose $(y_j^m)$ in $Y$ so that $y_j^m\to x_m$ as 
$j\to\infty$ and 
$$|g(x_m) - g(y_j^m)| > \varep_{k+1}-\varep \text{ for all } j\ . 
\leqno(29)$$ 
Hence, 
$$\eqalign{
\widetilde{\osc}_{k+1} f(x_m) &\ge \olim\limits_{j\to\infty} 
|f(x_m)- f(y_j^m)| + \osc_k f(y_j^m)\cr
&\ge \varep_{k+1}-\varep + \sum_{i=1}^k \varep_i\text{ by (28) and the 
induction hypotheses.} \cr}$$
Thus $\osc_{k+1} f(x) \ge \olim_{m\to\infty} \widetilde{\osc}_{k+1} f(x_m) 
\ge \sum_{i=1}^{k+1} \varep_i-\varep$. 
Since $\varep>0$ is arbitrary, (27) follows, and thus we have for all $n$ 
and $x$ that $\osc_n f(x)\ge W_n(x)$ (where $W_n$ is defined in the above 
Remark). 
It remains to prove 
$$\osc_n f(x) \le W_n(x)\text{ for all } x\in K\ . 
\leqno(30)$$ 

Again, this is obvious for $n=1$, since we need only take $\varep = \osc f(x)$; 
then $x\in \os_1 (f,\varep)$. 
Suppose $n\ge1$ and the statement proved for $n$. 
Fix $x$ and choose $(x_j)$ with $x_j\to x$ and 
$$\widetilde{\osc}_{n+1} f(x_j) \to \osc_{n+1} f(x)\text{ as } j\to\infty\ . 
\leqno(31)$$ 

Then choose for each $j$, a sequence $(y_m^j)$ with 
$$\widetilde{\osc}_{n+1} f(x_j) = \lim_{m\to\infty} |f(x_j)- f(y_m^j)| 
+ \osc_n f(y_m^j)\ . 
\leqno(32)$$ 

We may further assume, by passing to subsequences of $(x_j)$, and then 
of $(y_m^j)_{m=1}^\infty$ if necessary, that 
$$\lim_{m\to\infty} |f(x_j) - f(y_m^j)| \dfeq \delta_j\text{ and } 
\lim_{m\to\infty} \osc_n f(y_m^j) \dfeq \lambda_j \text{ exist,}$$ 
and that moreover $\delta \dfeq \lim_{j\to\infty} \delta_j$ and 
$\lambda = \lim_{j\to\infty} \lambda_j$ exist. 
Thus 
$$\osc_{n+1} f(x) = \delta+\lambda\ .
\leqno(33)$$ 

Now if $\delta=0$, then $\osc_{n+1} f(x) = \osc_nf(x)$ and so 
$\osc_{n+1} f(x) \le W_n (x) \le W_{n+1}(x)$. 
Similarly,  if $\lambda=0$, $\osc_{n+1}f(x) = \osc f(x)$ and we are done. 
So we assume $\delta,\lambda >0$. 
Now let $\varep>0$, $\varep <\min \{\delta,\lambda\}$. 
By passing to further subsequences, we may now assume that  
$$\eqalign{
\hbox{(i)}\qquad &\delta_j >\delta -\varep \text{ and } \lambda_j>\lambda
-\varep\text{ for all } j\ ,\cr 
\hbox{(ii)}\qquad &\osc_n f(y_m^j) >\lambda-\varep\text{ for all $j$ and $m$.} 
\cr}\leqno(34)$$

It then follows by the induction hypothesis, that for each $j$ and $m$, 
we may choose non-negative sequences $(\varep_i^{(j,m)})_{i=1}^n$ with 
$y_m^j \in \osc_n (f,(\varep_i^{j,m}))$ and 
$$\sum_{i=1}^n \varep_i^{(j,m)} > \lambda-\varep\ . 
\leqno(35)$$ 

On the other hand, since we have that $2n\|f\|_\infty \ge \osc_n f(y_m^j) 
\ge \sum_{i=1}^n \varep_i^{(j,m)}$, the sequences $(\varep_i^{j,m})$ are 
uniformly bounded. 
It then follows, by passing to further subsequences of $(y_m^j)$ and 
$(x_j)$ if necessary, that we may assume for each $i$, that 
$$\eqalign{
\hbox{(i)}\qquad&\lim_{m\to\infty} \varep_i^{j,m} \dfeq  \varep_i^j
\text{ exists,}\cr
\hbox{(ii)}\qquad&\lim_{j\to\infty} \varep_i^j \dfeq \varep_i \text{ exists.}
\cr} 
\leqno(36)$$ 

We now have, from (35), that 
$$\sum_{i=1}^n \varep_i \ge \lambda-\varep > 0\ . 
\leqno(37)$$ 

Now let $i_1<\cdots < i_k$ be the indices $j$ with $\varep_j >0$, and let 
$0<\eta < \min_{1\le j\le k} \varep_{i_j}$. 
By finally passing again to further subsequences, we may at last assume 
(using (36)) that $\varep_{i_\ell}^j > \varep_{i_\ell}-\eta$ for all 
$\ell$ and all $j$, and finally that 
$$\varep_{i_\ell}^{j,m} > \varep_{i_\ell} -\eta \text{ for all } 
1\le \ell\le k\ ,\ \text{ all $j$, all $m$.}
\leqno(38)$$ 

Now let $\mu_\ell = \varep_{i_\ell} -\eta$, $1\le \ell\le k$. 
Then we have that 
$$y_m^j \in \osc_k (f,(\mu_\ell))\text{ for all $j$ and $m$.} 
\leqno(39)$$ 
Indeed, this follows from the following observation: 
Fix $j$ and $m$, and let $\boldsymbol\varepsilon_r^{j,m} = 0$ 
if $r\ne i_\ell$ any $\ell$; otherwise if $r=i_\ell$, let 
$\boldsymbol\varepsilon_r^{j,m} = \mu_\ell$; then 
since $\boldsymbol\varepsilon_r^{j,m} \le \varep_r^{j,m}$ for all 
$1\le r\le n$ by (38), 
$$\osc_n (f,(\varep_i^{j,m})) \subset 
\osc_n (f,\boldsymbol\varepsilon_r^{j,m}) = 
\osc_k (f,(\mu_\ell))\ .
\leqno(40)$$ 

Now since $\osc_k (f,(\mu_\ell))$ is closed, then fixing $j$, we have since
$x_j = \lim_{n\to\infty} y_m^j$ and (39) holds, that 
$$x\in \osc_k (f,(\mu_\ell))\ .
\leqno(41)$$ 
Also 
$$\sum_{\ell=1}^k \mu_\ell \ge \lambda -\varep-k\eta
\leqno(42)$$ 
(To see (42), 
$$\sum_{\ell=1}^k \mu_\ell = \sum_{\ell=1}^k (\varep_{i_\ell}-\eta) 
= \sum_{i=1}^n \varep_i -k\eta 
\ge \lambda -\varep -k\eta \text{ by (37)).}$$ 
At last, set $\mu_{k+1} = \delta-\varep$. 
Since $|f(x_j) - f(y_m^j)| >\delta-\varep$ for all $m$, and (39) holds, we 
have that 
$$x_j \in \os_{k+1} (f,(\mu_i))\text{ for all $j$\ .} 
\leqno(43)$$ 
Finally, since $\osc_{k+1} (f,(\mu_i))$ is closed, by (43) we have 
that also $x\in \osc_{k+1}(f,(\mu_i))$. 
We have by (42) and the definition of $\mu_{k+1}$ that 
$$\eqalign{ \sum_{i=1}^{k+1}\mu_i &\ge \lambda +\delta -2\varep-k\eta\cr 
&= \osc_{n+1} f(x) - 2\varep -k\eta\ \text{ (by (33)).}\cr}$$ 
But $\varep>0$, $\eta>0$ were arbitrary (and $k\le n$), so we have indeed 
proved (30) for the ``$n+1$'' case.\qed
\enddemo 

We now draw several simple consequences of the lemma. 
The first one gives an alternate formula for computing $\osc_w f$. 
(The norm assertion in its statement is given as Lemma~3 of the Introduction.) 

\proclaim{Corollary 3.7} 
Let $f:K\to\complex$ be a bounded function. 

{\rm (a)} For all $x\in K$, 
$$\eqalign{
\widetilde{\osc}_w f(x) &=\sup \biggl\{ \sum_{i=1}^k \varep_i : 1\le k< 
\infty\ ,\ \varep_i>0\text{ for all } 1\le i\le k\ ,\cr
&\qquad \text{and }\ x\in \os_k (f,(\varep_i))\biggr\}\ .\cr}$$
{\rm (b)} 
$$\|\osc_w f\|_\infty  = \sup \biggl\{ \sum_{i=1}^k \varep_i : 
1\le k<\infty\ ,\ \varep_i>0\text{ for all $i$ and } 
\os_k (f,(\varep_i)) \ne\emptyset\biggr\}\ .$$
\endproclaim 

\demo{Proof} 
(a) follows immediately from Lemma 3.6 and the definition: 
$\widetilde{\osc}_w f(x) = \sup_{k<\infty} \osc_k(f(x)$. 
(b) follows immediately from (a) and the definition: 
$\osc_w f= U\widetilde{\osc}_w f$, whence 
$\|\osc_w f\|_\infty = \|\widetilde{\osc}_w f\|_\infty$.\qed 
\enddemo 

The next result shows in particular that functions of finite Baire-index 
belong to $D$; this is proved by alternate methods in \cite{CMR} 
(with a little more work, we also recapture in Section~4, Corollary~4.11, 
the result of \cite{CMR} that such functions are in $SD$). 

\proclaim{Corollary 3.8} 
Let $f:K\to\complex$ be of finite Baire index. 
Then $f\in D(K)$ and $i_D(f)\le i_B(f)$. 
In particular if $n=i_B(f)$ and $f$ is real-valued, 
$$\|f\|_D \le (2n+1) \|f\|_\infty\text{ and } \|f\|_{qD} \le 
n\|\osc f\|_\infty \le 2n\|f\|_\infty\ .
\leqno(44)$$ 
\endproclaim 

\demo{Remark} 
We show below (as a simple exercise) that these estimates are best 
possible (for any $K$ with $K^{(m)}\ne \emptyset$ for all $m=1,2,\ldots$). 
\enddemo 

\demo{Proof} 
Let then $n$ be as above, and suppose $(\varep_i)_{i=1}^k$ given with 
$\varep_i>0$ for all $i$ and $\os_k (f,(\varep_i))\ne \emptyset$. 
Let $\varep = \min_{1\le i\le k} \varep_i$. 
Then $\os_k (f,\varep) \ne\emptyset$, hence by definition of $i_B(f)$, 
$k\le n$.  
It then follows by Lemma~3.6 that $\osc_n f=\osc_{n+1} f$, so $i_D(f)\le n$. 
Thus by Theorem~3.2, if $f$ is real, 
$$\|f\|_D = \|\, |f| + \osc_n f\|_\infty \ \text{ and }\ 
\|f\|_{qD} = \|\osc_n f\|_\infty\ . 
\leqno(45)$$ 
The estimates in (44) are now immediate.\qed  
\enddemo 

\proclaim{Corollary 3.9} 
If some finite derived set of $K$ is empty, then every bounded function 
on $K$ belongs to $D(K)$. 
\endproclaim 

\demo{Proof} 
Suppose $K^{(n)} \ne\emptyset$, $K^{(n+1)} = \emptyset$. 
Now fixing $f:K\to\complex$ bounded and $\varep>0$, it follows easily by 
induction that $\os_j (f,\varep) \subset K^{(j)}$ for all $j=1,2,\ldots$. 
Hence $i_B f \le n$, so $f\in D(K)$.\qed 
\enddemo 

We next illustrate some of these results by computing the finite 
oscillations and $D$-norms of some simple functions. 
We first give the case of characteristic  functions  of sets. 

\proclaim{Proposition 3.10} 
Let $A$ be a non-clopen set in $\D(K)$, and $n= i(A)$; 
set $f=\chix_A$. Then for $1\le m\le n$, 
$$\left\{ \eqalign{
&\osc_m f(x) = j\text{ for all } x\in\partial^j A\sim\partial^{j+1}A\ ,
\qquad 0\le j<m\cr 
&\osc_m f(x) =m\text{ for all } x\in\partial^m A\cr}\right.
\leqno(46)$$ 
\endproclaim

It follows immediately, by Corollary 3.8, that then also $n=i_D(f)$. 
Indeed, we have that $i_D(f) \le i_B(f)=n$ by 3.8; 
but (46) shows that $\osc_{n-1} f\ne \osc_n f$. 
Moreover this, together with Theorem~3.2, proves Theorem~2.2(a) and gives 
another proof of Theorem~2.2(b). 
Indeed we have by 3.2 that $\|f\|_{qD} = \|\osc_n f\|_\infty = n$, while 
$$\|f\|_D = \|\, |f| + \osc f\|_\infty = n\text{ if } 
A\cap \partial^n A=\emptyset\ ;\ = n+1\text{ if } A\cap \partial^n A
\ne\emptyset\ .$$

\demo{Proof of 3.10} 
We prove (46) by induction on $m$. 
For notational convenience, set $K^j = \partial^j A$ for all $j$. 
First note that since $f$ is $\{0,1\}$-valued, 
$\osc_m f\le m$ for all $m=0,1,2,\ldots.$ 
(In fact it follows easily by induction that $\osc_m g\le m\|g\|_\infty$ 
for any non-negative function $g$ and all $m<\infty$.) 
Since each point of $K^1 = \partial A$ is a cluster point of $A$ and $\sim A$, 
$\osc f\mid K^1 \ge 1$; 
hence $\osc f\mid K^1 =1$, and of course 
$\osc f\mid (K^0 \sim K^1) =0$ (recall that $\partial^0 A=K$ by definition). 

Suppose (46) established for $1\le m<n$. 
Now it follows that $i(A\mid \sim K^{m+1})$, the index of $A\cap \sim K^{m+1}$ 
relative to the metric space $\sim K^{m+1}$, equals $m$. 
Hence $(\osc_{m+1} f) \mid \sim K^{m+1} =\osc_{m+1} (f\mid \sim K^{m+1}) 
= (\osc_mf) \mid \sim K^{m+1}$ and this satisfies (46). 
Now if $x\in K^{m+1} \cap A$, since $K^{m+1} = \partial (A\mid K^m)$, we may 
choose $(x_j)$ in $K^m\sim A$ with $x_j \to x$. 
But then $\widetilde{\osc}_{m+1} f(x)\ge \lim_{j\to\infty} |f(x_j) -f(x)| 
+ \osc_m f(x_j) = m+1$. 
Similarly if $x\in K^{m+1} \sim A$, choose $(x_n)$ in $K^m\cap A$ with 
$x_n\to x$ to again obtain $\widetilde{\osc}_{m+1} f(x) \ge m+1$. 
Since $\osc_{m+1} f\le m+1$, we have thus established that 
$\osc_{m+1} f=m+1$ on $K^{m+1}$, completing the proof.\qed 
\enddemo 

By Proposition 2.6, if $K^{(m)} \ne\emptyset$ for all $m$, we thus obtain 
that for all finite integers $n$, there exists an $f$ in $D(K)$ with 
$i_D f = n$. 
To see that the estimates in Corollary~3.8 are best possible in such a $K$, 
fix $n$, choose sets $K = K^0\supset  K^1 \supset \cdots\supset K^n 
\supset K^{n+1}=\emptyset$ with $K^j$ closed non-empty nowhere dense in 
$K^{j-1}$ for all $j=1,2,\ldots,n$, and now define $f$ on $K$ by 
$$f(x) = (-1)^j\ \text{ for }\ x\in K^j\sim K^{j+1}\ ,\qquad 
j=0,1,\ldots,n\ . 
\leqno(47)$$ 
Of course $f= 2\chix_A -1$, where $A= (K^0\sim K^1)\cup (K^2\sim K^3) \cup 
\cdots$. 
Thus since $i(A) = n$,  $i_Bf=n$, and $\osc_mf = 2\osc_m \chix_A$ for all $m$. 
So by what we have already proved, $i_Df=n$ and moreover $\osc_n f=2n$ 
on $K^n$. 
Hence $|f| + \osc_n f= 2n+1$ on $K^n$, and trivially $|f| + \osc_n f\le 2n+1$. 
Thus by Theorem~3.2, 
$$\left\{ \eqalign{
&\|f\|_D = \|\, |f| + \osc_n f\|_\infty = 2n+1\cr 
&\|f\|_{qD} = \|\osc_n f\|_\infty = 2n\cr 
&\text{and }\ \|\osc f\|_\infty =1\ .\cr}\right.
\leqno(48)$$
Thus for this $f$, since $\|f\|_\infty=1$, the inequalities in (44) are 
all equalities.

For another example, let $J[0,1]$ denote the space of all real-valued 
bounded functions on $[0,1]$ with only jump-discontinuities; i.e., all 
functions $f$ so that $f(x+)$, $f(x-)$, the ``right and left limits,'' 
exist at each $x$. 
As is well known, $J[0,1]$ is a Banach algebra under the sup-norm; 
if $f\in J[0,1]$ and $\varep>0$, then $\os (f,\varep)$ is a finite set, 
hence $i_Bf\le1$ for all such $f$. 
Evidently then $i_Bf=1=i_Df$ iff $f\in J[0,1]\sim C[0,1]$. 
Thus $J[0,1]\subset D[0,1]$. 
(This result, with an alternate proof, is due jointly to F.~Chaatit and the 
author; cf.\ \cite{C}.) 
It then follows that $\|f\|_D \le 3\|f\|_\infty$ and $\|f\|_{qD}\le 
2\|f\|_\infty$ for all $f\in J[0,1]$. 
Finally, we have for any $f\in J[0,1]$, $x\in [0,1]$, that 
$$\leqalignno{
\osc f(x) &= \max \bigl\{ |f(x) - f(x-)|, |f(x) - f(x+)|\bigr\}\ , &(49)\cr
\oosc f(x) &= \max \bigl\{|f(x) -f(x-)|, |f(x)-f(x+)|, 
|f(x+) - f(x-)| \bigr\}\ .&(50)\cr}$$ 

Evidently if $f\in J[0,1]$ and is right-continuous, we have that 
$\oosc f = \osc f$ and hence $\|f\|_{qD} = 2\inf\{\|f-\varphi\|_\infty 
: \varphi \in C[0,1]\}$. 
(However if e.g., $f= \chix [0,\frac12) - \chix (\frac12,1]$, then 
$\|\osc f\|_\infty =1$ and $\|\oosc f\|_\infty =2$; then 
$f$ has the same distance from $C[0,1]$, in both the $D$ and 
sup-norms.)

We give one last example, computing the $D$-norms for a  
natural class of simple $D$-functions. 

\proclaim{Proposition 3.11} 
Let $n\ge1$ and $K= K_0\supset \cdots\supset K_n$ be non-empty closed 
subsets of $K$ with $K_i$ nowhere dense in $K_{i-1}$ for all $1\le i\le n$; 
set $K_{n+1}=\emptyset$. 
Let $a_0,\ldots,a_n$ be given real numbers, and $f:K\to\real$ be the 
function with $f\mid (K_i\sim K_{i+1})\equiv a_i$ for all $0\le i\le n$. 
Then 
$$ \|f\|_D = \sum_{i=1}^n |a_i-a_{i-1}| + |a_n|\ \text{ and }\ 
\|f\|_{qD} = \sum_{i=1}^n |a_i - a_{i-1}|\ . 
\leqno(51)$$
\endproclaim 

\demo{Proof} 
Let $1\le j\le n$. We shall prove by induction that 
$$\left\{ \eqalign{&\osc_j f \le \sum_{i=1}^j |a_i - a_{i-1}|\text{ on } 
\sim K_{j+1}\ ,\cr
&\text{ with equality holding on }\ K_j\sim K_{j+1}\ .\cr}\right. 
\leqno(52)$$ 

The proof for $j=1$ is rather evident, for $f$ is continuous on $K_0\sim K_1$, 
and if $x\in K_1\sim K_2$, then choosing $(x_k)$ in $K_0\sim K_1$ with 
$x_k\to x$ we have that $\uosc f(x) =\lim_{k\to\infty} |f(x_k)-f(x)| 
= |a_1-a_0|$. 
On the other hand, it's clear that $\uosc f(x) \le |a_1-a_0|$, which shows 
immediately that $\osc f(x) = U\uosc f(x) = |a_1-a_0|$ also. 

Now suppose $j<n$ and (52) is proved for $j$. 
We seek to prove this for $j+1$. 
Again if we let $x\in K_{j+1}\sim K_{j+2}$, we may choose $(x_k)$ in 
$K_j\sim K_{j+1}$ with $x_k\to x$, and then 
$$\eqalign{ \widetilde{\osc}_{j+1} f(x) 
&\ge \lim_{k\to\infty} |f(x_k) - f(x)| + \osc_j f(x_k)\cr
&= |a_{j+1} - a_j| + \sum_{i=1}^j |a_i -a_{i-1}|\ .\cr} 
\leqno(53)$$ 

To obtain the reverse inequality, let $\lambda = \|\osc_{j+1} f\mid 
\sim K_{j+2}\|_\infty$, and assume $\lambda >0$. 
We must show 
$$\lambda \le \sum_{i=1}^{j+1} |a_i-a_{i-1}|\ . 
\leqno(54)$$ 

Since $f$ is simple, there must be a smallest $i$, $0\le i\le j$, and an 
$x\notin K_{j+2}$ with $\lambda = \widetilde{\osc}_{i+1} f(x)$. 
Then 
$$\osc_i f(x) <\lambda\ .
\leqno(55)$$ 
Indeed otherwise, again since $f$ is simple, there would exist a $y$ with 
$\osc_i f(x) = \lambda = \widetilde{\osc}_i f(y)$, contradicting the 
definition of $i$. 
Now if $x\notin K_{j+1}$ then since $i_B f\mid \sim K_{j+1} \le j$, 
$i_D f\le j$, and so $\osc_i f(x) = \osc_j f(x)$, whence (54) holds by 
our induction hypothesis (52). 
Thus suppose $x\in K_{j+1}$ and choose a sequence $(x_k)$ in $K$ with 
$x_k\to x$ and 
$$\widetilde{\osc}_{i+1} f(x) = \lim_{k\to\infty} |f(x_k) - f(x)| + 
\osc_i f(x_k)\ . 
\leqno(56)$$ 
Now without loss of generality, by passing to a subsequence, we may assume 
there is an $r$, $0\le r\le j+1$, with 
$$x_k\in K_r\sim K_{r+1}\ \text{ for all }\ k\ .$$ 
In fact, $r=j+1$ is impossible, for then we obtain that 
$\widetilde{\osc}_{i+1} f(x) \le \osc_i f(x)$, contradicting (55). 
Again, since $i_D f\mid \sim K_{r+1} \le r$, we have that $\osc_i f = 
\osc_r f$ on $\sim K_{r+1}$, and thus by (52) and (56), 
$$\eqalign{ 
\widetilde{\osc}_{i+1} f(x) 
& \le |a_{j+1} - a_r| + \sum_{\ell=1}^r |a_\ell - a_{\ell-1}|\cr 
&\le \sum_{\ell=r+1}^{j+1} |a_\ell - a_{\ell-1}|  
+ \sum_{\ell=1}^r |a_\ell - a_{\ell-1}|\ ,\cr}$$ 
proving (54). 
Of course (53) and (54) establish (52) for $j+1$. 

The conclusion of 3.11 now follows immediately from Theorem~3.2, for 
$i_Df\le n$ and hence $\|f\|_{qD} = \|\osc_n f\| = 
\sum_{i=1}^n |a_i-a_{i-1}|$ by (52), while 
$\|f\|_D = \|\, |f| + \osc_n f\|_\infty 
= \sum_{i=1}^n |a_i - a_{i-1}| + |a_n|$. 
(Indeed, on $K_j \sim K_{j+1}$, $\osc_n f= \osc_j f$, so 
$|f| +\osc_n f|_{K_j\sim K_{j+1}} = 
\sum_{i=1}^j |a_i -a_{i-1}| + |a_j| 
\le \sum_{i=1}^n |a_i - a_{i-1}| + |a_n|$.)\qed 
\enddemo 

\demo{Remarks}  
1. Of course (52) holds for complex numbers $a_0,\ldots,a_n$ as well. 
Now assuming (as we may) that $a_i\ne a_{i+1}$ for all $0\le i\le n-1$, 
then $i_B f=n$. 
It can be shown that $i_Df=n$ if and only if $a_i\notin co\{a_{i-1},a_{i+1}\}$ 
for all $1\le i\le n-1$. 

2. It follows from Theorem 2.3 that if $f\in D(K)$ is real-valued and 
$K= \bigcup_{i=1}^m W_i$, the $W_i$'s  closed, then 
$\|f\|_{D(K)} = \max_{1\le i\le m} \| f\mid W_i\|_{D(W_i)}$. 
(Cf.\ Lemma~4.21 and the remark following its proof.) 
Call a function satisfying the hypotheses of 3.11 a cell. 
It can be shown that if $f$ is a simple $D$-function, then there exist 
closed non-empty sets $W_1,\ldots, W_m$ with 
$K= \bigcup_{i=1}^m W_i$ and $f\mid W_i$ a cell, for all $i$. 
Thus in theory, one can compute the $D$-norm of an arbitrary simple 
$D$-function, using this fact and 3.11.
\enddemo 

\head \S4. Strong $D$-functions.\endhead 

We begin with some natural examples of strong $D$-functions, needed in 
the sequel. 
Our first result is also shown in \cite{CMR}; we give it again here, 
for completeness. 
(As before, $K$ denotes a given metric space.) 

\proclaim{Proposition 4.1} 
Every bounded continuous function on  $K$ is a strong $D$-function.
\endproclaim 

\demo{Proof} 
Obviously it suffices to prove that continuous bounded real-valued functions 
$f$ are in $SD(K)$. 
Let $f$ be such a function, and suppose without loss of generality that 
$\|f\|_\infty \le 1$. 
Let $\varep>0$. We shall show there is an upper semi-continuous simple 
function $\varphi$ with 
$$0\le f-\varphi <\varep\ . 
\leqno(1)$$ 
It follows, since then $f-\varphi$ is non-negative lower semi-continuous, 
that $\|f-\varphi\|_D= \|f-\varphi\|_\infty\le \varep$, proving our result. 

Choose $n$ with $\frac1n \le\varep$, let $-n\le j\le n$, and set $A_j = 
\{\omega\in K: \frac{j}n \le f(\omega) < \frac{j+1}n\}$. 
Now define $\varphi$ by 
$$\varphi = \sum_{j=-n}^n {j\over n} \chix_{A_j}\ . 
\leqno(2)$$ 
It is then trivial that $\varphi$ is a simple $D$-function and that (1) holds, 
so we need only verify that $\varphi$ is upper semi-continuous. 
Let then $x\in K$, and $(x_j)$ a sequence in $K$ with $x_j\to x$ such that 
$\lambda \dfeq \lim_j \varphi (x_j)$ exists. 
We must show that 
$$\lambda \le \varphi(x)\ . 
\leqno(3)$$ 

Since $(A_j)_{j=-n}^n$ is a partition of $K$, by passing to a subsequence, 
we may assume that there is a $j$ so that $x_k\in A_j$ for all $k$. 
But then for all $k$, $\varphi (x_k) = \frac{j}n$  and since 
$\frac{j}n \le f(x_k) < \frac{j+1}n$ for all $k$, 
$${j\over n} \le f(x) \le {j+1\over n}\ \text{ by the continuity of }\ f\ .
\leqno(4)$$ 
Thus $\lambda =\frac{j}n$ and $\varphi (x) = \frac{j}n$ or 
$\frac{j+1}n$, so (3) holds.\qed 
\enddemo 

We next give a useful class of functions of finite-index, containing $S(K)$,  
the space of simple $D$ functions on $K$.  

\proclaim{Proposition 4.2} 
Let $n\ge1$, and $K= K_0 \supset K_1 \supset\cdots\supset K_n\supset K_{n+1} 
= \emptyset$ be closed subsets of $K$, with $K_n\ne \emptyset$. 
Let $f:K\to \complex$ be such that $f\mid (K_i\sim K_{i+1})$ belongs to 
$C_b(K_i\sim K_{i+1})$ for all $i$. 
Then $f$ belongs to $SD(K)$ and $i_B f\le n$.
\endproclaim 

\demo{Proof} 
It is worth noting first that if $W$ is a DCS, then for any $g:W\to\complex$,  
$$g\in SD(W)\text{ if and only if } g\chix_W \in SD(K)\ . 
\leqno(5)$$ 
(Recall that $g\cdot\chix_W =0$ off $W$; $=g$ on $W$.) 

Indeed, it is evident that $g$ is a simple $D$-function on $W$ if and only 
if $g\chix_W$ is a simple $D$-function on $K$. 
Thus if $g\cdot \chix_W\in SD(K)$, choose $(f_n)$  simple $D$-functions  
on $K$ with $f_n\to g\cdot\chix_W$ in $D$-norm; then evidently 
$f_n\mid W\to g$ in $D(W)$-norm. 
Conversely, if $(f_n)$ is a sequence of simple $D$-functions on $W$ and 
$f_n\to g$ in $D(W)$, we have that $f_n\chix_W\to g\chix_W$ in $D(K)$, 
since $\|f_n\chix_W-g\chix_W\|_{D(K)} \le 2\|f_n-g\|_{D(W)}$ for all $n$ 
by Proposition~1.6. 

Now let $f$ be as in 4.2. 
We thus have by the preceding result that $f\in SD(K)$, since letting 
$\varphi_i = f\mid K_i \sim K_{i+1}$, then  $\varphi_i\in SD(K_i\sim 
K_{i+1})$, and $f= \sum_{i=0}^n \varphi_i \chix_{K_i\sim K_{i+1}}$. 

Now to show the index assertion, let $\varep >0$ be given. 
We then have by induction that 
$$\os_j (f,\varep) \subset K_j\ \text{ for all }\ j\le n\ . 
\leqno(6)$$ 
Indeed, the assertion is trivial for $j=0$. 
Suppose proved for $j<n$. 
But then $\os_{j+1}(f,\varep)= \os f\mid L,\varep)$ where 
$L= \os_j(f,\varep)$. Since $f$ is continuous on $K_j\sim K_{j+1}$, it 
is continuous on $L\sim K_{j+1}$, whence $\os_{j+1}(f,\varep) \subset 
K_{j+1}$. 

Of course (6) yields that $\os_{n+1}(f,\varep)=\emptyset$, since $f$ is 
assumed continuous on $K_n$. 
This proves that $i_B(f)\le n$.\qed
\enddemo 

\demo{Remark} 
It is evident that every simple $D$-function satisfies the hypotheses of 4.2. 
Indeed, let $f$ be such a function,  let $\lambda_1,\ldots,\lambda_k$ be the 
distinct values of $f$, and set $\varep = \min \{|\lambda_i - \lambda_j|
: i\ne j$, $1\le i,j\le k\}$. 
Then it follows (as noted in Section~1) that if $L\subset K$, $x\in L$, 
and $\osc f\mid L(x)<\varep$, $f$ is continuous at $x$. 
Now let $n=i_B (f,\varep)$, and let $K_j = \os_j(f,\varep)$ for 
$j=1,2,\ldots$. 
Then $K_{n+1}=\emptyset$ and $f\mid K_j \sim K_{j+1}$ is continuous for all 
$0\le j\le n$. 
Moreover then $i_B f=n$. 
\enddemo 

We next prove the characterization of $B_{1/4}$ given in Theorem~2 of the 
Introduction. 
Thanks to Lemma~3.6, this follows from the following result.  

\proclaim{Theorem 4.3} 
Let $f:K\to\complex$ be a given function. 
Then the following are equivalent. 

\iitem{\rm (a)} $f\in B_{1/4}(K)$. 
\iitem{\rm (b)} There exists a sequence $(\varphi_n)$ of simple $D$-functions 
with $\varphi_n\to f$ uniformly and $\sup \|\varphi_n\|_D<\infty$. 
\iitem{\rm (c)} $\osc_\omega f$ is bounded. 

Moreover when this occurs and $f$ is real-valued, 
$$\tfrac12 (\|f\|_\infty + \|\osc_\omega f\|_\infty) 
\le \|f\|_{B_{1/4}} \le \|f\|_\infty + 3\|\osc_\omega f\|_\infty \ .
\leqno(7)$$ 
\endproclaim 

To obtain Theorem 2, of the Introduction, simply note that by 
Corollary~3.7, $\beta= \|\osc_\omega f\|_\infty$, where $\beta$ is defined 
in the statement of Theorem~2.
We also note that for  $f\in B_{1/4}$, $f$ complex-valued 
$$\|\, |f| + \widetilde{\osc}_\omega f\|_\infty \le \|f\|_{B_{1/4}}\ . 
\leqno(8)$$ 
This follows directly from Lemma 1.8; hence the first inequality in (7) 
also holds for complex-valued functions. 
The argument below does not use Lemma~1.8, however. 

\demo{Remark} 
After writing the first draft of this paper, we learned of the following 
remarkable result of V.~Farmaki and A.~Louveau \cite{FL}. 

{\it If $f$ is a real-valued function on $K$, then} 
$$\|f\|_{B_{1/4}} = \big\|\,|f| + \widetilde\osc_\omega f\big\|_\infty\ .$$ 
(We  obtain this identity for strong $D$-functions $f$ in Corollary~4.6 
below). 
The proof in \cite{FL} is rather different than the argument for 4.3 given 
below; also our argument yields that if $\osc_\omega f$ is bounded (with 
$f$ real), then there exists a sequence $(\varphi_n)$ of simple $D$-functions 
with $\varphi_n\to f$ uniformly and $\|\varphi_n\|_D \le \|f\|_\infty 
+ 3\|\osc_\omega f\|_\infty$ for all $n$. 
Although this estimate is probably not optimal, it seems unlikely that one 
could choose such a sequence $(\varphi_n)$ with $\|\varphi_n\|_D 
\le \|\,|f| + \widetilde \osc_\omega f\|_\infty$ for all $n$. 
\enddemo 

\demo{Proof of 4.3} 
We may obviously assume that $f$ is real-valued. 
To see (8), suppose first $f$ is in $B_{1/4}$ and let 
$\lambda=\|f\|_{B_{1/4}}$. 
Let $n$ be a positive integer, $\varep>0$, and choose $\varphi$ in $D(K)$ with 
$$\|\varphi\|_D < \lambda +\varep\ \text{ and }\ 
\|\varphi-f\|_\infty  < {\varep\over n}\ . 
\leqno(9)$$ 

Then applying Proposition 3.1(b), 
$$\osc_n f - \osc_n \varphi \le \osc_n (f-\varphi) \le 2n\|\varphi-f\|_\infty 
< 2\varep\ . 
\leqno(10)$$ 
Thus 
$$\eqalign{ |f| + \osc_n f 
&\le |\varphi| + \osc_n\varphi +3\varep\ \text{ by (9), (10)}\cr 
&\le \|\varphi\|_D + 3\varep\ \text{ by Theorem 3.2}\cr 
&< \lambda +4 \varep \ \text{ by (9).}\cr}$$
Since $\varep>0$ is arbitrary and $|f|+\widetilde{\osc}_\omega f = 
\sup_n |f| +\osc_n f$, (8) now follows, and of course (8) yields the 
first inequality in (7). 

Suppose conversely that $\osc_\omega f$ is bounded and let $\mu = 
\|\osc_\omega f\|_\infty$. 
Now fix $\varep >0$ and set $n= i_B(f,\varep)$. Then 
$$n\varep \le \mu\ \text{ by Corollary 3.9.} 
\leqno(11)$$ 
Now let $K^j = \osc_j (f,\varep)$ for $j=0,1,2,\ldots$; thus 
$K^n\ne \emptyset$, $K^{n+1} = \emptyset$. 
By Proposition~1.17, we may choose a function $\varphi :K\to\real$ so 
that for all $j$, $0\le j\le n$, 
$$\varphi \mid (K^j\sim K^{j+1})\ \text{ is continuous and }\ 
|\varphi (x) - f(x)| <\varep\text{ for } x\in K^j\sim K^{j+1}\ . 
\leqno(12)$$ 

Evidently we thus have 
$$\|\varphi-f\|_\infty \le \varep\ .
\leqno(13)$$ 

Now by Proposition 4.2, $\varphi \in SD(K)$ and moreover $i_B(\varphi)\le n$, 
so by Corollary~3.8, $i_D(\varphi) \le n$, and thus by Theorem~3.2,  
$$\|\varphi\|_D = \|\, |\varphi| + \osc_n \varphi\|_\infty\ . 
\leqno(14)$$ 

Now we have that 
$$\leqalignno{
\osc_n\varphi \le \osc_n f + \osc_n(f-\varphi) 
&\le \mu +2n\varep\ \text{ by (13)}
&(15)\cr
&\le 3\mu\ \text{ by (11).}\cr}$$ 

Thus (13), (14) and (15) yield 
$$\|\varphi\|_D \le \|f\|_\infty + 3\|\osc_\omega f\|_\infty +\varep\ . 
\leqno(16)$$ 

Applying (13) and (16) for arbitrary $\varep$, we have thus established 
the existence of a sequence $(\varphi_j)$ in $SD(K)$ with $\varphi_j\to f$ 
uniformly and 
$$\olim\limits_{j\to \infty} \|\varphi_j \|_D \le \|f\|_\infty 
+ 3\|\osc_\omega f\|_\infty\ .
\leqno(17)$$ 
Of course this proves $f\in B_{1/4}$, and moreover yields the right hand side 
of (7). 
A simple density argument yields that in fact we may choose the $\varphi_j$'s 
to be simple $D$-functions, thus yielding (b) and completing the proof.\qed 
\enddemo 

\demo{Remark}
Define the quotient $B_{1/4}$-semi-norm, $\|\cdot\|_{qB_{1/4}}$, by 
$\|f\|_{qB_{1/4}} = \inf\{ \|f-\varphi\|_{B_{1/4}} : \varphi \in C_b(K)\}$. 
We then easily obtain that for $f\in B_{1/4}$, $\|f\|_{qB_{1/4}}$ is 
equivalent to $\|\osc_\omega f\|_\infty$. 
Indeed, the proof of Theorem~4.3 yields that for real $f$, $\|\osc_\omega f
\|_\infty \le \|f\|_{qB_{1/4}}$. 
On the other hand, we have by (7) that 
$$\eqalign{\|f\|_{qB_{1/4}} &\le \inf_{\varphi\in C_b(K)} 
\|f-\varphi\|_\infty + 3\|\osc_\omega f\|_\infty\cr 
&\le \|\osc f\|_\infty + 3\|\osc_\omega f\|_\infty 
\le 4\osc_\omega \|f\|_\infty\ .\cr}$$ 
That is, we have 
$$\|\osc_\omega f\|_\infty \le \|f\|_{qB_{1/4}} 
\le 4\|\osc_\omega f\|_\infty\ .$$ 
\enddemo 

We next give some useful oscillation invariants for $SD(K)$. 

\proclaim{Proposition 4.4} 
Let $f\in SD(K)$. 
Then $(\osc_nf)_{n=1}^\infty$ converges uniformly to $\osc_\omega f$.
\endproclaim 

We delay the proof, to draw some immediate consequences. 

\proclaim{Corollary 4.5} 
Let $f\in SD(K)$. Then 

\iitem{\rm (a)} $i_D f\le\omega$

\noindent and 

\iitem{\rm (b)} $\widetilde{\osc}_\omega f= \osc_\omega f$. 
\endproclaim 

\demo{Proof} 

(a) Let $\varep>0$ and choose $n$ with 
$$\osc_\omega f \le \osc_n f+\varep\ . 
\leqno(18)$$ 

It follows that fixing $x\in K$, then 
$$\widetilde{\osc}_{\omega+1} f(x) \le \widetilde{\osc}_{n+1} f(x)+\varep\ . 
\leqno(19)$$ 
Of course (19) yields that $\widetilde{\osc}_{\omega+1} f\le 
\osc_\omega f+\varep$; since $\varep>0$ is arbitrary, 
$\widetilde{\osc}_{\omega+1} f= \osc_\omega f\To \osc_{\omega+1} f 
= \osc_\omega f \To i_D f\le \omega$. 

(b) This is immediate from 4.4, since 
$\widetilde{\osc}_\omega f = \sup_n \osc_n f = \lim_{n\to\infty} \osc_n f$ 
point-wise. 

We shall use Corollary 4.5 later on, to construct some simple examples 
of functions in $D(K)\sim SD(K)$ for suitable $K$. 
The fact that $D(K)\sim SD(K)$ is non-empty in general, is obtained by 
different arguments in \cite{CMR}. 
Now it follows also by the results in \cite{HOR} that the $\|\cdot\|_D$ 
and $\|\cdot\|_{B_{1/4}}$ are not equivalent on $D$, in general. 
This also produces functions in $D\sim SD$, by the following result.
\enddemo 

\proclaim{Corollary 4.6} 
Let $f\in SD(K)$, $f$ real-valued. Then $\|f\|_D = \|f\|_{B_{1/4}}$. 
\endproclaim 

\demo{Proof} 
We have that $\|f\|_{B_{1/4}} \le \|f\|_D$ by definition. 
On the other hand, by the previous corollary and Theorem~3.2, 
$$\eqalign{\|f\|_D & = \|\, |f| +\osc_\omega f\|_\infty 
\ \text{ since }\ i_D f\le \omega\cr 
&= \|\, |f|+ \widetilde{\osc}_\omega f\|_\infty \le \|f\|_{B_{1/4}}\ ,\cr}$$ 
the last equality holding by the (8) (as shown in the proof 
of Theorem~4.3).\qed
\enddemo 

To prove Proposition 4.4, we first note that if $\alpha$ is any 
ordinal and $(f_n)$, $f$ are in $D$ with $f_n\to f$ in $D(K)$, then 
also $\osc_\alpha f\to \osc_\alpha f$ uniformly. 
Indeed, this follows immediately from the following simple result. 

\proclaim{Lemma 4.7} 
Let $f,g$ belong to $D(K)$, and $\alpha$ be a given ordinal. 
Then $\|\osc_\alpha f-\osc_\alpha g \|_\infty \le \|f-g\|_D$ if $f,g$ 
are real-valued, while $\|\osc_\alpha f-\osc_\alpha g\|_\infty 
\le 2\|f-g\|_D$ in general.
\endproclaim 

\demo{Proof} 
By Proposition 3.1(b), $\osc_\alpha f\le \osc_\alpha g + \osc_\alpha (f-g)$ 
and so $\osc_\alpha g\le \osc_\alpha f+ \osc_\alpha (g-f)$, 
also $\osc_\alpha (f-g) = \osc_\alpha (g-f)$, whence 
$$|\osc_\alpha f-\osc_\alpha g | \le \osc_\alpha (f-g)\ .
\leqno(20)$$ 
Lemma 4.7 now follows immediately from Theorem 3.2, since 
$\osc_\alpha\varphi \le \|\varphi\|_D$ for real-valued $\varphi$, 
$\osc_\alpha\varphi \le 2\|\varphi\|_D$ for complex-valued $\varphi$ 
(cf.\ (22) in Section~3).\qed 
\enddemo 

\demo{Proof of Proposition 4.4} 

Let $\varep>0$ and choose $\varphi$ a simple $D$-function with 
$$\|\varphi -f\|_D <\varep\ . 
\leqno(21)$$ 

Then by Lemma 4.7, 
$$\|\osc_\alpha \varphi - \osc_\alpha f\|_\infty < 2\varep
\ \text{ for any ordinal }\ \alpha\ . 
\leqno(22)$$ 

As noted in the remark following Proposition 4.2, $\varphi$ is of 
finite Baire-index. 
Thus if $n\ge i_B\varphi$, we have that $\osc_n\varphi = \osc_\omega \varphi$. 
Applying (22) for $\alpha=n$, $\alpha=\omega$, we obtain via the triangle 
inequality that 
$$\|\osc_n f-\osc_\omega f \|_\infty \le 4\varep\ . 
\leqno(23)$$ 
This proves 4.4.\qed 
\enddemo 

The next result yields that $SD$ is the span of its semi-continuous members. 
(This is Theorem~5a of the Introduction.) 
The proof uses the quantitative information in Theorem~3.2. 

\proclaim{Proposition 4.8} 
Let $f\in SD(K)$, $\varep>0$, $f$ real-valued. 
There exist non-negative $u,v$ lower semi-continuous functions belonging 
to $SD$ so that 
$$f= u-v\ \text{ and }\ \|u+v\|_\infty \le \|f\|_D +\varep\ . 
\leqno(24)$$ 
\endproclaim 

We first require the corresponding result for simple functions. 

\proclaim{Lemma 4.9} 
Let $f$ be a simple $D$-function. 

{\rm (a)} $\osc_nf$ is simple for all $n=0,1,2,\ldots$. 

{\rm (b)} If $f$ is real-valued, there exist simple non-negative 
lower semi-continuous functions $u$ and $v$ with 
$$f=u-v\ \text{ and }\ \|u+v\|_\infty = \|f\|_D\ . 
\leqno(25)$$ 
\endproclaim 

\demo{Proof} 

(a): Let $\lambda_1,\ldots,\lambda_k$ be the distinct values of $f$. 
Let $W= \{|\lambda_i - \lambda_j| : 1\le i,j\le k\}$. 
Let $L_j = \{w_1+\cdots+ w_r :1\le r\le j$ and $w_i\in W$ for all $i\}$ 
if $j\ge1$; let $L_0 = \{0\}$. 
Evidently $L_j$ is a finite set, for all $j$. 
We then have that 
$$\osc_j f\ \text{ is valued in $L_j$ for all $j$.} 
\leqno(26)$$ 
This is trivial for $j=0$; suppose the result proved for $j$. 
Let $x\in K$, and choose $(x_j)$ in $K$, $x_j\to x$, with 
$$\widetilde{\osc}_{j+1} f(x) = \lim_{n\to\infty} |f(x_n)-f(x)| 
+ \osc_j f(x_n)\ .$$ 
By passing to a subsequence, we may chose $w\in W$ and $z\in L_j$ 
so that $|f(x_n)-f(x)| = w$ and $\osc_j f(x_n) = z$ for all $n$. 
Evidently then $\widetilde{\osc}_{j+1} f(x) = w+z$, and this belongs 
to $L_{j+1}$. 
It is now evident that also $U\widetilde{\osc}_{j+1} f= \osc_{j+1} f$ is 
valued in $L_{j+1}$. 
Hence (26) holds and thus (a) is proved. 

(b): Since $f$ is of finite Baire index, there is an $n<\infty$ with 
$i_D f=n$. Thus by Theorem~3.2, letting $\lambda =\|f\|_D$, we have that 
$f=u-v$ and $\|u+v\|_\infty = \lambda$, where 
$u= {\lambda+f-\osc_n f\over2}$, $v= {\lambda-f-\osc_n f \over 2}$, 
and $u,v$ are lower semi-continuous non-negative. 
Now $u$ and $v$ are simple functions by part~(a), proving (b).\qed 
\enddemo 

\demo{Proof of Proposition 4.8} 
Recall that $S(K)$ denotes the family of simple $D$-functions on $K$. 
Since $SD(K) = \overline{S(K)}$ by definition, a standard density argument 
shows that given $\varep>0$ and $f\in SD(K)$, we may choose $(f_n)$ in 
$S(K)$ with 
$$\sum \|f_n\|_D < \|f\|_D + \varep\ \text{ and }\ 
f= \sum f_n\  
\leqno(27)$$ 
(where the series in (27) converges in $D(K)$). 
By Lemma~4.9(b), for each $n$ we may choose $u_n,v_n\ge0$ simple lower 
semi-continuous with 
$$f_n = u_n - v_n\ \text{ and }\ \|f_n\|_D = \|u_n+v_n\|_\infty\ . 
\leqno(28)$$ 
Now set $u=\sum u_n$ and $v=\sum v_n$. 
Since the series $\sum u_n$ and $\sum v_n$ converge uniformly, $u,v$ are 
non-negative lower semi-continuous, and of course $f=u-v$, 
and for any $x\in K$,  
$$u(x) + v(x) = \sum (u_n+v_n)(x) \le \sum \|f_n\|_D 
< \|f\|_D + \varep \ \text{ by (43).}$$ 
Hence $\|u+v\|_\infty \le \|f\|_D +\varep$. 
Finally, we have that $\sum u_n$, $\sum v_n$ converge to $u,v$ respectively, 
in the $D$-norm. 
Indeed, fixing $n$, then  $u-u_n$ is a non-negative lower semi-continuous 
function, since this is the uniform limit of $\sum_{j=n+1}^m u_j$ as 
$m\to\infty$. 
But then $\|u-u_n\|_D = \|u-u_n\|_\infty 
\le \sum_{j=n+1}^\infty \|u_j\|_\infty \to 0$ as 
$n\to \infty$ (by (27) and (28)); 
the argument for $v$ is identical. 
Thus since $\sum_{j=1}^n u_j$ and $\sum_{j=1}^n v_j$ are simple  
for all $n$, $u,v$ belong  to $SD$.\qed 
\enddemo 

\demo{Remark} 
If $f\in SD(K)$ and $\osc_\omega f$ is also in $SD(K)$, then Theorem~3.2 
and Corollary~4.5a yield Proposition~4.8, with in fact the functions $u,v$ 
in its statement chosen with $\|f\|_D = \|u+v\|_\infty$. 
However it can be seen that for any compact metric space $K$ with 
$K^{(\omega)} \ne\emptyset$, there exists an $f:K\to\real$ with 
$i_B(f)=1$, yet $\osc_f = \uosc f= \oosc f$ ($=\osc_\omega f$) 
is {\it not\/} strong-$D$.
\enddemo 

We next assemble some tools to prove that $SD$ is a complex lattice. 
We require the following structural lemma, which is obtained in 
\cite{CMR}. 

\proclaim{Lemma 4.10} 
Let $f$ and $g$ belong to $B_{1/4}(K)$, $\varep>0$. 
Then $i_B(f+g,\varep) \le i_B(f,\frac{\varep}2) + i_B(g,\frac{\varep}2)$. 
\endproclaim 

(We establish a generalization of this result later, in Lemma~4.20, 
in order to characterize $SD$ intrinsically.) 
It follows immediately from 4.10 that 
{\it if $f$ and $g$ are of finite Baire index, so is $f+g$ and} 
$$i_B(f+g) \le i_B(f) + i_B(g)\ . 
\leqno(29)$$ 

Proposition 4.2, together with (29), easily yields that functions of 
finite index are strong $D$. 
This result is obtained in \cite{CMR} by (29) and other methods. 

\proclaim{Corollary 4.11} 
Every function of finite Baire index belongs to $SD$. 
\endproclaim 

\demo{Proof} 
Assume that $f$ is real valued on $K$, of finite Baire index, 
let $n=i_B(f)$, and let $\varep>0$. 
Setting $K^j = \os_j (f,\varep)$ for all $j$, then $K^{n+1} =\emptyset$, 
and by Proposition~1.17, we may choose a function $\varphi$ on $K$ so 
that for all $0\le j\le n$, 
$$\varphi \mid (K_j\sim K_{j+1})\ \text{ is continuous and } 
|\varphi-f|<\varep \text{ on } K_j\sim K_{j+1}\ . 
\leqno(30)$$ 
By Proposition 4.2, we have that $\varphi \in SD(K)$ and $i_B(\varphi)\le n$. 
Hence by Lemma~4.10, applying (29), $i_B(\varphi-f)\le 2n$. 
Thus by Corollary~3.8, $i_D(\varphi-f)\le2n$, and thus  by (3.44), 
$\|\varphi-f\|_D \le (4n+1)\varep$. 
Since $\varep>0$ is arbitrary, the result is proved.\qed 
\enddemo 

We next recall a class of functions containing $SD$, which we need. 

\demo{Definition} 
$B_{1/2}^0 (K)$ denotes the family of all bounded functions $f:K\to\complex$ 
so that 
$$\lim_{\varep\to0} \varep i_B (f,\varep) = 0\ . 
\leqno(31)$$ 
\enddemo 

We need the following result, established in \cite{CMR} and reproved here for 
completeness. 

\proclaim{Proposition 4.12} 

\iitem{\rm (a)} $SD(K)\subset B_{1/2}^0 (K)$. 

\iitem{\rm (b)} If $f\in B_{1/2}^0 (K)$ and $f$ is semi-continuous, then 
$f\in SD(K)$. 

\iitem{\rm (c)} $B_{1/2}^0 (K)$ is a linear space. Moreover 
$|f|\in B_{1/2}^0 (K)$ provided $f\in B_{1/2}^0(K)$. 
\endproclaim 

\demo{Proof} 
We first show (c). 
If $f,g\in B_{1/2}^0 (K)$, then using Lemma~4.10, 
$$\lim_{\varep\to0} \varep i_B (f+g,\varep) 
\le 2\lim_{\delta\to0} \delta i_B(f,\delta) + 
2\lim_{\delta\to0} \delta i_B(g,\delta) =0\ .$$ 
Thus $f+g\in B_{1/2}^0 (K)$. 
If $\lambda$ is a non-zero scalar and $\varep>0$, then by induction 
we see that $\os_j (\lambda f,\varep) = \os_j (f,{\varep\over |\lambda|})$ 
for all $j$, hence $i_B (\lambda f,\varep) = i_B (f,{\varep\over|\lambda|})$, 
and so $\lim_{\varep\to0} \varep i_B (\lambda f,\varep) 
= |\lambda| \lim_{\delta\to0} \delta i_B(f,\delta)=0$. 
Finally, $\os_j (|f|,\varep) \subset \os_j (f,\varep)$ for all $j$, 
hence $i_B (|f|,\varep)\le i_B (f,\varep)$, whence $|f| \in 
B_{1/2}^0 (K)$, proving (c). 
To prove (a), let $f\in SD(K)$, assume without loss of generality 
that $f$ is real, let $\eta>0$, and choose $g$ a simple $D$-function 
with $\|f-g\|_D<\eta$. 
It then follows by Lemma~1.8 (or Theorem~3.2 and Lemma~3.6) that 
$$\varep i_B (f-g,\varep) <\eta\ \text{ for all }\ \varep>0\ .
\leqno(32)$$ 
Since $g$ is a simple $D$-function, $g$ has finite index; say $\mu = i_B(g)$. 
Then by Lemma~4.10, for any $\varep>0$, 
$$\eqalign{\varep i_B (f,\varep) 
& \le \varep i_B \Bigl(f-g ,{\varep\over2}\Bigr) + \varep i_B
\Bigl( g,{\varep\over2}\Bigr)\cr 
&\le 2\eta +\varep\mu\ \text{ by (32) and the definition of $\mu$.}\cr}$$ 
Hence $\olim_{\varep\to0} \varep i_B(f,\varep) \le 2\eta$. 
Since $\eta >0$ is arbitrary, (31) holds. 

Finally, to prove (b), suppose without loss of generality that $f$ is 
upper semi-continuous, let $\eta >0$, and choose $\varep>0$ so that 
$$\varep i_B (f,\varep) < \eta\ \text{ (with $\varep<\eta$).} 
\leqno(33)$$ 

Let then $n=i_B(f,\varep)$ and set $K^j = \os_j (f,\varep)$ for all $j$. 
Thus $K^n \ne\emptyset$, $K^{n+1} = \emptyset$; since for all $j$, 
$\osc f\mid (K^j\sim K^{j+1}) <\varep$, we may choose for each $j$ a 
continuous function $\varphi_j$ on $K^j\sim K^{j+1}$ with 
$|\varphi_j -f|<\varep$ on $K^j\sim K^{j+1}$. 

Now set $g=\sum_{j=0}^n \varphi_j \chix_{K^j\sim K^{j+1}}$. 
By Proposition~4.2, $g\in SD(K)$. 
Fixing $j$ and letting $W= K^j\sim K^{j+1}$, then evidently $f-g$ is 
upper semi-continuous on $W$; hence 
$$\|(f-g)\mid W\|_{D(W)} \le 3\|f-g\|_\infty \le 3\varep\ . 
\leqno(34)$$ 
Thus by Corollary 1.8, 
$$\|(f-g)\chix_W\|_{D(K)} \le 6\varep\ . 
\leqno(35)$$ 
Thus 
$$\eqalign{ \|f-g\|_D 
& = \Big\|\sum_{j=0}^n (f-g) \chix_{K^j\sim K^{j+1}}\Big\|_D\cr 
&\le \sum_{j=0}^n \|(f-g) \chix_{K^j\sim K^{j+1}} \|_D\cr 
&\le 6n\varep +6 \varep\cr 
&< 7\eta\ \text{ by (33).}\cr}$$ 
Since $\eta>0$ is arbitrary, we have proved $f\in SD(K)$.\qed 
\enddemo 

We need one last rather delicate structural result. 

\proclaim{Lemma 4.13} 
Let $f$ belong to $SD(K)$. 
There exists a non-negative upper semi-continuous function $F$, 
belonging to $SD(K)$, so that $F+|f|$ is upper semi-continuous. 
\endproclaim 

We can now easily prove that $SD(K)$ is a complex function lattice, 
completing the proof of Theorem~5 of the Introduction. 

\proclaim{Theorem 4.14} 
Let $f\in SD(K)$. Then $|f|\in SD(K)$.
\endproclaim 

\demo{Proof} 
Let $f\in SD$. Hence by (a) and (c) of Proposition~4.12, 
$|f|\in B_{1/2}^0(K)$. 
Choosing $F$ as in Lemma~4.13, $F\in B_{1/2}^0(K)$ and hence 
$F+ |f| \in B_{1/2}^0 (K)$ by 4.12. 
But then since $F+|f|$ is upper semi-continuous, 
$F+|f| \in SD$ by 4.12(b), so $|f|\in SD$.\qed 
\enddemo 

\demo{Remark} 
Of course 4.14 is equivalent to the statement that if $f,g$ are real 
functions in $SD$, then $\sqrt{f^2+g^2}$ is in $SD$. 
An indirect proof seems necessary, in view of the fact that only Lipschitz 
functions operate on $SD$ (Proposition~2.8). 
We also  give an alternate proof of 4.14 
in the remarks at the end of this section, 
using the structural characterization of $SD$ given there. 
\enddemo 

\demo{Proof of Lemma 4.13} 
Let $f$ and $\varep$ be as in the statement. 
By the argument in the Remark following the proof of Theorem~3.2, it 
suffices to construct $F$ a non-negative strong $D$ upper semi-continuous 
function so that 
$$F+\Re \mu f\ \text{ is upper semi-continuous for all $\mu$ with 
$|\mu|=1$.} 
\leqno(36)$$ 
Let $\varep>0$, and choose $(f_n)$ simple $D$-functions with 
$$\sum \|f_n\|_D < \infty \quad\text{and}\quad 
f=\sum f_n\ . 
\leqno(37)$$ 
Now for each $n$, let $F_n = \osc_\omega f_n$. 
Since $f_n$ is simple, it is of finite index, and hence $F_n=\osc_{m_n}f_n$ 
for some $m_n<\infty$; thus by Lemma~4.9, $F_n$ is simple. 
Now it follows by Theorem~3.2 that 
$$\|F_n+|f_n|\, \|_\infty \le 2\|f_n\|_D\ . 
\leqno(38)$$ 
Thus we have that 
$$\sum \|F_n\|_\infty <\infty\ . 
\leqno(39)$$ 
It follows that the series $\sum F_n$ converges in $D$-norm, to a function 
$F$ say. 
Indeed, we have for all $n$ 
that $\|F_n\|_D\le 2\|F_n\|_\infty$ since the $F_n$'s are  
non-negative upper semi-continuous. 
Hence $F\in SD$, and $F$ is non-negative upper semi-continuous, since it 
is a uniform limit of such functions. 
Now since $\sum f_n$ converges uniformly to $f$, then given $\mu$ a 
scalar with $|\mu|=1$, $\sum \Re \mu f_n$ converges uniformly to $\Re \mu f$.  
But by the Remark following the proof of Theorem~3.2, we have that 
$F_n+ \Re \mu f_n$ is upper semi-continuous for all $n$. 
Now by (37) and (38), $\sum F_n + \Re \mu f_n$ converges uniformly to 
$F+\Re \mu f$; thus $F+\Re \mu f$ is upper semi-continuous, being a uniform 
limit of such functions. 
Hence (36) holds, completing the proof.\qed 
\enddemo 

\demo{Remark} 
Suppose $f$ is real-valued, and $\varep>0$ is given. 
Then with a little more care in the proof, using the fact  that then 
$\|F_n + |f_n|\,\|_\infty = \|f_n\|_D$ for all $n$, we may choose $F$ 
satisfying the conclusion of the Lemma with $\|F+|f|\,\|_\infty \le 
\|f\|_D +\varep$. 
The proof also constructs our $F$ so that $F\pm f$ are both upper 
semi-continuous. 
Thus setting $\lambda = \|F+|f|\,\|_\infty$, 
$u= {\lambda -F+f\over2}$, $v= {\lambda -F-f\over2}$, 
we have that $u,v\ge0$ are lower semi-continuous $SD$-functions with 
$f= u-v$ and $\|u+v\|_\infty <\|f\|_D +\varep$; that is, we recapture 
Proposition~4.8.  
\enddemo 

We next give several examples of functions in $D\sim SD$, illustrating 
the invariants for $SD$ given above. 

\proclaim{Example 1} 
A bounded upper semi-continuous function which is not strong-$D$.
\endproclaim 

Let $K^{(n)}\ne\emptyset$ for all $n$. 
By the discussion in Remark~1 after Corollary~2.7, we may choose 
$U_1,U_2,\ldots,$ disjoint open sets, and for each $n$, a set $A_n\subset 
U_n$ with $i(A_n\mid U_n) = \|\chix_{A_n}\|_D =n $. 
Let $g= \sum_{n=1}^\infty {\chix_{A_n}\over n}$ and $f=\osc_\omega g$. 
Then by Proposition~3.10, for each $n$, we may choose sets 
$K_n^0 \supset K_n^1 \supset\cdots\supset K_n^{n_1}$ with $K_n^0 = U_n$, 
$K_n^i$ relatively closed nowhere dense in $K_n^{i-1}$ for all 
$1\le i\le n$, $K_n^n\ne\emptyset$, $K_n^{n+1}= \emptyset$, so that 
$$\osc_\omega g= \osc_n g = {j\over n} \text{ on } K_n^j\sim K_n^{j+1}
\ \text{ for all }\ 0\le j\le n\ . 
\leqno(40)$$ 
It follows that $K_n^n \subset \os_n (f,\frac1n)$, hence 
$i_B(f,\frac1n)\ge n$ for all $n$, so $f$ fails (31), and so $f\notin SD$ 
by Proposition~4.12. 

Alternatively, we may argue directly that $g$ itself is not in $SD$ by 
showing that $(\osc_ng)$ fails to converge uniformly; since $g=u-v$ 
for some non-negative upper semi-continuous functions, we have that either 
$u$ or $v$ cannot be strong $D$. 

\proclaim{Example 2} 
A $D$-function $f$ with $\widetilde{\osc}_\omega f\ne \osc_\omega f$, 
and $i_D f=\omega$. 
\endproclaim 

Suppose $K$ is compact with $K^{(\omega)} \ne\emptyset$, and let 
$p\in K^{(\omega)}$. Choose disjoint open subsets $U_1,U_2,\ldots$ of $K$ 
with $p\notin \bigcup_{n=1}^\infty \bar U_n$ so that 
$$\dist (p,U_n) \to 0\ \text{ as }\ n\to\infty 
\leqno(41)$$ 
and such that for all $n$, there is a set $A_n$ with 
$$\bar A_n \subset U_n\text{ and } i(A_n\mid U_n) = n= \|\chix_{A_n}\|_D\ .$$ 
Now set $f= \sum_{n=1}^\infty \chix_{A_n/n}$. 
Then $\|f\|_D =1$, so $\osc_\alpha f\le 1$ for all ordinals $\alpha$. 

Now it follows that $A\dfeq \bigcup_{n=1}^\infty \bar A_n \cup \{p\}$ 
is closed and $f$ is zero off $A$, hence also $\osc f|_{\sim A}\equiv 0$. 
If $x\in A$, $x\ne p$, then $x\in U_n$ for some $n$, and then 
$\osc_n f(x) = \osc_{n+1} f(x)$. 
Hence 
$$i_D f|_{\sim p} \le \omega\ . 
\leqno(42)$$ 
Now we have that 
$$\osc_k f(p) =0\ \text{ for }\ k=1,2,\ldots\ .
\leqno(43)$$ 

Indeed we have that 
$$\lim_{\scriptstyle y\to p\atop \scriptstyle y\ne p} \osc_k f(y) 
= 0\ \text{ for any }\ k\ . 
\leqno(44)$$ 

For, fix $k$, let $\varep>0$, choose $N$ with $\frac{k}N<\varep$, 
then choose $\delta>0$ so that $\rho (y,p) <\delta$ and 
$y\in \bigcup_{j=1}^\infty U_j$ implies $j\ge N$, where $\rho$ is the 
metric on $K$. Then if $0<\rho (y,p) < \delta$, $y\in \bigcup_{j=1}^\infty 
U_j$, $\osc_k f\le \frac{k}N <\varep$, while otherwise $\osc_kf(y)=0$. 

We now easily obtain (43) by induction and the continuity of $f$ at $p$. 
Indeed, (43) holds immediately for $k=1$. 
Suppose proved for $k$. 
But then by (44), 
$$\widetilde{\osc}_{k+1} f(p) \le \lim_{y\to p} |f(y) - f(p)| 
+ \lim_{y\to p} \osc_k f(y)=0\ .$$ 
Again by (44) for ``$k$'' $= k+1$, we obtain that 
$\lim_{y\to p} \widetilde{\osc}_{k+1} f(y) =0$, so 
$\osc_{k+1} f(p) =0$. 

Now (43) immediately yields  that 
$$\widetilde{\osc}_\omega f(p) =0\ . 
\leqno(45)$$ 

However for each $n$, we may choose $x_n\in U_n$ with $\osc_\omega f(x_n) = 
\osc_n f(x_n) \ge1$. 
Since then $x_n\to p$, $\osc_\omega f(p)\ge1$, so $\widetilde{\osc}_\omega 
f \ne \osc_\omega f$. 
But $\osc_\alpha f\le1$ everywhere, for any $\alpha$. 
So 
$$\osc_\omega f(p) = \osc_{\omega+1}(p) =1\ .
\leqno(46)$$ 
Thus by (42) and (46), $\osc_\omega f(x) = \osc_{\omega+1} f(x)$ for all $x$ 
proving $i_D f=\omega$. 

\proclaim{Example 3} 
A $D$-function $g$ with $i_D g=\omega+1$.
\endproclaim 

Let $K$ and $f$ be as in the preceding example and let $g= f+\chix \{p\}$. 
We then have immediately, by (44), that $\lim_{y\to p,\ y\ne p} \osc_k g(y)=0$ 
for all $k$. 

It then follows easily by induction that 
$$\osc_k g(p) = 1 \ \text{ for all } \ k\ . 
\leqno(47)$$ 
Now we obtain that $\widetilde{\osc}_\omega g= \osc_\omega f$, so since 
$\widetilde{\osc}_\omega g$ is upper semi-continuous, 
$\widetilde{\osc}_\omega g= \osc_\omega g$. 
Now if again $x_n \in U_n$ with $\osc_\omega f(x_n) = 
\osc_n f(x_n) =1$, then 
$$\leqalignno{
\widetilde{\osc}_{\omega+1} g(p) & \ge \lim_{n\to\infty} 
|g(x_n) - g(p)| + \osc_\omega g(x_n) 
&(48)\cr 
&= 1+1 =2\ .\cr}$$ 
But we easily have that $\widetilde{\osc}_{\omega+1} g\le 2$, whence also 
$\osc_{\omega+1}g\le2$, so $\osc_{\omega+1} g(p)=2$, showing $i_Dg\ge 
\omega+1$. Finally, since $i_D g|_{\sim p} = \omega$ and 
$\osc_\omega g|_{\sim p} \le1$, we obtain that 
$\widetilde{\osc}_{\omega+2} g(p) \le2$, proving that $i_D g=\omega+1$. 

Of course Examples 2 and 3 both produce functions in $D\sim SD$, by 
Corollary~4.5. 
(It is shown in \cite{R2} that for all $\alpha <\omega_1$, there exists 
a $D$-function $f:[0,1]\to\real$ with $i_D f=\alpha$; an analogous 
result for the positive oscillations was obtained previously 
in \cite{KL}.) 

\proclaim{Example 4} 
A function in $B_{1/2}^0 (K) \sim D(K)$. 
\endproclaim

Our construction is similar to one in \cite{HOR}. 
First fix $n$ and $K_0 \supset \cdots \supset K_{n+1}$ with $K_0 = K$, 
$K_i$ closed nowhere dense in $K_{i-1}$, $1\le i\le n$, $K_n \ne\emptyset$, 
$K_{n+1}=\emptyset$. 
Now let $a_j = (-1)^j/(j+1)$, $0\le j\le n$, then define $f=f_n$ by 
$f= a_j $ on $K_j\sim K_{j+1}$ for all  $0\le j\le n$. 
We then have by Proposition~3.11 that 
$$\|f\|_D = \sum_{j=1}^n \left( {1\over j} + {1\over j+1}\right) + 
{1\over n+1} \sim 2\log n\ . 
\leqno(49)$$ 
However we have 
$$\varep i (f,\varep) \le 3\ \text{ for all }\ \varep>0\ . 
\leqno(50)$$ 
Indeed, if $(\varep_i) =(\frac1i + \frac1{i+1})$, then $\os_j (f,(\varep_i)) 
= K_j$ for all $j\ge 1$. 
Hence if $\frac1j +\frac1{j+1} \le \varep < \frac1{j-1} +\frac1j$, $j>1$, 
then $\os_j (f,\varep) \subset K_j$, but $\osc f| K_j \le 
\frac1{j+1} +\frac1{j+2}$, so $\os_{j+1} (f,\varep)= \emptyset$, and 
$\varep i(f,\varep) < (\frac1{j-1} +\frac1j)\cdot j\le 3$ 
while $\os_1(f,\varep) = \emptyset$ if $\varep >\frac32$. 

Now again assume $K$ is compact with $K^{(\omega)} \ne\emptyset$; let 
$p\in K^{(\omega)}$, and again choose disjoint open subsets $U_1,U_2,\dots$ 
of $K$, with $p\notin \bigcup_{n=1}^\infty \bar U_n$, satisfying (41), 
so that also $U_n^{(n)} \ge n+2$ for all $n$. 

It then follows that for each $n$ we may choose closed sets 
$K_0^n \supset \cdots \supset K_n^n$ with $K_j^n$ nowhere dense in $K_{j-1}^n$ 
for all $0\le j\le n$, $K_n^n\ne \emptyset$, and also $K_0^n$ a closed 
nowhere dense subset of $U_n$. 
Now let $f_n$ be the function on $U_n$ with $f_n = 0$ on 
$U_n\sim K_0^n$, $f_n = (-1)^j /j+1$ on $K_j^n \sim K_{j+1}^n$, 
$0\le j\le n$ (with $K_{n+1}^n =\emptyset$). 
Then it follows by (49) and (50) that 
$$\|f_n\|_{D(U_n)} \sim \log n\ ,\qquad 
\varep i(f_n|U_n,\varep) \le 3\ \text{ for all }\ \varep 
\leqno(51)$$ 
(and again $i(f_n|U_n,\varep) = 0$ if $\varep>\frac32$). 

Now let $g= \sum_{n=2}^\infty (f_n/\sqrt{\log n}) \chix_{U_n}$. 
It follows immediately from (51) that $g\notin D(K)$. 
However if we fix $k$ and let $g_k = \sum_{n=k}^\infty 
(f_n/\sqrt{\log n})\chix_{U_n}$, then $\varep i(g_k,\varep) 
\le 3/\sqrt{\log k}$. 
But then since $\sum_{n=2}^{k-1} (f_n/\sqrt{\log n}) \chix_{U_n}$ in $D(K)$, 
we have that $\olim_{\varep\to0} \varep i (g,\varep) \le 6/\sqrt{\log k}$, 
whence $g\in B_{1/2}^0 (K)$. 

\proclaim{Example 5} 
A function in $B_{1/2}^0 (K)\cap (D(K)\sim SD(K))$. 
\endproclaim 

Let $K$, $(U_n)$, and $(f_n)$ as in the preceding example, and now set 
$f= \sum_{n=2}^\infty (f_n/\log n) \chix_{U_n}$. 
Then $f\in D$ by localization and (49). 
The fact that $f\in B_{1/2}^0 (K)$ follows from the argument for 
Example~4. 
To see  that $f\notin SD$, we need only show (by Proposition~4.4) that 
$(\osc_n f)$ does not converge uniformly. 
Now fixing $k$, then $\osc_k f|U_n \le  (2k/\log n) \to 0$ as $n\to\infty$. 
However by Proposition~3.10, $\|\osc_n f_n|U_n\|_\infty \sim 2\log n$, 
so $\|\osc_n f|U_n\|_\infty \sim 2$. 

We pass now to an intrinsic criterion for distinguishing strong $D$-functions.  
We first need an analogue of the finite oscillation sets, for general 
ordinals. 

\demo{Definition}  
Let $(\alpha_1,\ldots,\alpha_n)$ be non-zero ordinals and $(\varep_1,\ldots,
\varep_n)$ be non-negative numbers. Define the sets 
$\os_j (f,(\alpha_i),(\varep_i))$ inductively as follows: 
$\os_1(f,(\alpha_i),(\varep_i)) = \{x:\osc_{\alpha_1} f(\alpha) \ge 
\varep_1\}$. 
$\os_{j+1} (f,(\alpha_i),(\varep_i)) = \{x:\osc_{\alpha_{j+1}} f|L(x) 
\ge \varep_{j+1}\}$, where $L= \os_j (f,(\alpha_i),(\varep_i))$. 
In case $\alpha_i = \alpha$ and $\varep_i=\varep$ for all $i$, set  
$\os_j(f,\alpha,\varep) = \os_j(f,(\alpha_i),(\varep_i))$. 
Also, for convenience, set $\os_0 (f,\alpha,\varep) =K$. 
\enddemo 

We then have the following result, generalizing Lemma~3.6. 

\proclaim{Lemma 4.16} 
Let $f:K\to\complex$ be given. Then for all $n$, non-zero ordinals 
$\alpha_1,\ldots,\alpha_n$, and $x\in K$, 
$$\osc_{\alpha_1+\cdots+\alpha_n} f(x) = 
\sup\biggl\{\sum_{i=1}^n \varep_i :0 \le \varep_i\ \text{ and }\ 
x\in \os_n (f,(\alpha_i),(\varep_i))\ . 
\leqno(51)$$ 
\endproclaim 

We only need here the fact that $\osc_{\alpha_1+\cdots+\alpha_n}f(x)$ 
dominates the right side of (51). 
The proof of the other estimate may be found in \cite{R2}. 

\proclaim{Sublemma 4.17} 
Let ordinals $\gamma,\beta$, 
$L$ a non-empty subset of $K$, and $\delta>0$ be given. 
If $\osc_\gamma f\ge \delta$ on $L$, then 
$$\osc_{\gamma+\beta} f(x) \ge \delta + \osc_\beta (f|L)(x)
\ \text{ for all }\ x\in L\ . 
\leqno(52)$$ 
\endproclaim 

\demo{Proof} 
By induction on $\beta$. 
This is trivial for $\beta=0$. 
Suppose proved for $\beta$, and let $x\in L$. Then 
$$\eqalign{\widetilde{\osc}_{\gamma+\beta+1} f(x) 
& \ge \olim\limits_{y\to x} |f(y)-f(x)| + \osc_{\gamma+\beta} f(y)\cr 
&\ge \delta +\olim\limits_{\scriptstyle y\to x\atop\scriptstyle y\in L} 
|f(y) - f(x)| + \osc_\beta f|L(y)\cr 
&= \delta + \widetilde{\osc}_{\beta+1} (f|L)(x)\ .\cr}$$ 

Evidently taking upper semi-continuous envelopes now yields 
$\osc_{\gamma+\beta+1} f(x) \ge \delta + \osc_{\beta+1} f|L(x)$. 
The proof for $\beta$ a limit (with (52) holding for all $\beta'<\beta$) is 
immediate.\qed 

We now prove the needed half of 4.16, by showing 
$$\text{if } x\in \osc_n(f,(\alpha_i),(\varep_i))\ ,\ \text{ then }\ 
\osc_{\alpha_1+\cdots+\alpha_n} f(x) \ge \sum_{i=1}^n \varep_i\ . 
\leqno(53)$$ 
We show this by induction on $n$. 
The statement is trivial for $n=1$. 
Suppose proved for $n$, and let $x\in \os_{n+1}(f,(\alpha_i),(\varep_i))$. 
Now setting $L= \os_n(f,(\alpha_i),(\varep_i))$ and 
$\delta = \sum_{i=1}^n \varep_i$, we have that  
$\osc_{\alpha_1+\cdots + \alpha_n} f\ge \delta$ on $L$, by the 
induction hypothesis. 
Hence by Sublemma~4.17, 
$$\eqalignno{ 
\osc_{\alpha_1+\cdots+\alpha_{n+1}} f(x) 
&\ge \delta + \osc_{\alpha_n} f|L(x)\cr 
&\ge \delta + \varep_{n+1} &\qed\cr}$$
\enddemo 

We next define the $\alpha$, $\varep$-index of a function, for $\alpha$ a 
given  ordinal, $\varep>0$. 

\demo{Definition} 
$i(f,\alpha,\varep) = \sup \{ n\ge 0 : \os_n (f,\alpha,\varep) \ne
\emptyset\}$. 
\enddemo 

\proclaim{Corollary 4.18} 
$\varep i (f,\alpha,\varep) \le \|\osc_{\alpha\cdot\omega} f\|_\infty$. 
\endproclaim 

\demo{Proof} 
This follows immediately from (53).  
\enddemo 

We may now formulate the desired criterion. 

\proclaim{Theorem 4.19} 
Let $f:K\to\complex$ be a given bounded function. 
The following are equivalent. 

\iitem{\rm (a)} $f\in SD(K)$. 

\iitem{\rm (b)} {\rm (i)}\quad $\lim_{\varep\to0} \varep i(f,\omega,\varep) 
=0$ and 

\iitem{} {\rm (ii)}\quad $i_D f|W \le \omega$ for all closed $W\subset K$. 
\endproclaim 

In order to prove this, we need the following analogue of 4.10. 

\proclaim{Lemma 4.20} 
For any $\alpha$, and functions $f,g$, 
$$i(f+g,\alpha,\varep) \le i\Bigl(f,\alpha, {\varep\over2}\Bigr) 
+ i\Bigl( g,\alpha,{\varep\over2}\Bigr)\ .$$ 
\endproclaim 

\demo{Remark} 
We only need this for $\alpha=\omega$. 
\enddemo 

We prove 4.20 below, after first using it to give the 

\demo{Proof of Theorem 4.19} 

(a) $\To$ (b). 
$f\in SD$ implies $f\in SD(W)$ for any closed $W\subset K$, so (b)(ii) 
follows immediately from Corollary~4.5(a). 
Now suppose first $f$ is a simple $D$-function and let $n= i_B(f)$. 

As we have seen before, there exist closed non-empty subsets $K= 
K_0\supset \cdots \supset K_n$ with $f|_{K_i\sim K_{i+1}}$ continuous 
for all $i$, $0\le i\le n$ (where $K_{n+1} = \emptyset$). 
It now follows easily that for any $\varep>0$, 
$$i(f,\omega,\varep) \le n\ . 
\leqno(54)$$ 
Indeed, since $f$ is continuous on $K_0\sim K_1$, an open set, 
$\osc_\omega f= 0$ on $K_0 \sim K_1$, so 
$\os_1 (f,\omega,\varep) \subset K_1$. 
Assuming we have shown that $\os_j(f,\omega,\varep) \subset K_j$, 
then again since $\osc_\omega f|K_j = 0$ on $K_j\sim K_{j+1}$, 
$\os_{j+1}(f,\omega,\varep) \subset K_{j+1}$. 
Thus $\os_{n+1}(f,\omega,\varep)=\emptyset$. 

Now assuming $f\in SD$, $f$ real, let $\eta>0$, and choose $\varphi$ a 
simple $D$-function with 
$$\|\varphi -f\|_D <\eta\ .
\leqno(55)$$ 
Let $n=i_B(\varphi)$. 
Then we have that for $\varep>0$, 
$$\eqalign{ \varep i(f,\omega,\varep) 
& \le \varep i\bigl( \varphi -f,\omega,{\varep\over2}\Bigr) 
+ \varep i \Bigl( \varphi,\omega,{\varep\over2}\Bigr)\ 
\text{ (by Lemma 4.20)}\cr 
&\le 2\|\osc_\omega (\varphi-f)\|_\infty + \varep n \ 
\text{ (by Corollary 4.18)}\cr 
&\le 2\eta + \varep n\ \text{ by (55).}\cr}$$ 
Thus $\olim_{\varep\to0} \varep i (f,\omega,\varep) \le 2\eta$, 
proving (b)(i) since $\eta >0$ is arbitrary. 

(b) $\To$ (a). 
Assume without loss of generality that $f$ is real, and  let $\eta >0$. 
Choose $\varep>0$ so that 
$$\varep n<\eta\ ,\ \text{ where }\ n= i(f,\omega,\varep)\ . 
\leqno(56)$$ 
Let then $K_j = \os_{\omega,j} (f,\varep)$ for $j=0,1,2,\ldots$. 
Thus $K_0 \supset K_1\supset \cdots \supset K_n \ne \emptyset$ and 
$K_{n+1}=\emptyset$. 
Now fix $j$; since $\osc f|_{K_j\sim K_{j+1}} \le \osc_\omega 
f|_{K_j\sim K_{j+1}} <\varep$, 
we may choose $\varphi_j \in C_b (K_j\sim K_{j+1})$ with 
$$|\varphi_j (x) - f(x)| \le \varep\text{ for all } x\in K_j\sim K_{j+1} \ . 
\leqno(57)$$ 
Since $K_j\sim K_{j+1}$ is a relatively open subset of $K_j$, it follows 
from (57), (b)(ii) and Theorem~3.2 that 
$$\leqalignno{
\|\varphi_j - f|_{K_j\sim K_{j+1}}\|_D 
& \le \varep + \|\osc_\omega (\varphi_j -f)|_{K_j\sim K_{j+1}} \|_\infty 
&(58)\cr 
& = \varep + \|\osc_\omega f|_{K_j\sim K_{j+1}} \|_\infty \cr 
&\le \varep+\varep\ .\cr}$$
Now setting $\varphi = \sum_{j=0}^n \varphi_j\cdot \chix_{K_j\sim K_{j+1}}$, 
then $\varphi\in  SD$ and 
$$\leqalignno{
\|\varphi -f\|_D & \le \sum_{j=0}^n 
\|(\varphi_j -f) \chix_{K_j\sim K_{j+1}}\|_D
&(59)\cr
&\le 2\sum_{j=0}^n \|(\varphi_j -f)|_{K_j \sim K_{j+1}}\|_D\cr 
&\le 4(n+1)\varep\ \text{ by (58)}\cr 
&\le 4\eta +\varep\ \text{ by (56)}\cr 
&\le 5\eta\ .\cr}$$ 
Since $\eta>0$ is arbitrary, we obtain that $f\in SD$.\qed 
\enddemo 

It remains to prove Lemma 4.20. 
The proof is practically the same as the argument for Theorem~2.8(a) of 
\cite{CMR}, but we give it here for completeness. 

\proclaim{Lemma 4.21} 
Let $\alpha$ be an ordinal, $W_1,\ldots,W_n$ be closed non-empty sets 
with $K= \bigcup_{i=1}^n W_i$, and $f:K\to \complex$ be a bounded function. 
Then 
$$\osc_\alpha f= \max_{1\le i<n} \osc_\alpha f|W_i \chix_{W_i}\ . 
\leqno(60)$$ 
\endproclaim 

\demo{Proof} 
This is easily established by induction. 
Thus,  suppose proved for $\alpha$, let $x\in K$, and choose $(x_n)$ in $K$ 
with $x_n\to x$ and $\widetilde{\osc}_{\alpha+1}f(x) = \lim_{n\to\infty} 
|f(x_n) - f(x)| + \osc_\alpha f(x_n)$. 
After passing to a subsequence, we may assume there is an $i$ with 
$x_n\in W_i$ and $\osc_\alpha f(x_n) = \osc_\alpha f|K_i(x_n)$ for all $n$. 
Then since $W_i$ is closed, $x\in W_i$, and 
$\widetilde{\osc}_{\alpha+1} f(x) \le \widetilde{\osc}_{\alpha+1} 
f|W_i (x) \le \max_i \widetilde{\osc}_{\alpha+1}  f|W_i(x)$. 
(60) now follows for $\alpha+1$, by taking upper semi-continuous envelopes. 
We omit the even simpler proof for limit $\alpha$.\qed 
\enddemo 

\demo{Remark} 
It follows immediately from Theorem 3.2 and Lemma~4.21 that 
{\it if $W_1,\ldots,W_n$ are closed non-empty sets with 
$K= \bigcup_{i=1}^n W_i$ and $f\in D(K)$, $f$ real-valued, then} 
$$\|f\|_D = \max_i \|f|W_i\|_{D(W_i)} \quad ,\quad 
\|f\|_{qD} = \max_i \|f|W_i\|_{qD(W_i)}\ .$$ 
\enddemo 

\demo{Proof of Lemma 4.20} 
Let $f,g$ be as in 4.20 and $\varep>0$ be given. 
For each $n=1,2,\ldots$ and $\btheta = (\theta_1,\ldots,\theta_n)$ with 
$\theta_i=0$ or 1 for all $1\le i\le n$, we define closed subsets  
$L(\btheta)$ of $K$ as follows: 
$$L(0) = \Bigl\{ x\in K:\osc_\alpha f(x) \ge {\varep\over2}\Bigr\}\quad ;
\quad 
L(1) = \Bigl\{ x\in K :\osc_\alpha (x) \ge {\varep\over2}\Bigr\}\ . 
\leqno(61)$$ 
If $n\ge1$ and $L(\btheta) = L(\theta_1,\ldots,\theta_n)$ is defined, let 
$$\left\{ \eqalign{ 
&L(\theta_1,\ldots,\theta_{n+1}) 
= \Bigl\{ x\in L(\btheta) :\osc_\alpha f\mid L(\btheta) \ge{\varep\over2} 
\Bigr\} \ \text{ if }\  \theta_{n+1} =0\cr 
&L(\theta_1,\ldots,\theta_{n+1}) = 
\Bigl\{ x\in L(\btheta) :\osc_\alpha g\mid L(\btheta) \ge {\varep\over2}
\Bigr\} \ \text{ if }\ \theta_{n+1} = 1\ .\cr}\right.
\leqno(62)$$ 

These sets are closed, since $\osc f$, $\osc g$ are upper semi-continuous 
functions. We then have for all $n$ that  
$$\os_n (f+g,\alpha,\varep) \subset \bigcup_{\btheta \in \{0,1\}^n} 
L(\btheta)\ . 
\leqno(63)$$ 

We prove this by induction on $n$. 
Now for $n=1$, since $\osc_\alpha (f+g) \le \osc_\alpha f+\osc_\alpha g$, 
we then have that $\osc_\alpha (f+g)(x)\ge\varep$ implies 
$\osc_\alpha f(x) \ge {\varep\over2}$ or $\osc_\alpha g(x)\ge{\varep\over2}$; 
this gives $\os_1 (f+g,1,\varep) \subset  L(0) \cup L(1)$. 
Suppose (63) is proved for $n$, and suppose $K_n = \os_n (f+g,\alpha,
\varep)$ and $x\in \os_{\alpha,n+1} (f+g,\alpha,\varep)$. 
Thus $\osc_\alpha (f+g)\mid K_n (x)\ge \varep$. 
By Lemma~4.21 and (63), we may then choose $\theta \in \{0,1\}^n$ with 
$x\in K_n\cap L(\theta)$ and  
$$\eqalign{\osc_\alpha (f+g)\mid K_n(x) 
& = \osc_\alpha  (f+g)\mid K_n\cap L(\theta) (x)\cr 
&\le \osc_\alpha (f+g) \mid L(\theta) (x)\cr 
&\le \osc_\alpha f \mid L(\theta)(x) + \osc_\alpha g\mid L(\theta)(x)\ .\cr}$$ 
It follows immediately that $x\in L(\theta_1,\ldots,\theta_n,0)\cup 
L(\theta_1,\ldots,\theta_n,1)$; thus (62) holds at $n+1$. 

Next, fix $n$ and $\theta \in \{0,1\}^n$. Let 
$$j= j(\theta) = \card \{1\le i\le n: \theta_i =0\}\quad ,\quad 
k= k(\theta) = \card \{1\le i\le n:\theta_i=1\}\ . 
\leqno(64)$$ 
Then we claim 
$$L(\theta) \subset \os_j \Bigl( f,\alpha,{\varep\over2}\Bigr) \cap 
\os_k \Bigl( g,\alpha ,{\varep\over2}\bigr)\ . 
\leqno(65)$$ 

Again we prove this by induction on $n$. 
The case $n=1$ is trivial, by the definitions of $L(0)$ and $L(1)$. 
Now suppose (65) is proved for $n$, and $(\theta_1,\ldots,\theta_{n+1})$ 
is given; let $j= j(\theta_1,\ldots,\theta_n)$ and 
$k= k(\theta_1,\ldots,\theta_n)$. 
Now if $\theta_{n+1}=0$, then $j(\theta_1,\ldots,\theta_{n+1}) = j+1$ 
and $k(\theta_1,\ldots,\theta_{n+1}) = k$; then by (65), 
$L(\theta_1,\ldots,\theta_{n+1})\subset L(\theta_1,\ldots,\theta_n)\subset 
\os_k (f,\alpha,\frac{\varep}2)$ and by definition and (65), 
$$\eqalign{ L(\theta_1,\ldots,\theta_{n+1}) 
& \subset \left\{ x\in \os_j \Bigl(f,\alpha,{\varep\over2}\Bigr)  : 
\osc_\alpha f \mid \os_j \bigl( f,\alpha,{\varep\over2}\Bigr) (x) \ge
{\varep\over2}\right\}\cr 
& = \os_{j+1} \Bigl( f,\alpha,{\varep\over2}\Bigr)\ .\cr}$$ 
Of course if $\theta_{n+1}=1$, we obtain by the same reasoning that 
$L(\theta_1,\ldots,\theta_{n+1})\subset \os_j (f,\alpha,{\varep\over2}) 
\cap \os_{k+1} (g,\alpha,{\varep\over2})$ and $j=j(\theta_1,\ldots, 
\theta_{n+1})$, $k+1=k(\theta_1,\ldots,\theta_{n+1})$; thus (65) is proved 
for $n+1$, and so established for all $n$ by induction.

Now suppose, for a given $n$, that $\os_n(f+g,\alpha,\varep) \ne\emptyset$. 
Then by (63), there is a $\theta\in \{0,1\}^n$ with $L(\theta)\ne\emptyset$. 
Thus letting $j$ and $k$ be as in (64), we have by (65) that 
$\os_j (f,\alpha,{\varep\over2})\ne\emptyset$ and 
$\os_k (g,\alpha,{\varep\over2})\ne \emptyset$. 
But then $n=j+k\le i(f,\alpha,{\varep\over2}) + i(g,\alpha,{\varep\over2})$. 
Lemma~4.20 is thus established.\qed 
\enddemo 

\demo{Remark} 
The proof of Theorem 4.19 yields a generalization for functions of 
arbitrary $D$-index, showing that $SD$ occupies a special place in $D$. 
We define, for $\varep>0$ and $f:K\to\real$, a sequence of sets 
$\os_1^\infty (f,\varep), \os_2^\infty (f,\varep),\ldots$ by 
induction as follows: 
$\os_1^\infty (f,\varep) = \{x:\|f\|_{qD(x)} \ge\varep\}$, and 
$\os_{n+1}^\infty (f,\varep) = \{x\in L:\|f\mid L\|_{qD(x)} \ge\varep\}$ 
where $L= \os_n^\infty (f,\varep)$. 
($\|f\|_{qD(x)}$ is defined preceding Corollary~3.5.) 
Now let 
$i_{qD} (f,\varep) = \sup \{n:\os_n^\infty (f,\varep) \ne\emptyset\}$. 
It follows from Theorem~3.2 that in fact there are ordinals 
$\alpha_1,\alpha_2,\ldots$ so that $\os_n^\infty (f,\varep) = \os_n 
(f,(\alpha_i),\varep)$ for all $n$; again by Corollary~4.18 and Theorem~3.2 
we obtain that
$$\varep i_{qD} (f,\varep) \le \|\osc_{\sum \alpha_i} f\|_\infty \le 
\|f\|_{qD}\ .$$ 
The proof of Theorem~4.19 now yields the following rather surprising result: 
\enddemo 

\proclaim{Theorem} 
Let $f$ be a bounded function on $K$. 
Then $f$ is a strong $D$-function if and only if $\lim_{\varep\to0} 
\varep i_{qD} (f,\varep)=0$. 
\endproclaim 

This theorem yields that $f\in SD(K)$ implies $|f|\in SD(K)$ 
(Theorem~4.14 above). 
Indeed, it follows easily from the proof of Proposition~1.16, that 
if $f:K\to \complex$ is a given bounded function, then 
$\|\, |f|\,\|_{qD(x)} \le \|f\|_{qD(x)}$ for all $x\in K$. 
But then $i_{qD} (|f|,\varep) \le i_{qD} (f,\varep)$ for all 
$\varep >0$, so $f\in SD$ implies $\lim_{\varep\to0} \varep i_{qD} 
(|f|,\varep) = 0$, thus $|f|\in SD$.

\enddocument